\documentclass[10pt]{article}

\usepackage[dvips]{epsfig} 
\usepackage{amssymb,amsmath,amsthm,mathrsfs} 
\usepackage{graphicx}
\usepackage{lineno}

\usepackage[authoryear,sort&compress]{natbib}
\usepackage{url}
\usepackage{enumerate}
\usepackage{colordvi,color,pspicture}
\usepackage{psfrag}

\newcommand{\E}[1]{\mathbf{E}\,#1}

\newcommand{\Var}[1]{\mathbf{Var}\,#1}

\newcommand{\C}{\mathcal{C}}

\newcommand{\U}{\mathcal{U}}

\newcommand{\g}{\gamma}
\newcommand{\G}{\Gamma}
\newcommand{\V}{\mathcal{V}}
\newcommand{\A}{\mathcal{A}}
\newcommand{\Y}{\mathcal{Y}}
\newcommand{\y}{\mathsf{y}}
\newcommand{\X}{\mathcal{X}}
\newcommand{\M}{\mathcal{M}}

\newcommand{\mS}{\mathcal{S}}

\newcommand{\NY}{N_{\Y}}
\newcommand{\NS}{N_S}

\newcommand{\NPE}{N_{PE}}
\newcommand{\TY}{T(\Y_3)}

\newcommand{\R}{\mathbb{R}}

\newcommand{\RS}{\mathscr{R}_S}
\newcommand{\UT}{\U(T(\Y_3))}
\newcommand{\Tr}{\mathscr T_r}
\newcommand{\msPr}{\mathscr P_r}
\newcommand{\ve}{\varepsilon}
\newcommand{\ah}{\widehat{\alpha}}
\newcommand{\bh}{\widehat{\beta}}
\newcommand{\bgam}{\text{\pmb{\Large$\gamma$}}}

\DeclareMathOperator{\argmin}{argmin}

\DeclareMathOperator{\support}{support}
 
\DeclareMathOperator{\BIN}{BIN} \DeclareMathOperator{\BER}{BER}

\newtheorem{theorem}{Theorem}[section]

\theoremstyle{definition}

\theoremstyle{remark}
\newtheorem{conjecture}[theorem]{Conjecture}
\newtheorem{remark}[theorem]{Remark}

\begin{document}


\title{Technical Report \# KU-EC-09-6:\\
Spatial Clustering Tests Based on Domination Number of a New Random Digraph Family}
\author{
Elvan Ceyhan
\thanks{Address:
Department of Mathematics, Ko\c{c} University, 34450 Sar{\i}yer, Istanbul, Turkey.
e-mail: elceyhan@ku.edu.tr, tel:+90 (212) 338-1845, fax: +90 (212) 338-1559.
}
}

\date{\today}
\maketitle


\begin{abstract}
\noindent
We use the domination number of a parametrized random digraph family
called proportional-edge \emph{proximity catch digraphs}
(PCDs) for testing multivariate spatial point patterns.
This digraph family 
is based on relative positions of data points from various classes.
We extend the results on the distribution of the domination number of proportional-edge PCDs,
and use the domination number as a statistic for testing segregation and association
against complete spatial randomness.
We demonstrate that the domination number of the PCD
has binomial distribution when size of one class is fixed
while the size of the other
(whose points constitute the vertices of the digraph) tends to infinity and
asymptotic normality when sizes of both classes tend to infinity.
We evaluate the finite sample performance of the test by Monte Carlo simulations,
prove the consistency of the test under the alternatives,
and suggest corrections for the support restriction on the class of points of interest
and for small samples.
We find the optimal parameters for testing each of the segregation and association alternatives.
Furthermore, the methodology discussed in this article is valid for data in
higher dimensions also.
\end{abstract}



\newpage


\section{Introduction}
\label{sec:introduction}
In statistical literature, the problem of clustering 
received considerable attention.
The spatial interaction between two or more classes has important implications especially for plant species.
See, e.g., \cite{pielou:1961}, \cite{dixon:1994,dixon:NNCTEco2002},
\cite{stoyan:2000}, and \cite{perry:2006}.
Recently, a new clustering test based on the relative allocation of points
from two or more classes has been developed.
The method is based on a graph-theoretic approach and is used to test the spatial pattern of complete
spatial randomness (CSR) against segregation or association.
Rather than the pattern of points from one-class with respect to the ground,
the patterns of points from one class with respect to points from other classes are investigated.
\emph{CSR} is roughly defined as the lack of
spatial interaction between the points in a given study area.
\emph{Segregation} is the pattern in which points of one class tend to
cluster together, i.e., form one-class clumps.
On the other hand, \emph{association} is the pattern in which
the points of one class tend to occur more frequently around points from the other class.
For convenience and generality,
we call the different types of points as ``classes",
but the class can be replaced by any characteristic of an observation at a particular location.
For example, the pattern of spatial segregation has been investigated
for species (\cite{diggle:2003}), age classes of plants (\cite{hamill:1986})
and sexes of dioecious plants (\cite{nanami:1999}).

Many methods to analyze spatial clustering have been proposed in the literature (\cite{kulldorff:2006}).
These include Ripley's $K$ or $L$-functions (\cite{ripley:2004}),
comparison of NN distances (\cite{diggle:2003}, \cite{cuzick:1990}),
and analysis of nearest neighbor contingency tables (NNCTs)
which are constructed using the NN frequencies of classes
(\cite{pielou:1961} and \cite{dixon:1994, dixon:NNCTEco2002}).
The tests (i.e., inference) based on Ripley's $K$ or $L$-functions are only appropriate
when the null pattern can be assumed to be the CSR independence pattern,
but not if the null pattern is the RL of points from an inhomogeneous Poisson pattern (\cite{kulldorff:2006}).
But, there are also variants of $K(t)$ that explicitly correct for inhomogeneity
(see \cite{baddeley:2000b}).
Cuzick and Edward's $k$-NN tests are designed for testing bivariate spatial interaction
and mostly used for spatial clustering of cases or controls in epidemiology.
Diggle's $D$-function is a modified version of Ripley's $K$-function (\cite{diggle:2003})
and is appropriate for the case in which the null pattern is the RL
of points where the points are a realization from any arbitrary point pattern.
Ripley's and Diggle's functions are designed to analyze
univariate or bivariate spatial interaction at various scales (i.e., inter-point distances).

In recent years,
the use of mathematical graphs has also gained popularity in spatial analysis (\cite{roberts:2000})
providing a way to move beyond Euclidean metrics for spatial analysis.
Although only recently introduced to landscape ecology, graph theory is well suited to
ecological applications concerned with connectivity or movement (\cite{minor:2007}).
Conventional graphs do not explicitly maintain geographic reference,
reducing utility of other geo-spatial information.
\cite{fall:2007} introduce spatial graphs that
integrate a geometric reference system that ties patches and paths
to specific spatial locations and spatial dimensions
thereby preserving the relevant spatial information.
However, after a graph is constructed using spatial data,
usually the scale is lost (see for instance, \cite{su:2007}).
Many concepts in spatial ecology depend on the idea of spatial adjacency which
requires information on the close vicinity of an object.
Graph theory conveniently can be used to express
and communicate adjacency information allowing one to compute
meaningful quantities related to spatial point pattern.
Adding vertex and edge properties to graphs extends the problem domain to network modeling (\cite{keitt:2007}).
\cite{wu:2008} propose a new measure based on graph theory and spatial interaction,
which reflects intra-patch and inter-patch relationships by
quantifying contiguity within patches and potential contiguity among patches.
\cite{friedman:1983} also propose a graph-theoretic method
to measure multivariate association,
but their method is not designed to analyze spatial interaction
between two or more classes;
instead it is an extension of generalized correlation coefficient
(such as Spearman's $\rho$ or Kendall's $\tau$)
to measure multivariate (possibly nonlinear) correlation.

The graph-theoretic method we use to test spatial randomness is based on \emph{proximity catch digraphs} (PCDs) which are
a special type of \emph{proximity graphs} introduced by \cite{toussaint:1980}.
A \emph{digraph} is a directed graph with
vertices $V$ and arcs (directed edges) each of which is from one
vertex to another based on a binary relation.
Then the pair $(p,q) \in V \times V$ is an ordered pair
which stands for an arc from vertex $p$ to vertex $q$ in $V$.
For example, \emph{nearest neighbor (di)graph} which is defined by placing
an arc between each vertex and its nearest neighbor is a proximity digraph (\cite{paterson:1992}).
The nearest neighbor digraph has the vertex set $V$
and $(p,q)$ as an arc iff $q$ is a nearest neighbor of $p$.
The domination number of PCDs is first investigated for data in one
Delaunay triangle (in $\R^2$) and the analysis is generalized to
data in multiple Delaunay triangles.
Some trivial proofs are omitted and
shorter proofs are given in the main body of the article.
\emph{Data-random digraphs} are directed graphs in which
each vertex corresponds to a data point,
and arcs are defined in terms of some bivariate function on the data.
\cite{priebe:2001} introduced a data random digraph
called \emph{class cover catch digraph} (CCCD) in $\mathbb{R}$ and
extended it to multiple dimensions.
In this model, the vertices correspond to data points from a
single class $\X$ and the definition of the arcs utilizes the
other class $\Y$.
For each $x_i \in \X$ a radius is defined as $r_i=\min_{\y \in \Y} d(x_i,\y)$.
There is an arc from $x_i$ to $x_j$ if $d(x_i,x_j)<r_i$;
that is, the (open) sphere of radius $r_i$ ``catches'' $x_j$.
\cite{devinney:2002a}, \cite{marchette:2003}, \cite{priebe:2003b}, and \cite{priebe:2003a}
demonstrated relatively good performance of CCCDs in classification.
Their methods involve \emph{data reduction} (\emph{condensing}) by
using approximate minimum dominating sets as \emph{prototype sets}
(since finding the exact minimum dominating set is an NP-hard
problem in general --- e.g., for CCCD in multiple dimensions --- (see \cite{devinney:2006}).
For the domination number of CCCDs for one-dimensional data, a SLLN result is proved in (\cite{devinney:2002b}),
and this result is extended by \cite{wiermanSLLN:2008}; furthermore,
a CLT is also proved by \cite{xiangCLT:2009}.
The asymptotic distribution of the domination number of CCCDs for non-uniform data
in $\mathbb{R}$ is also calculated in a rather general setting (\cite{ceyhan:dom-num-CCCD-NonUnif}).
Although intuitively appealing and easy to extend to higher dimensions,
finding the minimum dominating set of CCCD is an NP-hard problem and
the distribution of the domination number of CCCDs is not
analytically tractable for $d>1$.
This drawback has motivated us to define new types of proximity maps.
As alternatives to CCCD,
\cite{ceyhan:CS-JSM-2003} introduced an (unparametrized) type of PCDs called
\emph{central similarity PCDs};
\cite{ceyhan:dom-num-NPE-SPL} also introduced another parametrized family of PCDs
called \emph{proportional-edge PCDs}
and used the domination number of this PCD with a fixed parameter for testing spatial patterns.
The domination number approach is appropriate when at least one of the classes is sufficiently large.
The relative (arc) density of these PCDs are also used for testing the spatial patterns in
(\cite{ceyhan:arc-density-PE}) and (\cite{ceyhan:arc-density-CS}).
These new PCDs are designed to have better distributional and mathematical properties.
These new families are both applicable to pattern classification also.
\cite{ceyhan:CS-JSM-2003} introduced the central similarity
proximity maps and the associated PCDs, and \cite{ceyhan:arc-density-CS}
computed the asymptotic distribution of the relative (arc) density of the
parametrized version of the central similarity PCDs and applied the
method to testing spatial patterns.
\cite{ceyhan:dom-num-NPE-SPL} introduced proportional-edge PCD with expansion parameter $r$,
where the distribution of the domination number of proportional-edge PCD with $r=3/2$ is used in
testing spatial patterns of segregation or association.
\cite{ceyhan:arc-density-PE} computed the asymptotic distribution of the
relative density of the proportional-edge PCD and used it for the same purpose.
\cite{ceyhan:dom-num-NPE-MASA} derived the asymptotic distribution
of the domination number of proportional-edge PCDs for uniform  data.
An extensive treatment of the
PCDs based on Delaunay tessellations is available in \cite{ceyhan:Phd-thesis}.

In this article,
we investigate the use of the domination number of proportional-edge PCDs,
whose asymptotic distribution was computed in (\cite{ceyhan:dom-num-NPE-MASA})
for testing spatial patterns of segregation and association.
Furthermore, we extend this result for the whole
range of the expansion parameter in a more general setting.
By construction, in our PCDs, the further an $\X$ point is from $\Y$ points,
it will be more likely to have more arcs to other $\X$ points,
hence the domination number will be more likely to be smaller.
This probabilistic behavior lends the domination number as
a statistic for testing spatial segregation or association.
In addition to the mathematical tractability and applicability to
testing spatial patterns and classification,
this new family of PCDs
is more flexible as it allows choosing an optimal parameter for
testing against various types of spatial point patterns.

We define proximity maps and the associated PCDs in Section \ref{sec:prox-maps-and-PCDs},
present the asymptotic distribution of the domination number for uniform data in one triangle
and in multiple triangles in Section \ref{sec:asy-gam-NYr},
describe the alternative patterns of segregation and association in Section \ref{sec:alternatives},
present the Monte Carlo simulation analysis to assess the empirical size
and power performance in Section \ref{sec:monte-carlo-sim},
suggest an adjustment for data points from the class of interest
which are outside the convex hull of data from the other class in Section \ref{sec:conv-hull-correction},
suggest a correction method for small sample sizes of the class of interest in Section \ref{sec:small-sample-correction},
provide an example data set in Section \ref{sec:example},
and describe the extension of proportional-edge PCDs to higher dimensions in Section \ref{sec:NYr-higher-D}.
We also provide the guidelines in using this test in Section \ref{sec:disc-conc}.

\section{Proximity Maps and the Associated PCDs}
\label{sec:prox-maps-and-PCDs}
Our PCDs are based on the proximity maps which are defined in a
fairly general setting. Let $(\Omega,\M)$ be a measurable space.
The \emph{proximity map} $N(\cdot)$ is defined as
$N:\Omega \rightarrow 2^\Omega$, where $2^\Omega$ is the power set of $\Omega$.
The \emph{proximity region} associated with $x \in \Omega$, denoted
$N(x)$, is the image of $ x \in \Omega$ under $N(\cdot)$.
The points in $N(x)$ are thought of as being ``closer" to $x \in \Omega$ than
are the points in $\Omega \setminus N(x)$.
Hence the term ``proximity" in the name \emph{proximity catch digraph}.
The $\G_1$-region $\G_1(\cdot)=\G_1(\cdot,\Y):\Omega \rightarrow 2^\Omega$
associates the region $\G_1(x):=\{z \in \Omega: x \in \NY(z)\}$ with each point $x \in \Omega$.
Proximity maps are the building blocks of the \emph{proximity graphs} of \cite{toussaint:1980};
an extensive survey on proximity maps and graphs is available in (\cite{jaromczyk:1992}).

The \emph{proximity catch digraph} $D$ has the vertex set
$\V=\bigl\{ p_1,\ldots,p_n \bigr\}$; and the arc set $\A$ is defined
by $(p_i,p_j) \in \A$ iff $p_j \in N(p_i)$ for $i\not=j$.
Notice that the proximity catch digraph $D$ depends on the \emph{proximity}
map $N(\cdot)$ and if $p_j \in N(p_i)$, then we call the region $N(p_i)$ (and the point $p_i$) \emph{catches} point $p_j$.
Hence the term ``catch" in the name \emph{proximity catch digraph}.
If arcs of the form $(p_i,p_i)$ (i.e., loops) were allowed,
$D$ would have been called a \emph{pseudodigraph} according to some authors
(see, e.g., \cite{chartrand:1996}).

In a digraph $D=(\V,\A)$, a vertex $v \in \V$ \emph{dominates}
itself and all vertices of the form $\{u: (v,u) \in \A\}$.
A \emph{dominating set} $S_D$ for the digraph $D$ is a subset of $\V$
such that each vertex $v \in \V$ is dominated by a vertex in $S_D$.
A \emph{minimum dominating set} $S^*_{D}$ is a dominating set of
minimum cardinality and the \emph{domination number} $\g(D)$ is
defined as $\g(D):=|S^*_{D}|$ (see, e.g., \cite{lee:1998}) where
$|\cdot|$ denotes the set cardinality functional.
See \cite{chartrand:1996} and \cite{west:2001} for more on graphs and digraphs.
If a minimum dominating set is of size one, we call it a \emph{dominating point}.
Note that for $|\V|=n>0$, $1 \le \g(D) \le n$,
since $\V$ itself is always a dominating set.

\begin{figure}[ht]
\centering
\rotatebox{-90}{ \resizebox{3. in}{!}{ \includegraphics{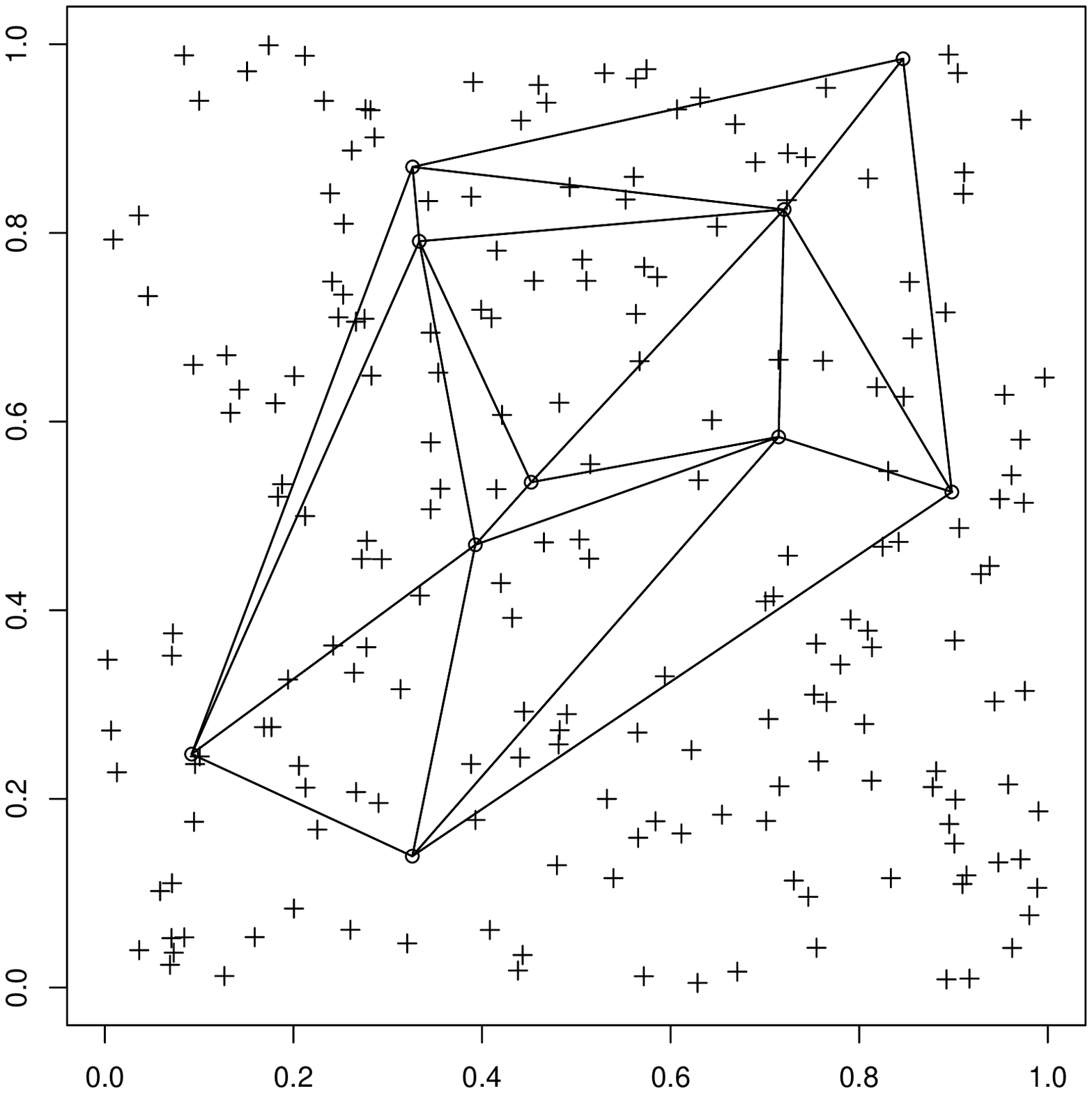}}}
\rotatebox{-90}{ \resizebox{3. in}{!}{\includegraphics{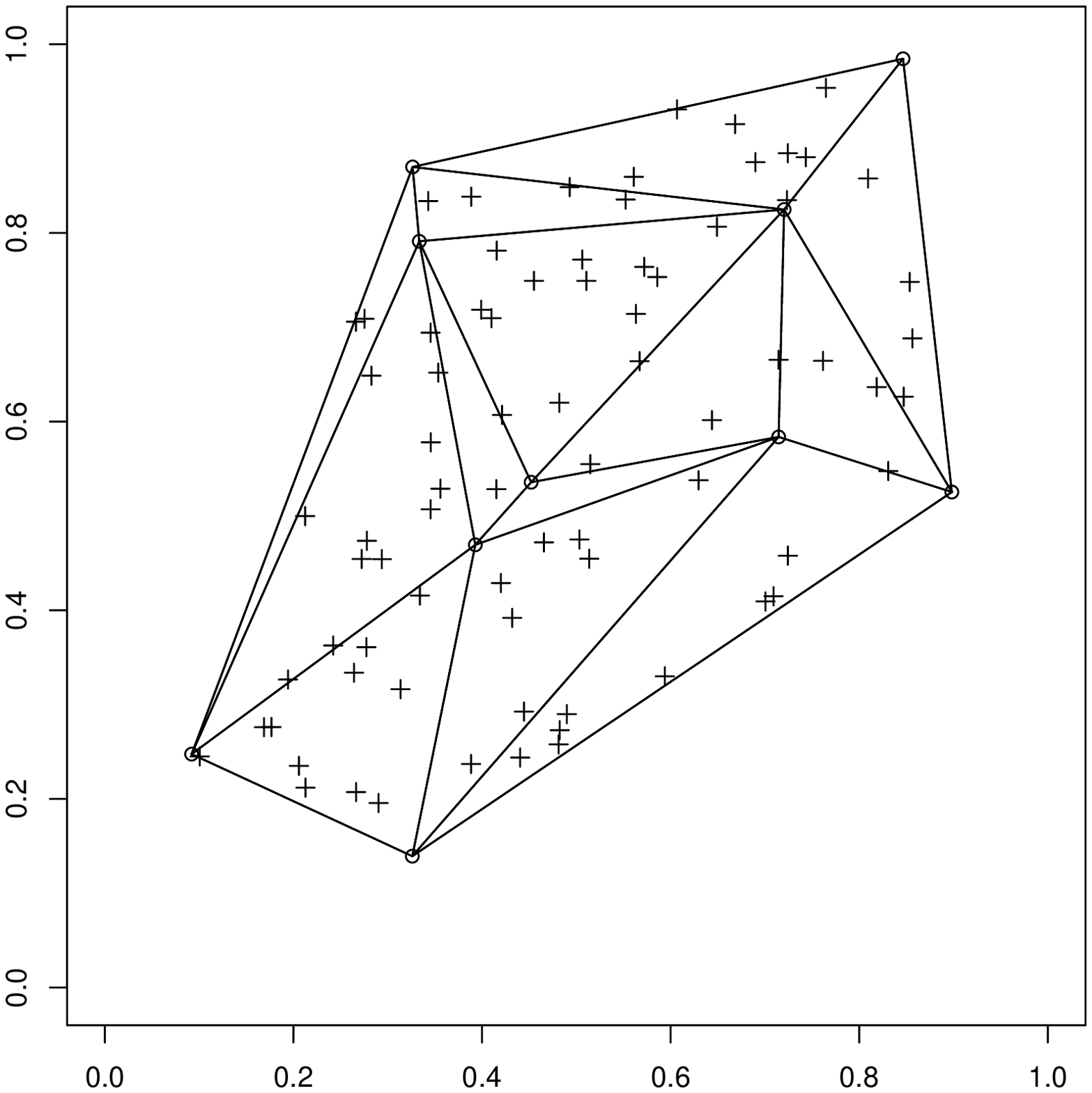}}}
\caption{ \label{fig:deltri}
In left, plotted is a realization of $200$ $\X$ points (pluses, $+$) and
the Delaunay triangulation based on 10 $\Y$ points (circles, $\circ$).
In right, plotted is the 77 $\X$ points which are in the convex hull of $\Y$ points.
Both $\X_n$ and $\Y_m$ are random samples from $\U \bigl((0,1)\times (0,1)\bigr)$,
the uniform distribution on the unit square.}
\end{figure}

We construct the proximity regions using two data sets $\X_n$ and $\Y_m$
of sizes $n$ and $m$ from classes $\X$ and $\Y$, respectively.
Given $\Y_m \subseteq \Omega$, the {\em proximity map} $\NY(\cdot):
\Omega \rightarrow 2^{\Omega}$ associates a {\em proximity region}
$\NY(x) \subseteq \Omega$ with each point $x \in \Omega$.
The region $\NY(x)$ is defined in terms of the distance between $x$ and $\Y_m$.
More specifically, our proportional-edge proximity maps will be based on
the relative position of points from $\X_n$ with respect to the
Delaunay tessellation of $\Y_m$.
In this article, a triangle refers
to the closed region bounded by its edges.
See Figure \ref{fig:deltri} for an example with $n=200$ $\X$ points iid
$\U \bigl((0,1)\times (0,1)\bigr)$, the uniform distribution on the unit
square and the Delaunay triangulation (which yields 13 triangles) is based on $m=10$ $\Y$
points which are also iid $\U \bigl((0,1)\times (0,1)\bigr)$
and 77 of these $\X$ points are inside the convex hull of $\Y$ points.

If $\X_n=\bigl\{ X_1,\ldots,X_n \bigr\}$ is a set of $\Omega$-valued
random variables then $\NY(X_i)$ and $\G_1(X_i)$ are random sets.
If $X_i$ are iid then so are the random sets $\NY(X_i)$.
The same holds for $\G_1(X_i)$.
We define the \emph{data-random proximity catch digraph} $D$
--- associated with $\NY(\cdot)$ --- with vertex set
$\X_n=\{X_1,\cdots,X_n\}$ and arc set $\A$ by
 $$(X_i,X_j) \in \A \iff X_j \in \NY(X_i).$$
Since this relationship is not symmetric, a digraph is used rather than a graph.
The random digraph $D$ depends on the (joint)
distribution of $X_i$ and on the map $\NY(\cdot)$.
For $\X_n=\bigl\{ X_1,\cdots,X_n \bigr\}$, a set of iid random variables
from $F$, the domination number of the associated data-random
PCD based on the proximity map $N(\cdot)$,
denoted $\g(\X_n,N)$, is the minimum number of point(s) that dominate
all points in $\X_n$.
The random variable $\g(\X_n,N)$ depends explicitly on $\X_n$ and $N(\cdot)$
and implicitly on $F$.
Furthermore, in general, the distribution, hence the expectation $\E[\g(\X_n,N)]$, depends on
$n$, $F$, and $N$; $ 1 \leq \E[\g(\X_n,N)] \le n.$
In general, the variance of $\g(\X_n,N)$ satisfies,
$1 \le \Var[\g(\X_n,N)] \le n^2/4$.
For example, the CCCD of \cite{priebe:2001} can be viewed as an
example of PCDs and is briefly discussed in the next section.
We use some of the properties of CCCD in $\R$ as guidelines in
defining PCDs in higher dimensions.

\subsection{Spherical Proximity Maps}
\label{sec:spherical-PCD}
\cite{priebe:2001} introduced the \emph{class cover catch digraphs}
(CCCDs) and gave the exact and the asymptotic distribution
of the domination number of the CCCD based on two sets,
$\X_n$ and $\Y_m$, which are of sizes $n$ and $m$,
from classes, $\X$ and $\Y$, respectively,
and are sets of iid random variables from uniform
distribution on a compact interval in $\R$.

Let $\Y_m=\left \{\y_1,\ldots,\y_m \right\} \subset \R$.
Then the proximity map associated with CCCD is defined as the open
ball $\NS(x):=B(x,r(x))$ for all $x \in \R$, where $r(x):=\min_{\y \in \Y_m}d(x,\y)$ with $d(x,y)$ being the
Euclidean distance between $x$ and $y$ (\cite{priebe:2001}).
That is, there is an arc from $X_i$ to $X_j$ iff there exists an open ball centered at $X_i$
which is ``pure" (or contains no elements) of $\Y_m$ in its
interior, and simultaneously contains (or ``catches") point $X_j$.
We consider the closed ball, $\overline{B}(x,r(x))$ for $\NS(x)$ in this article.
Then for $x \in \Y_m$, we have $\NS(x)=\{x\}$.
Notice that a ball is a sphere in higher dimensions, hence the notation $N_S$.
Furthermore, dependence on $\Y_m$ is through $r(x)$.

A natural extension of the proximity region $\NS(x)$ to $\R^d$ with
$d>1$ is obtained as $N_S(x):=B(x,r(x))$ where
$r(x):=\min_{\y \in \Y_m} d(x,\y)$ which
is called the \emph{spherical proximity map}.
The spherical proximity map $N_S(x)$ is well-defined for all $x \in \R^d$
provided that $\Y_m \not= \emptyset$.
Extensions to $\R^2$ and higher dimensions with the spherical proximity map
--- with applications in classification --- are investigated by
\cite{devinney:2002a}, \cite{marchette:2003}, \cite{priebe:2003b, priebe:2003a},
and \cite{devinney:2006}.

\subsection{The Proportional-Edge Proximity Maps}
\label{sec:PE-PCD}
First, we describe the construction of the
$r$-factor proximity maps and regions, then state some of its basic
properties and introduce some auxiliary tools.
Note that in $\R$ the CCCDs are based on the intervals
whose end points are from class $\Y$.
$I_i=\left(\y_{(i-1):m},\y_{i:m} \right)$ for $i=0,\ldots,(m+1)$ with
$\y_{0:m}=-\infty$ and $\y_{(m+1):m}=\infty$,
where $\y_{i:m}$ is the $i^{th}$ order statistic in $\Y_m$.
This interval partitioning can
be viewed as the \emph{Delaunay tessellation} of $\R$ based on $\Y_m$. So
in higher dimensions, we use the Delaunay triangulation based
on $\Y_m$ to partition the support.

Let $\Y_m=\left \{\y_1,\ldots,\y_m \right\}$ be $m$ points in
general position in $\R^d$ and $T_i$ be the $i^{th}$ Delaunay cell
for $i=1,\ldots,J_m$, where $J_m$ is the number of Delaunay cells.
Let $\X_n$ be a set of iid random variables from distribution $F$ in
$\R^d$ with support $\mS(F) \subseteq \C_H(\Y_m)$
where $\C_H(\Y_m)$ stands for the convex hull of $\Y_m$.
In particular, for illustrative purposes, we focus on $\R^2$ where
a Delaunay tessellation is a \emph{triangulation}, provided that no more
than three points in $\Y_m$ are cocircular (i.e., lie on the same circle).
Furthermore, for simplicity,
let $\Y_3=\{\y_1,\y_2,\y_3\}$ be three non-collinear points
in $\R^2$ and $\TY=T(\y_1,\y_2,\y_3)$ be the triangle
with vertices $\Y_3$.
Let $\X_n$ be a set of iid random variables from $F$ with
support $\mS(F) \subseteq \TY$.
If $F=\U(\TY)$, a composition of translation,
rotation, reflections, and scaling
will take any given triangle $\TY$
to the basic triangle $T_b=T((0,0),(1,0),(c_1,c_2))$
with $0 < c_1 \le 1/2$, $c_2>0$,
and $(1-c_1)^2+c_2^2 \le 1$, preserving uniformity.
That is, if $X \sim \U(\TY)$ is transformed in the same manner to,
say $X'$, then we have $X' \sim \U(T_b)$.
In fact this will hold for any distribution $F$
up to scale.

\begin{figure} [ht]
\centering
\scalebox{.35}{\input{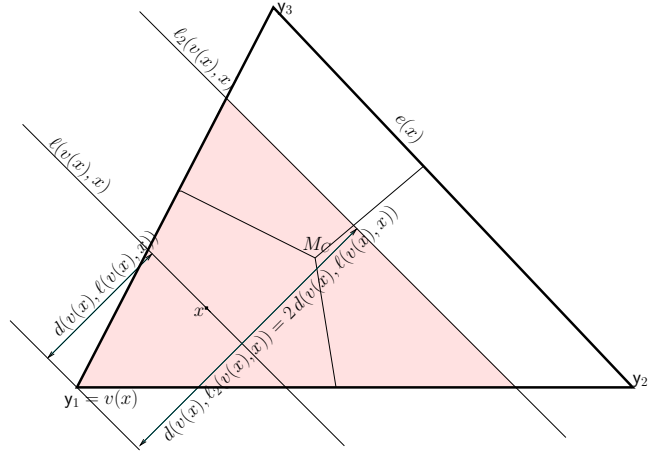}}
\caption{
\label{fig:prox-map-def}
Construction of proportional-edge proximity region, $\NPE^{r=2}(x,M_C)$ (shaded region)
for an $x$ in the CM-vertex region for $\y_1$, $R_{M_C}(\y_1)$.}
\end{figure}

\begin{figure}
\begin{center}
\psfrag{A}{\scriptsize{$\y_1$}}
\psfrag{B}{\scriptsize{$\y_2$}}
\psfrag{C}{\scriptsize{$\y_3$}}
\psfrag{CM}{\scriptsize{$M_C$}}
\psfrag{IC}{\scriptsize{$M_{I}$}}
\psfrag{x}{}
\psfrag{P1}{\scriptsize{$M_1$}}
\psfrag{P2}{\scriptsize{$M_2$}}
\psfrag{P3}{\scriptsize{$M_3$}}
\psfrag{R(A)}{\scriptsize{$R_{M_C}(\y_1)$} }
\psfrag{R(B)}{\scriptsize{$R_{M_C}(\y_2)$} }
\psfrag{R(C)}{\scriptsize{$R_{M_C}(\y_3)$} }
 \epsfig{figure=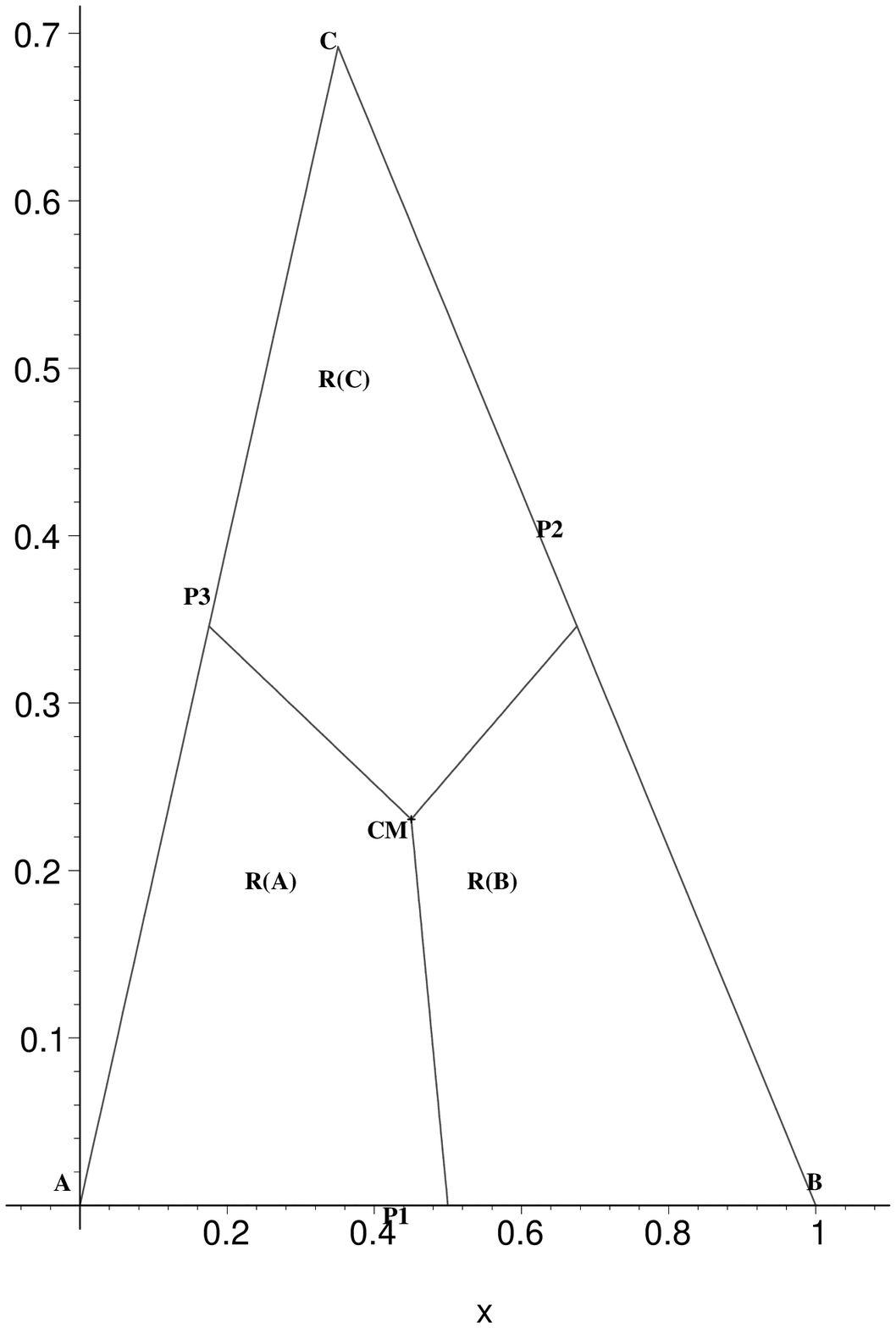, height=140pt,width=200pt}
\psfrag{1}{\scriptsize{$\y_2$}}
\psfrag{R(A)}{\scriptsize{$R_{M_I}(\y_1)$} }
\psfrag{R(B)}{\scriptsize{$R_{M_I}(\y_2)$} }
\psfrag{R(C)}{\scriptsize{$R_{M_I}(\y_3)$} }
 \epsfig{figure=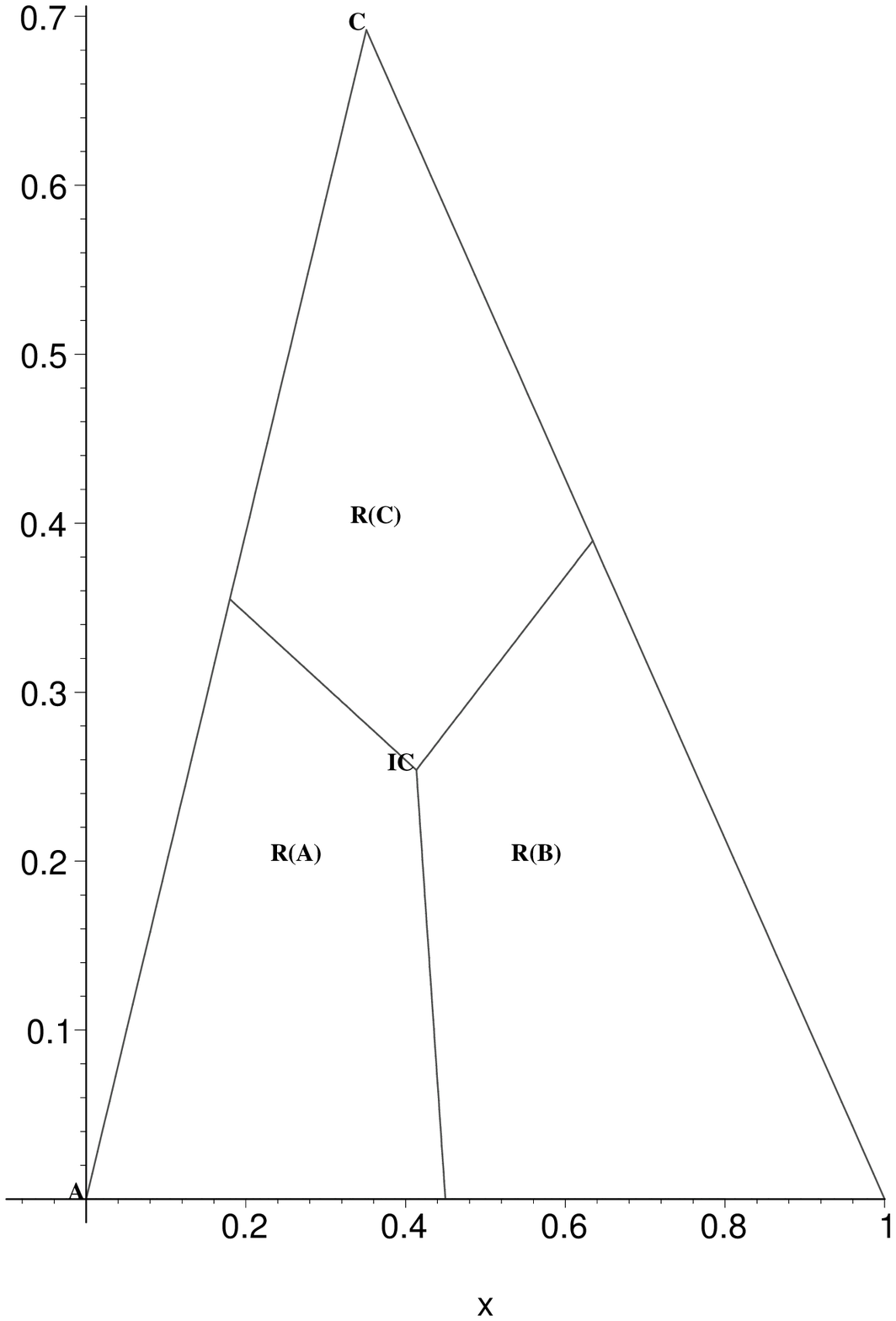, height=140pt,width=200pt}
\end{center}
\caption{
 \label{fig:vert-regions}
The vertex regions constructed with center of mass $M=M_C$ (left) and
incenter $M=M_I$ (right) using the line segments on the line joining each vertex in $\Y_3$ to $M$. }
\end{figure}

For $r \in [1,\infty]$, define $\NPE^r(\cdot,M):=N(\cdot,M;r,\Y_3)$
to be the \emph{(parametrized) proportional-edge proximity map} with
$M$-vertex regions as follows (see also Figure \ref{fig:prox-map-def}
with $M=M_C$ and $r=2$).
For $x \in \TY \setminus \Y_3$, let $v(x)
\in \Y_3$ be the vertex whose region contains $x$; i.e., $x \in
R_M(v(x))$. In this article \emph{$M$-vertex regions} are
constructed by the lines joining any point $M \in \R^2 \setminus
\Y_3$ to a point on each of the edges of $\TY$.
Preferably, $M$ is selected to be in the interior of the triangle $\TY^o$.
For such an $M$, the corresponding vertex regions can be defined using the line
segment joining $M$ to $e_j$, which lies on the line joining $\y_j$ to $M$;
e.g., see Figure \ref{fig:vert-regions} (left) for vertex
regions based on center of mass $M_C$, and Figure \ref{fig:vert-regions} (right) for vertex
regions based on incenter $M_I$.
With $M_C$, the lines joining $M$ and $\Y_3$ are the \emph{median lines},
that cross edges at $M_j$ for $j=1,2,3$.
\emph{$M$-vertex regions}, among many possibilities, can also be defined by the orthogonal
projections from $M$ to the edges.
See \cite{ceyhan:Phd-thesis} for a more general definition.
The vertex regions in Figure \ref{fig:prox-map-def} are
center of mass vertex regions (i.e., $CM$-vertex regions).
If $x$ falls on the boundary of two $M$-vertex regions,
we assign $v(x)$ arbitrarily.
Let $e(x)$ be the edge of $\TY$ opposite of $v(x)$.
Let $\ell(v(x),x)$ be the line parallel to $e(x)$ and passes through $x$.
Let $d(v(x),\ell(v(x),x))$ be the Euclidean (perpendicular)
distance from $v(x)$ to $\ell(v(x),x)$.
For $r \in [1,\infty)$,
let $\ell_r(v(x),x)$ be the line parallel to $e(x)$ such that
$$
d(v(x),\ell_r(v(x),x)) = r\,d(v(x),\ell(v(x),x))\\
\text{ and }\\
d(\ell(v(x),x),\ell_r(v(x),x)) < d(v(x),\ell_r(v(x),x)).
$$
Let $T_r(x)$ be the triangle similar to and with the same
orientation as $\TY$ having $v(x)$ as a vertex and $\ell_r(v(x),x)$
as the opposite edge.
Then the \emph{proportional-edge proximity region} $\NPE^r(x,M)$ is defined to be $T_r(x) \cap \TY$.
Notice that $\ell(v(x),x)$ divides the edges of $T_r(x)$ (other than
the one lies on $\ell_r(v(x),x)$) proportionally with the factor $r$.
Hence the name \emph{proportional-edge proximity region}.

Notice that $r \ge 1$ implies $x \in \NPE^r(x,M)$ for all $x \in \TY$.
Furthermore,
$\lim_{r \rightarrow \infty} \NPE^r(x,M) = \TY$ for all $x \in \TY \setminus \Y_3$,
so we define $\NPE^{\infty}(x,M) = \TY$ for all such $x$.
For $x \in \Y_3$, we define $\NPE^r(x,M) = \{x\}$ for all $r \in [1,\infty]$.

\begin{figure}[ht]
\centering
\rotatebox{-90}{ \resizebox{2. in}{!}{\includegraphics{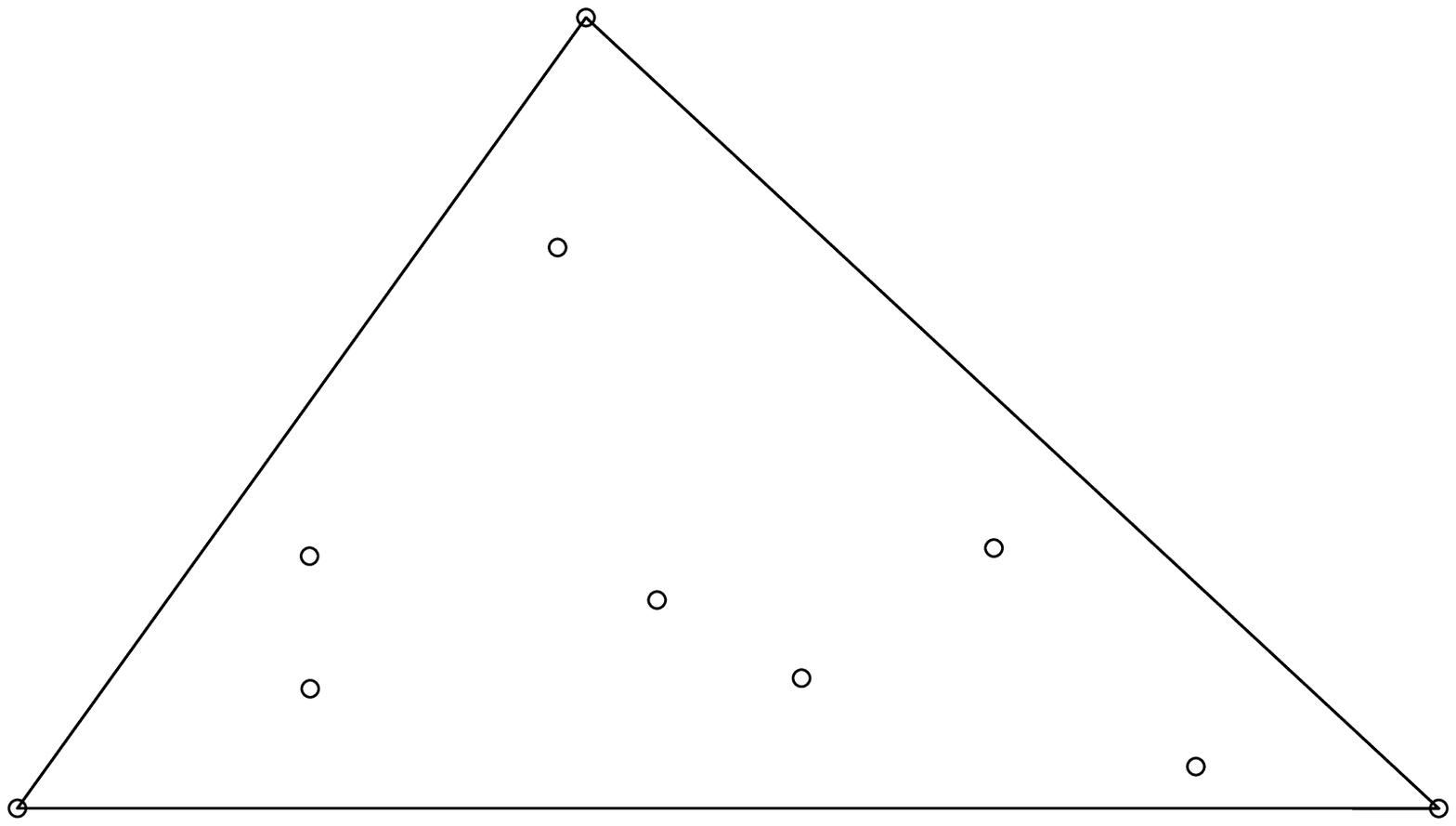}}}
\rotatebox{-90}{ \resizebox{2. in}{!}{\includegraphics{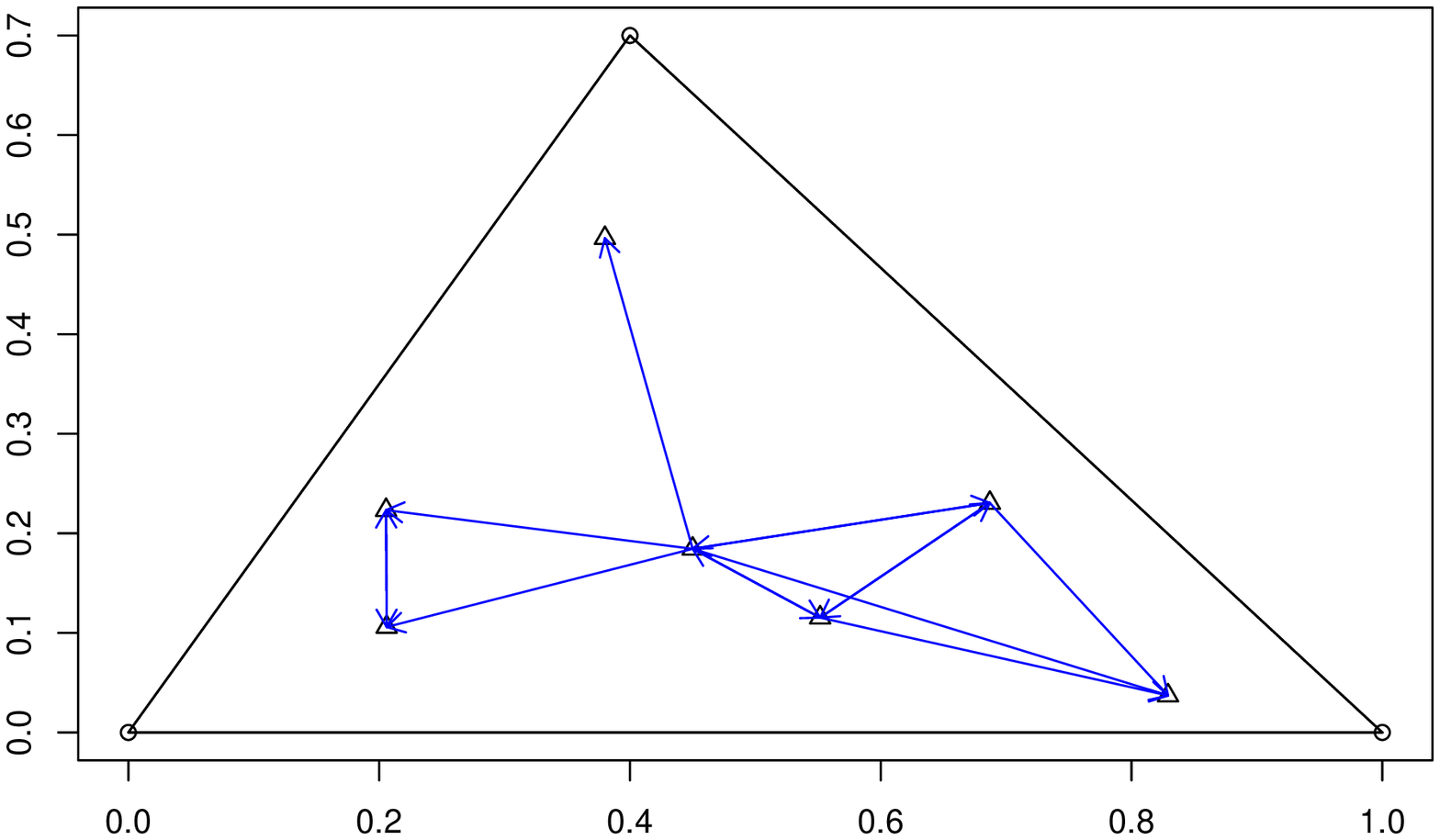}}}
\rotatebox{-90}{ \resizebox{2. in}{!}{\includegraphics{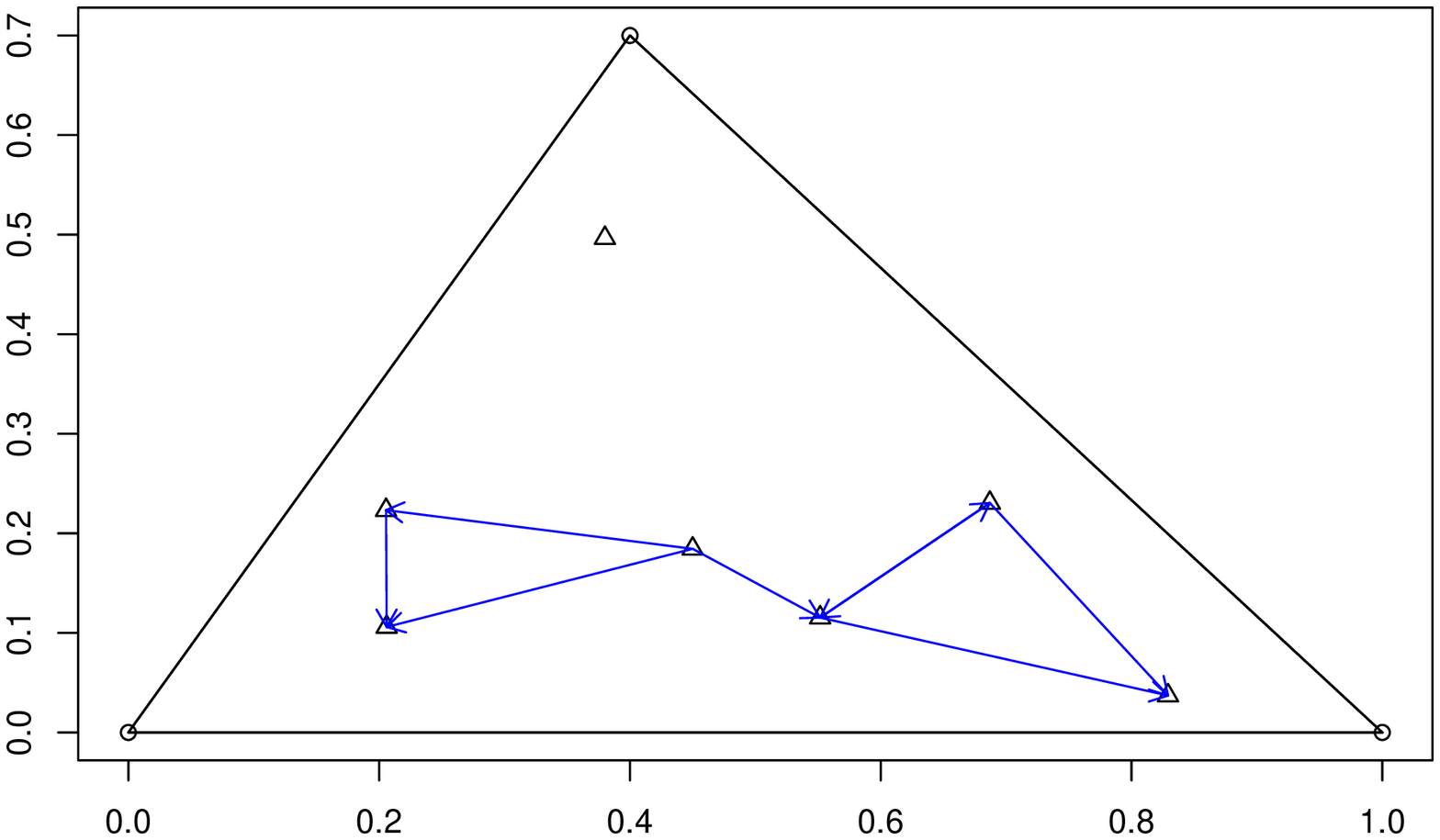}}}
\caption{
\label{fig:one-tri-arcs}
A realization of 7 $\X$ points generated iid $\U(\TY)$, the uniform distribution on $\TY$, (top left)
and the corresponding arcs of proportional-edge PCD with $M=M_C$ for $r=3/2$ (top right) and $r=5/4$ (bottom). }
\end{figure}

The proportional-edge PCD has vertices $\X_n$ and arcs $(x_i,x_j)$ iff $x_j \in \NPE^r(x_i,M)$.
See Figure \ref{fig:one-tri-arcs} for a realization of $\X_n$ with $n=7$ in one triangle (i.e., $m=3$).
For $r=3/2$, the number of arcs is 12 and the domination number $\bgam \left(\X_n,\NPE^{r=3/2}\right)=1$;
and for $r=5/4$, the number of arcs is 9 and $\bgam \left(\X_n,\NPE^{r=5/4}\right)=3$.
By construction, note that as $x$ gets closer to $M$ (or equivalently
further away from the vertices in vertex regions), $\NPE^r(x,M)$
increases in area, hence it is more likely for the outdegree of $x$ to increase.
So if more $\X$ points are around the center $M$,
then it is more likely for the domination number $\g(\X_n,\NPE^r)$ to decrease;
on the other hand, if more $\X$ points are around the vertices $\Y_3$,
then the regions get smaller,
hence it is more likely for the outdegree for such points to be smaller,
thereby implying $\g(\X_n,\NPE^r)$ to increase.
We exploit this probabilistic behavior of $\g(\X_n,\NPE^r)$ in testing spatial patterns
of segregation and association.

Note also that, $\NPE^r(x,M)$ can be viewed as a \emph{homothetic transformation
(enlargement)} with $r \ge 1$ applied on a translation of the region $\NPE^{r=1}(x,M)$.
Furthermore, this transformation is also an \emph{affine similarity transformation}.

\subsection{Some Auxiliary Tools Associated with PCDs}
\label{sec:auxiliary-tools}
\begin{figure} [ht]
    \centering
    \scalebox{.35}{\input{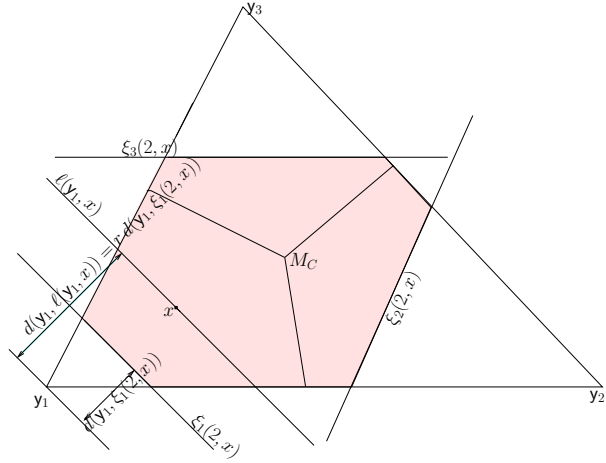}}
    \caption{Construction of the $\G_1$-region, $\G_1^2(x,M_C)$ (shaded region). }
\label{fig:g1-region-def}
    \end{figure}

First, notice that $\NPE^r(x,M)$ is similar to $\TY$ with the similarity ratio being equal to
$$ \psfrag{r}{\normalsize{$r$}} \frac{\min \Bigl(d\bigl(v(x),\,e(x)\bigr),r\,d\bigl(v(x),\,
\ell(v(x),x)\bigr)\Bigr)} {d(v(x),\,e(x))}.$$ 
To define the $\G_1$\emph{-region}, let $\xi_i(x)$ be the line such that
$\xi_i(x)\cap T\left(\Y_3 \right) \not=\emptyset$ and $r\,d(\y_i,\xi_i(x))=d(\y_i,\ell(\y_i,x))$  for $i=1,2,3$.
See also Figure \ref{fig:g1-region-def}.
Then $\G_1^r(x,M)=\bigcup_{i=1}^3 \bigl(\G_1^r(x,M)\cap R_M(\y_i) \bigr)$ where
$\G_1^r(x,M)\cap R_M(\y_i)=\{z \in R_M(\y_i): d(\y_i,\ell(\y_i,z)) \ge d(\y_i,\xi_i(x)\}$, for $i=1,2,3$.
Notice that $r \ge 1$ implies $x \in \G_1^r(x,M)$.
Furthermore, 
$\lim_{r \rightarrow \infty} \G_1^r(x,M) = T\left(\Y_3 \right)$
for all $x \in T\left(\Y_3 \right) \setminus \Y$,
and so we define 
$\G_1^{\infty}(x,M) = T\left(\Y_3 \right)$ for all such $x$.

For $X_i \stackrel{iid}{\sim} F$,
with the additional assumption
that the non-degenerate two-dimensional
probability density function $f$ exists
with $\support(f) \subseteq T\left(\Y_3 \right)$,
implies that the special case in the construction of $\NPE^r$
---$X$ falls on the boundary of two vertex regions ---
occurs with probability zero.
Note that for such an $F$, $\NPE^r(x,M)$ is a triangle a.s.
and $\G^r_1(x,M)$ is a star-shaped (not necessarily convex) polygon.

Let $X_e:=\argmin_{X \in \X_n}d(X,e)$ be the (closest) \emph{edge extremum} for edge $e$
(i.e., closest point among $\X_n$ to edge $e$).
Then it is easily seen that $\G^r_1(\X_n,M)=\bigcap_{i=1}^3\G^r_1(X_{e_i},M)$,
where $e_i$ is the edge opposite vertex $\y_i$, for $i=1,2,3$.
So $\G^r_1(\X_n,M)\cap R_M(\y_i)=\{z \in R_M(\y_i):\; d(\y_i, \ell(\y_i,z) \ge d(\y_i,\xi_i(X_{e_i}))\}$, for $i=1,2,3$.

Let the domination number be $\g_n(r,F,M):=\g_n(\X_n;F,\NPE^r)$ and
$X_{[i,1]}:=\argmin_{X\in \X_n \cap R_M(\y_i)}d(X,e_i)$.
Then $\g_n(r,M) \le 3$ with probability 1,
since $\X_n \cap R_M(\y_i) \subset \NPE^r\left( X_{[i,1]},M \right)$ for each of $i=1,2,3$.
Thus
$$1 \le \E[\g_n(r,F,M)] \le 3~~ \text{ and } ~~0 \le \Var[\g_n(r,F,M)] \le 9/4.$$

\begin{figure}
\begin{center}
\psfrag{A}{\scriptsize{$\y_1$}}
\psfrag{B}{\scriptsize{$\y_2$}}
\psfrag{C}{\scriptsize{$\y_3$}}
\psfrag{CC}{}
\psfrag{IC}{}
\psfrag{x}{\scriptsize{$x$}}
\psfrag{D}{}
\psfrag{E}{}
\psfrag{F}{}
\psfrag{S}{}
\psfrag{Ac}{\scriptsize{$q_1(r,x)$} }
\psfrag{Ab}{}
\psfrag{Ba}{}
\psfrag{Bc}{\scriptsize{$q_2(r,x)$} }
\psfrag{Ca}{}
\psfrag{Cb}{\scriptsize{$q_3(r,x)$} }
 \epsfig{figure=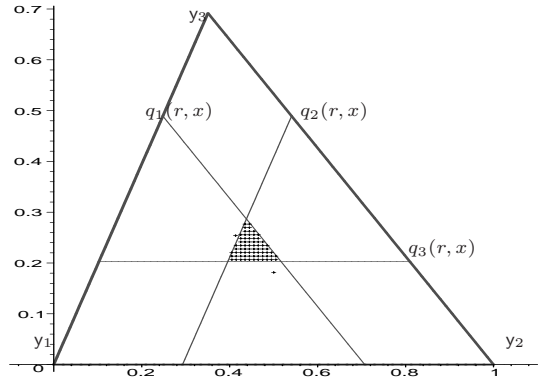, height=140pt , width=200pt}
\end{center}
\caption{The triangle $\Tr$ with $r=\sqrt{2}$ (the hatched region).
} \label{fig:Tr-RS-NDA-CC}
\end{figure}

In $\TY$, drawing the lines $q_i(r,x)$ such that
$d(\y_i,e_i)=r\,d(\y_i,q_i(r,x))$ for $i\in \{1,2,3\}$  yields another triangle,
denoted as $\Tr$, for $r<3/2$.
See Figure \ref{fig:Tr-RS-NDA-CC} for $\Tr$ with $r=\sqrt{2}$.

The functional form of $\Tr$ in the basic triangle $T_b$ is given by
{\small
\begin{align}
\label{eqn:T^r-def}
& \Tr =T(t_1(r),t_2(r),t_3(r))= \left \{(x,y) \in T_b:
y \ge \frac{c_2\,(r-1)}{r};\; y \le \frac{c_2\,(1-r\,x)}{r\,(1-c_1)};\;
y \le \frac{c_2\,(r\,(x-1)+1)}{r\,c_1} \right\}\\
&=T\Biggl( \left(\frac{(r-1)\,(1+c_1)}{r},\frac{c_2\,(r-1)}{r} \right),
\left(\frac{2-r+c_1\,(r-1)}{r},\frac{c_2\,(r-1)}{r} \right),
\left(\frac{c_1\,(2-r)+r-1}{r},\frac{c_2\,(2-r)}{r} \right) \Biggr) \nonumber.
\end{align}
}

In the standard equilateral triangle, this functional form becomes:
$$
\Tr = T\Biggl( \left(\frac{3\,(r-1)}{2\,r},\frac{\sqrt{3}(r-1)}{2\,r} \right),
\left(\frac{3-r}{2\,r},\frac{\sqrt{3}(r-1)}{2\,r} \right),
\left(\frac{1}{2},\frac{\sqrt{3}(2-r)}{r} \right) \Biggr).
$$

\begin{figure}
\begin{center}
\psfrag{A}{\scriptsize{$\y_1$}}
 \psfrag{B}{\scriptsize{$\y_2$}}
\psfrag{C}{\scriptsize{$\y_3$}}
 \psfrag{CC}{\scriptsize{$M_{CC}$}}
 \psfrag{x}{}
 \psfrag{y}{}
  \psfrag{M1}{}
 \psfrag{M2}{}
 \psfrag{M3}{}
 \epsfig{figure=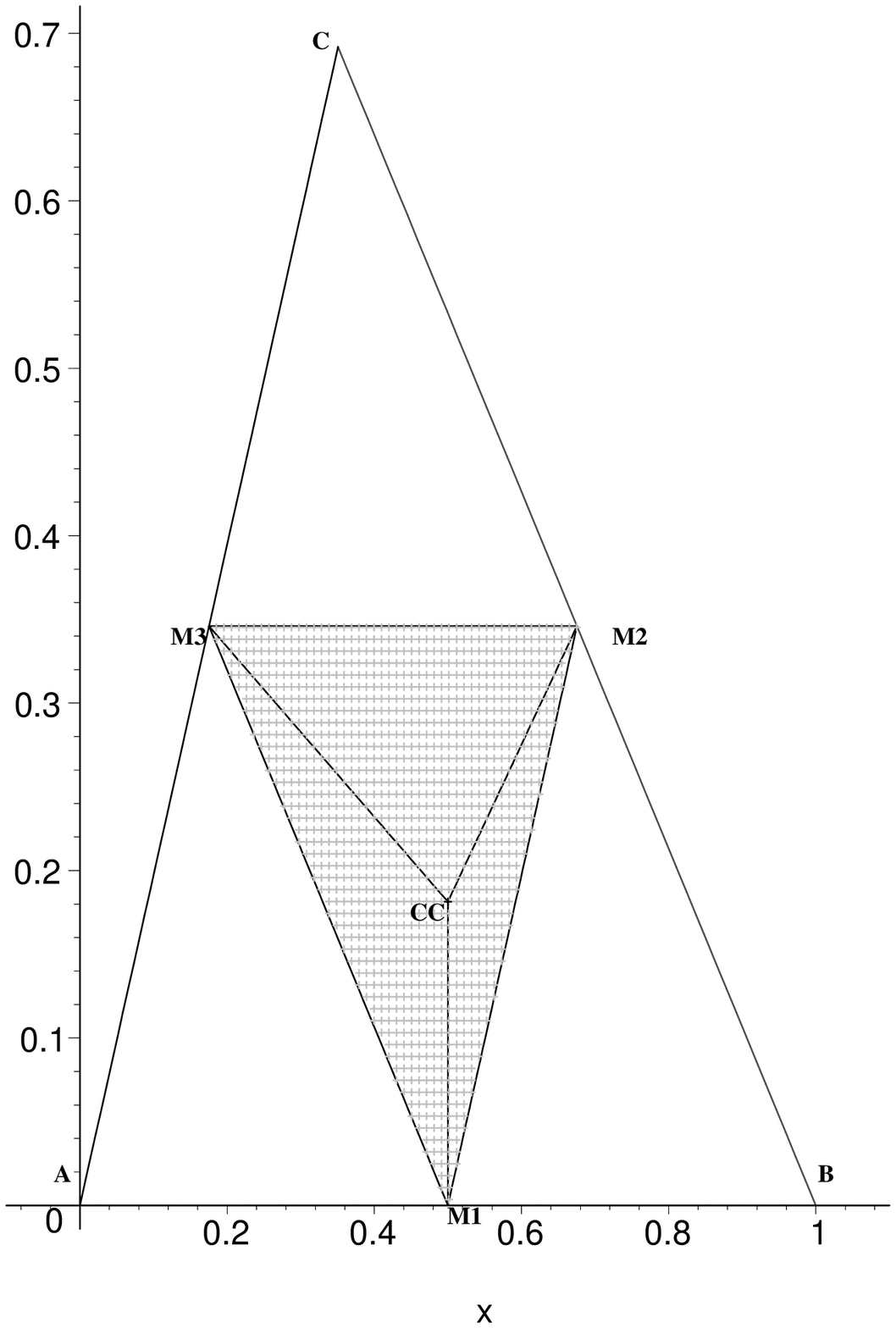, height=140pt,width=200pt}
 \epsfig{figure=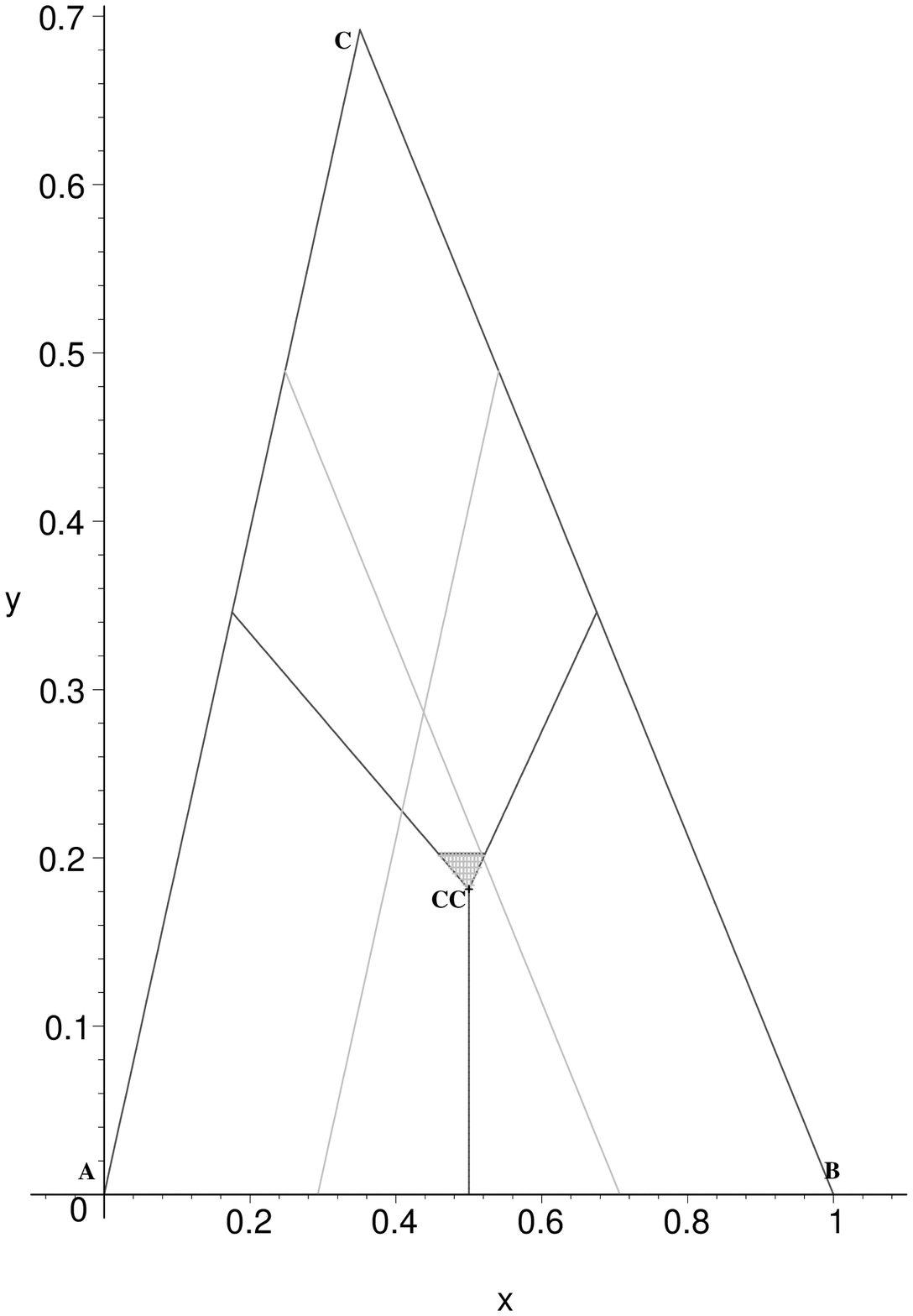, height=140pt,width=200pt}
\end{center}
\caption{
 \label{fig:SS-regions}
The superset regions (the shaded regions) constructed with
circumcenter $M_{CC}$ with $r=\sqrt{2}$ (left) and $r=2$ (right)
with vertex regions constructed with orthogonal projections to the
edges. }
\end{figure}

There is a crucial difference between the triangles $\Tr$ and $T(M_1,M_2,M_3)$.
More specifically $T(M_1,M_2,M_3) \subseteq \RS(r,M)$ for all $M$ and $r \ge 2$, but
$(\Tr)^o$ and $\RS(r,M)$ are disjoint for all $M$ and $r$.
So if $M \in (\Tr)^o$, then $\RS(r,M)=\emptyset$;
if $M \in \partial(\Tr)$, then $\RS(r,M)=\{M\}$; and if $M \not\in \Tr$,
then $\RS(r,M)$ has positive area.
See Figure \ref{fig:SS-regions} for two examples of superset regions
with $M$ that corresponds to circumcenter $M_{CC}$ in this triangle
and the vertex regions are constructed using orthogonal projections.
For $r=2$, note that $\Tr=\emptyset$ and the superset region is
$T(M_1,M_2,M_3)$ (see Figure \ref{fig:SS-regions} (left)),
while for $r=\sqrt{2}$, $\Tr^o$ and $\RS(r=\sqrt{2},M)^o$ are disjoint
(see Figure \ref{fig:SS-regions} (right))

The triangle $\Tr$ given in Equation \eqref{eqn:T^r-def} plays an important role in the
distribution of the domination number of the proportional-edge PCDs.

\section{The Asymptotic Distribution of Domination Number for Uniform Data}
\label{sec:asy-gam-NYr}

\subsection{The One-Triangle Case}
\label{sec:asy-gam-NYr-one-tri}
For simplicity, we consider $\X$ points iid uniform in one triangle only.
The null hypothesis we consider is a type of
{\em complete spatial randomness}; that is,
$$H_o: X_i \stackrel{iid}{\sim} \U(T\left(\Y_3 \right)) \text{ for } i=1,2,\ldots,n,$$
where $\U(T\left(\Y_3 \right))$ is the uniform distribution on $T\left(\Y_3 \right)$.
If it is desired to have the sample size be a random variable,
we may consider a spatial Poisson point process on $T\left(\Y_3 \right)$
as our null hypothesis.
Let $\g_n(r,M):=\g\left(\X_n,\UT,\NPE^r,M \right)$ be the domination number
of the PCD based on $\NPE^r$ with $\X_n$, a set of iid random
variables from $\UT$, with $M$-vertex regions.

We present a ``geometry invariance" result for $\NPE^r(\cdot,M)$
where \emph{$M$-vertex regions} are constructed using the line segment joining $M$
to edge $e_i$ on the line joining $\y_i$ to $M$,
rather than the orthogonal projections from $M$ to the edges.
This invariance property will simplify the notation in
our subsequent analysis by allowing us to consider the special case
of the (standard) equilateral triangle.

\begin{theorem}
\label{thm:geo-inv-NYr}
(Geometry Invariance Property)
Suppose $\X_n$ is a set of iid random variables from $\UT$.
Then for any $r \in [1,\infty]$ the distribution of $\g_n(r,M)$ is
independent of $\Y_3$ and hence the geometry of $\TY$.
\end{theorem}

\noindent
{\bf Proof:}
See \cite{ceyhan:dom-num-NPE-MASA} for the proof. $\blacksquare$

Note that geometry invariance of $\g_n(r=\infty,M)$ follows trivially
for all $\X_n$ from any $F$ with support in $\TY \setminus \Y_3$, since for $r=\infty$,
we have $\g_n(r=\infty,M)=1$ a.s.
Based on Theorem \ref{thm:geo-inv-NYr} we may assume that $\TY$ is a
standard equilateral triangle with $\Y_3= \{ (0,0),(1,0),\left( 1/2,\sqrt{3}/2 \right)\}$
for $\NPE^r(\cdot,M)$ with $M$-vertex regions.

\begin{remark}
Notice that, we proved the geometry invariance property for $\NPE^r(\cdot)$
where $M$-vertex regions are defined with the lines joining $\Y_3$ to $M$.
On the other hand, if we use the orthogonal projections from
$M$ to the edges, the vertex regions (hence $\NPE^r$) will depend on
the geometry of the triangle.
That is, the orthogonal projections
from $M$ to the edges will not be mapped to the orthogonal
projections in the standard equilateral triangle.
Hence with the
orthogonal projections, the exact and
asymptotic distribution of $\g_n(r,M)$ will depend on $c_1,c_2$ of $T_b$,
so one needs to do the calculations for each possible combination of $c_1,c_2$.
$\square$
\end{remark}

\begin{figure}
\begin{center}
\psfrag{r}{\scriptsize{$r$}}
\psfrag{pr}{\scriptsize{$p_r$}}
\epsfig{figure=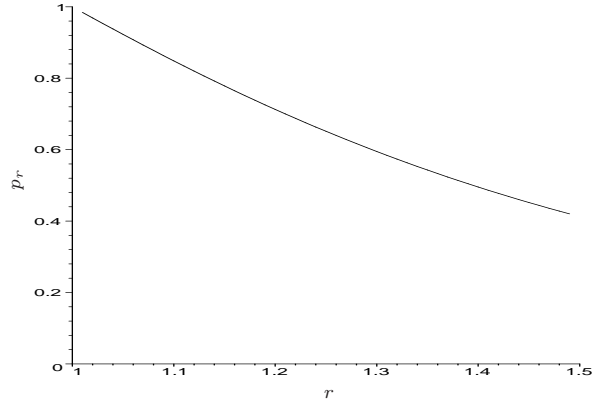, height=150pt, width=220pt}
\end{center}
\caption{ Plotted is the probability $p_r=\lim_{n\rightarrow
\infty}P\left( \g_n(r,M)=2 \right)$ given in Equation \eqref{eqn:p_r-form} as a function of $r$ for $r \in [1,3/2)$ and $M
\in \{t_1(r),t_2(r),t_3(r)\}$.} \label{fig:pglt2ofr}
\end{figure}

The domination number $\g_n(r,M)$ of the PCD has the following
asymptotic distribution (\cite{ceyhan:dom-num-NPE-MASA}).
As $n \rightarrow \infty$,
\begin{equation}
\label{eqn:asymptotic-NYr}
\g_n(r,M) \stackrel{\mathcal L}{\longrightarrow}
\left\lbrace \begin{array}{ll}
       2+\BER(1-p_r)& \text{for $r \in [1,3/2)$ and $M \in \{t_1(r),t_2(r),t_3(r)\}$,}\\
       1            & \text{for $r>3/2$ and $M \in \TY^o$,}\\
       3            & \text{for $r \in [1,3/2)$ and $M \in \Tr\setminus \{t_1(r),t_2(r),t_3(r)\}$,}\\
\end{array} \right.
\end{equation}
where $\stackrel{\mathcal L}{\longrightarrow}$ stands for ``convergence in law"
and $\BER(p)$ stands for Bernoulli distribution with probability of
success $p$, $\Tr$ and $t_i(r)$ are defined in Equation \eqref{eqn:T^r-def}, and for $r \in [1,3/2)$ and $M \in
\{t_1(r),t_2(r),t_3(r)\}$,
\begin{equation}
\label{eqn:p_r-form}
p_r=\int_0^{\infty}\int_0^{\infty}\frac {64\,r^2}{9\,(r-1)^2}\,w_1\,w_3\,\exp\left(\frac{4\,r}
{3\,(r-1)}\,(w_1^2+w_3^2+2\,r\,(r-1)\,w_1\,w_3)\right)\,dw_3w_1,
\end{equation}
and for $r=3/2$ and $M=M_C=\left\{(1/2,\sqrt{3}/6)\right\}$, $p_r \approx 0.7413$,
which is not computed as in Equation \eqref{eqn:p_r-form};
for its computation, see \cite{ceyhan:dom-num-NPE-SPL}.
For example, for $r=5/4$ and
$M \in \left\{t_1(r)=\left(3/10,\sqrt{3}/10\right),t_2(r)=\\
\left(7/10,\sqrt{3}/10\right),t_3(r)=\left(1/2,3\sqrt{3}/5\right)\right\}$,
$p_r \approx 0.6514$.
See Figure \ref{fig:pglt2ofr} for the plot of the numerically computed values
(i.e., the values computed by numerical integration) of $p_r$ as a function of $r$
according to Equation \eqref{eqn:p_r-form}.
Notice that in the nondegenerate case in \eqref{eqn:asymptotic-NYr},
$\E[\g_n(r,M)]=3-p_r$ and $\Var[\g_n(r,M)]=p_r(1-p_r)$.

In Equation (\ref{eqn:asymptotic-NYr}), the first line is referred
as the non-degenerate case, the second and third lines are referred
as degenerate cases with a.s. limits 1 and 3, respectively.

\begin{table}[ht]
\begin{center}
\begin{tabular}{cc}

\begin{tabular}{|c|c|c|c|c|c|}
\hline
\multicolumn{6}{|c|}{$r=2$ and $M=M_C$} \\
\hline
$ k \diagdown n$  & 10 & 20 & 30 & 50 & 100\\
\hline
1 & 961 & 1000 & 1000 & 1000 & 1000 \\
\hline
2 & 34 & 0 & 0 & 0 & 0  \\
\hline
3 & 5 & 0 & 0 & 0 & 0 \\
\hline
\end{tabular}
&
\begin{tabular}{|c|c|c|c|c|c|}
\hline
\multicolumn{6}{|c|}{$r=5/4$ and $M=M_C$} \\
\hline
$ k \diagdown n$  & 10 & 20 & 30 & 50 & 100\\
\hline
1 & 9 & 0 & 0 & 0 & 0 \\
\hline
2 & 293 & 110 & 30 & 8 & 0  \\
\hline
3 & 698 & 890 & 970 & 992 & 1000 \\
\hline
\end{tabular}

\end{tabular}
\end{center}
\caption{
\label{tab:numerical-gamma-NYr=2,5/4}
The number of $\g_n(r,M)=k$ out of $N=1000$ Monte Carlo replicates with $M=M_C$
and $r=2$ (left) and $r=5/4$ (right).
Here, ``$r=2$ and $M=M_C$" is an example of the case ``$r>3/2$ and $M \in \TY^o$",
and ``$r=5/4$ and $M=M_C$" is an example of the case ``$r \in [1,3/2)$ and $M \in \Tr\setminus \{t_1(r),t_2(r),t_3(r)\}$".
}
\end{table}

\begin{table}[ht]
\begin{center}
\begin{tabular}{|c|c|c|c|c|c|c|c|c|}
\hline
\multicolumn{9}{|c|}{$r=5/4$ and $M=\left( 3/5,\sqrt{3}/10\right)$} \\
\hline
$ k \diagdown n$  & 10 & 20 & 30 & 50 & 100 & 500 & 1000 & 2000\\
\hline
1 & 118 & 60 & 51 & 39 & 15 & 1 & 2 & 1 \\
\hline
2 & 462 & 409 & 361 & 299 & 258 & 100 & 57 & 29  \\
\hline
3 & 420 & 531 & 588 & 662 & 727 & 899 & 941 & 970\\
\hline
\end{tabular}

\begin{tabular}{|c|c|c|c|c|c|c|c|c|}
\hline
\multicolumn{9}{|c|}{$r=5/4$ and $M=\left(7/10,\sqrt{3}/10 \right)$} \\
\hline
$ k \diagdown n$  & 10 & 20 & 30 & 50 & 100 & 500 & 1000 & 2000\\
\hline
1 & 174 & 118 & 82 & 61 & 22 & 5 & 1 & 1\\
\hline
2 & 532 & 526 & 548 & 561 & 611 & 617 & 633 & 649  \\
\hline
3 & 294 & 356 & 370 & 378 & 367 & 378 & 366 & 350\\
\hline
\end{tabular}
\end{center}
\caption{
\label{tab:numerical-gamma-NYr=5/4}
The number of $\g_n(r,M)=k$ out of $N=1000$ Monte Carlo replicates with $r=5/4$
and $M=\left( 3/5,\sqrt{3}/10\right)$ (top) and $M=\left(7/10,\sqrt{3}/10 \right)$ (bottom).
Here ``$r=5/4$ and $M=(3/5,\sqrt{3}/10)$" is an example of the case
``$r \in [1,3/2)$ and $M \in \Tr\setminus \{t_1(r),t_2(r),t_3(r)\}$"
with $M$ being on the line segment joining $t_1(r)$ and $t_2(r)$,
and ``$r=5/4$ and $M=(7/10,\sqrt{3}/10)$" is an example of the case
``$r \in [1,3/2)$ and $M \in \Tr\setminus \{t_1(r),t_2(r),t_3(r)\}$"
with $M=t_2(r)$.
}
\end{table}

We also estimate the distribution of $\g_n(r,M)$ for various values of
$n$, $r$, and $M$ using Monte Carlo simulations.
At each Monte Carlo replication,
we generate $n$ points iid $\U(\TY)$ and compute the value of $\g_n(r,M)$.
The frequencies of $\g_n(r,M)=k$ out of $N=1000$ Monte Carlo replicates are presented
in Tables \ref{tab:numerical-gamma-NYr=2,5/4}, \ref{tab:numerical-gamma-NYr=5/4},
and \ref{tab:numerical-gamma-NYr=3/2}.
Notice that in Table \ref{tab:numerical-gamma-NYr=2,5/4} (left)
``$r=2$ and $M=M_C$" is an example of the case ``$r>3/2$ and $M \in \TY^o$",
in Table \ref{tab:numerical-gamma-NYr=2,5/4} (right)
``$r=5/4$ and $M=M_C$" is an example of the case
``$r \in [1,3/2)$ and $M \in \Tr\setminus \{t_1(r),t_2(r),t_3(r)\}$";
in Table \ref{tab:numerical-gamma-NYr=5/4} (top)
``$r=5/4$ and $M=(3/5,\sqrt{3}/10)$" is an example of the case
``$r \in [1,3/2)$ and $M \in \Tr\setminus \{t_1(r),t_2(r),t_3(r)\}$"
with $M$ being on the line segment joining $t_1(r)$ and $t_2(r)$,
in Table \ref{tab:numerical-gamma-NYr=5/4} (bottom)
``$r=5/4$ and $M=(7/10,\sqrt{3}/10)$" is an example of the case
``$r \in [1,3/2)$ and $M \in \Tr\setminus \{t_1(r),t_2(r),t_3(r)\}$"
with $M=t_2(r)$;
and in Table \ref{tab:numerical-gamma-NYr=3/2},
``$r=3/2$ and $M=M_C$" is an example of the case discussed in (\cite{ceyhan:dom-num-NPE-SPL}).
Notice that as the sample size $n$ increases, the values on these tables
get closer and closer to the expected values under their asymptotic distribution.

\begin{table}[ht]
\begin{center}
\begin{tabular}{|c|c|c|c|c|c|c|c|c|}
\hline
\multicolumn{9}{|c|}{$r=3/2$ and $M=M_C$} \\
\hline
$ k \diagdown n$  & 10 & 20 & 30 &  50 & 100 & 500 & 1000 & 2000\\
\hline
1  & 151 & 82 & 61 & 50 & 27 & 2 & 3 & 1\\
\hline
2 &  602 & 636 & 688 & 693 & 718 & 753  & 729 & 749\\
\hline
3  & 247 & 282 & 251 & 257 & 255 & 245 & 268 & 250 \\
\hline
\end{tabular}
\end{center}
\caption{\label{tab:numerical-gamma-NYr=3/2}
The number of $\g_n(3/2,M_C)=k$ out of $N=1000$ Monte Carlo replicates.
Here ``$r=3/2$ and $M=M_C$" is an example of the case discussed in (\cite{ceyhan:dom-num-NPE-SPL}).}
\end{table}

\begin{theorem}
\label{thm:stoch-order-1tri}
Let $\g_n(r,M) = \g(\X_n;\U(T\left(\Y_3 \right)),\NPE^r,M)$.
Then $r_1<r_2$ implies $\g_n(r_2,M) <^{ST} \g_n(r_1,M)$
where $<^{ST}$ stands for ``stochastically smaller than".
\end{theorem}

\noindent
{\bf Proof:}
Suppose $r_1<r_2$.
Then $P(\g_n(r_2,M)=1)>P(\g_n(r_1,M)=1)$ and $P(\g_n(r_2,M)=2)>P(\g_n(r_1,M)=2)$ and $P(\g_n(r_2,M)=3)<P(\g_n(r_1,M)=3)$.
Hence the desired result follows. $\blacksquare$

\subsection{The Multi-Triangle Case}
\label{sec:NYr-mult-tri}
In this section, we present the asymptotic
distribution of the domination number of proportional-edge PCDs in multiple Delaunay triangles.
Suppose $\Y_m=\{\y_1,\y_2,\ldots,\y_m\} \subset
\R^2$ be a set of $m$ points in general position with $m>3$ and no
more than 3 points are cocircular.
Then there are $J_m>1$ Delaunay triangles each of which is denoted as $T_j$ (\cite{okabe:2000}).
We wish to investigate
\begin{equation}
\label{eqn:null-pattern-mult-tri}
H_o: X_i \stackrel{iid}{\sim} \U(C_H(\Y_m))  \text{ for } i=1,2,\ldots,n
\end{equation}
against segregation and association alternatives
(see Section \ref{sec:alternatives}).
Figure \ref{fig:deldata} (middle) presents a realization of 1000 observations
independent and identically distributed according to $\U(C_H(\Y_m))$ for $m=10$ and $J_m=13$.

\begin{figure}[ht]
\centering
\rotatebox{-90}{ \resizebox{3. in}{!}{\includegraphics{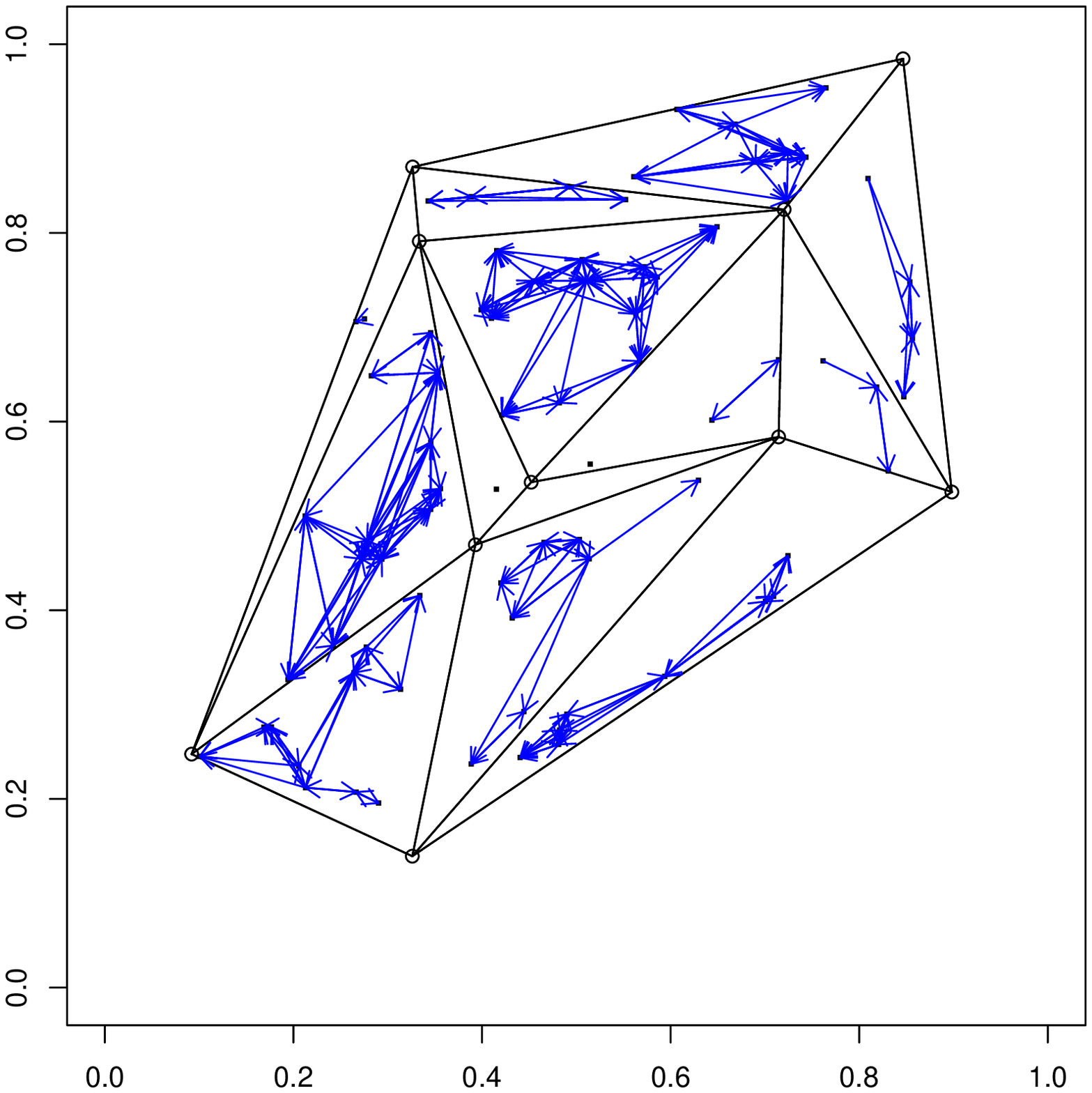}}}
\rotatebox{-90}{ \resizebox{3. in}{!}{\includegraphics{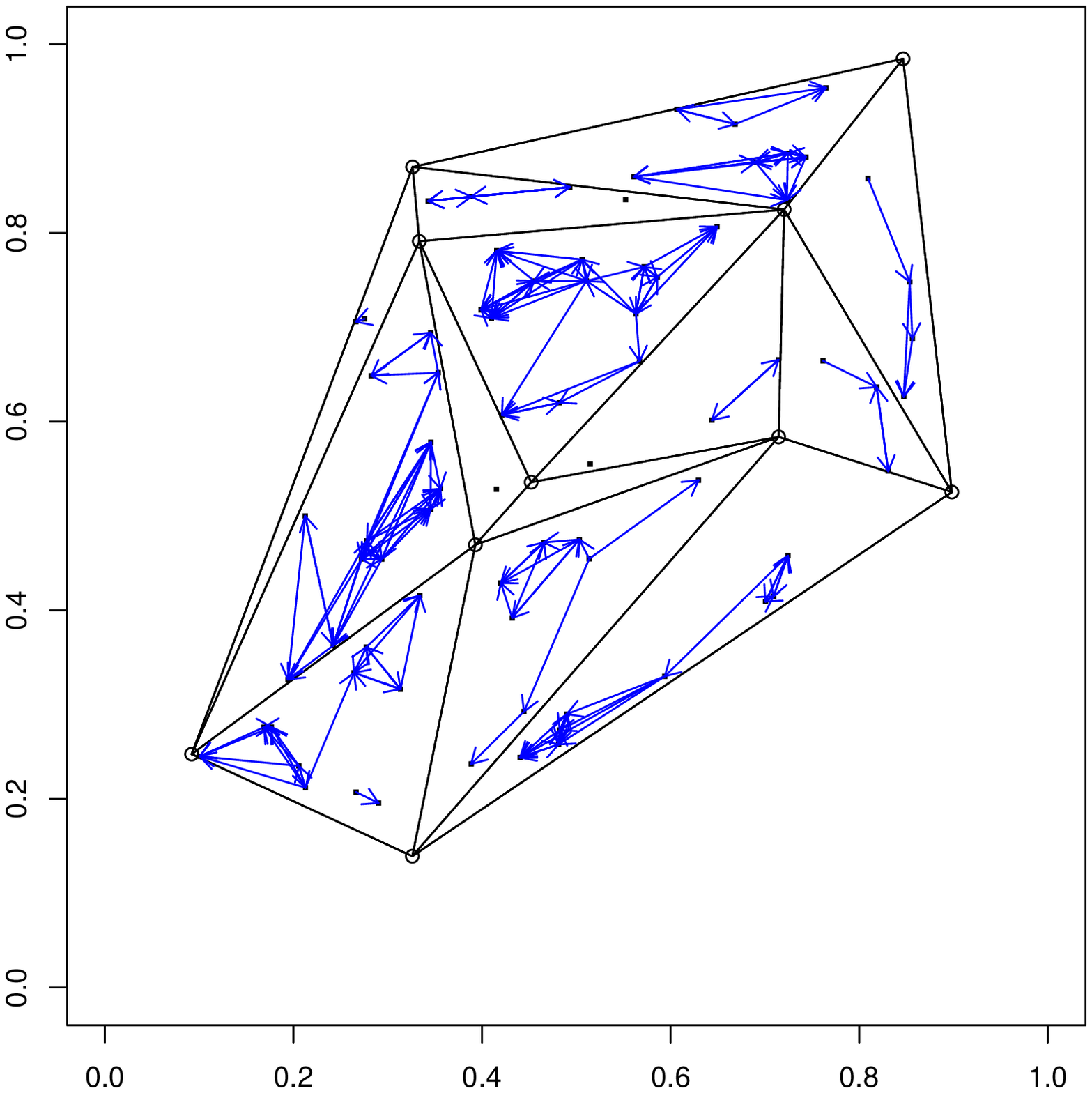}}}
\caption{
\label{fig:multi-tri-arcs}
The arcs for the 77 $\X$ points (dots, $\centerdot$) in the convex hull of $\Y$ points (circles, $\circ$)
given in Figure \ref{fig:deltri} for the proportional-edge PCD with $M=M_C$ for $r=3/2$ (left) and $r=5/4$ (right).}
\end{figure}

Let $M^j$ be the point in $T_j$ that corresponds to $M$ in $T_e$,
$\Tr^j$ be the triangle that corresponds to $\Tr$ in $T_e$,
and $t_i^j(r)$ be the vertices of $\Tr^j$ that correspond to $t_i(r)$ in $T_e$ for $i \in \{1,2,3\}$.
Moreover, let $n_j:=|\X_n \cap T_j|$,
the number of $\X$ points in Delaunay triangle $T_j$.
The digraph $D$ is constructed using $\NPE^r(\cdot,M^j)$ as described above,
where the three points in $\Y_m$ defining the
Delaunay triangle $T_j$ are used as $\Y_{m(j)}$.
Then we have $\ge J_m$ disconnected sub-digraphs.
For $\X_n \subset \C_H(\Y_m)$,
let $\g_{n_j}(r,M^j)$ be the domination number of the digraph
induced by vertices of $T_j$ and $\X_n \cap T_j$.
Then the domination number of the proportional-edge PCD in $J_m$ triangles is
$$\g_{n,m}(r,M)=\sum_{j=1}^{J_m}\g_{n_j}\left(r,M^j\right).$$
See Figure \ref{fig:multi-tri-arcs} for two examples of the proportional edge PCDs based on
the 77 $\X$ points that are in $\C_H(\Y_m)$ out of the 200 $\X$ points plotted in Figure \ref{fig:deltri}.
The arcs are constructed for $M=M_C$ with $r=3/2$ (left) and $r=5/4$ (right)
and the corresponding domination number values are
$\g_{n,10}(3/2,M_C)=22$ and $\g_{n,10}(5/4,M_C)=26$.
Suppose $\X_n$ is a set of iid random variables from $\U(\C_H(\Y_m))$, the uniform distribution on convex
hull of $\Y_m$ and we construct the proportional-edge PCDs using the points
$M^j$ that correspond to $M$ in $T_e$. Then
For fixed $m$ (or fixed $J_m$), as $n \rightarrow \infty$, so does each $n_j$.
Furthermore, as $n \rightarrow \infty$, each component $\g_{n_j}(r,M^j)$ become independent.
Therefore using Equation \eqref{eqn:asymptotic-NYr},
we can obtain the asymptotic distribution of $\g_{n,m}(r,M)$.
For fixed $J_m$, as $n \rightarrow \infty$,\\

\noindent
\begin{equation}
\label{eqn:asymptotic-NYr-Jm}
\g_{n,m}(r,M) \stackrel{\mathcal L}{\longrightarrow}
\left\lbrace
\begin{array}{ll}
       2\,J_m+\BIN(J_m,1-p_r)& \text{for $M^j \in \left\{t^j_1(r),t^j_2(r),t^j_3(r)\right\}$ and $r \in [1,3/2]$,}\\
       J_m                   & \text{for $r>3/2$ and for all $M^j \not= \Y_3$,}\\
       3\,J_m                & \text{for $M \in \Tr^j\setminus \left\{t^j_1(r),t^j_2(r),t^j_3(r)\right\}$ and $r \in [1,3/2)$,}\\
\end{array} \right.
\end{equation}
\noindent where $\BIN(n,p)$ stands for binomial distribution with
$n$ trials and probability of success $p$, for $r \in [1,3/2)$ and
$M \in \{t_1(r),t_2(r),t_3(r)\}$, $p_r$ is given in Equation \eqref{eqn:p_r-form}
and $j=1,2,\ldots,J_m$.
Observe that in the nondegenerate case in Equation \eqref{eqn:asymptotic-NYr-Jm},
we have $\E[\g_{n,m}(r,M)]=J_m(3-p_r)$ and $\Var[\g_{n,m}(r,M)]=J_m\,p_r\,(1-p_r)$.

\begin{theorem}
\label{thm:asy-normality}
(\textbf{Asymptotic Normality})
Suppose $n_j$ and $J_m$ are sufficiently large with $n_j \gg J_M$.
Then the asymptotic null distribution of the mean domination number (per triangle)
$\displaystyle \overline G(r,M):=\frac{1}{J_m}\,\sum_{j=1}^{J_m} \g_{n_j}(r,M)=\frac{\g_{n,m}(r,M)}{J_m}$ is
approximately normal;
i.e., for large $n_j \gg J_M$
$$\displaystyle \overline G(r,M) \stackrel {\text{approx}}{\Large\sim} \mathcal N (\mu,\sigma^2/J_m)$$
where $\mu=3-p_r$ and $\sigma^2=p_r(1-p_r)/J_m$.
\end{theorem}

\noindent
{\bf Proof:} For fixed $J_m$ sufficiently large and each $n_j$ sufficiently large
with $n=\sum_{j=1}^{J_m} n_j \gg J_m$,
$\g_{n_j}(r,M)$ are approximately independent identically distributed as in Equation \eqref{eqn:asymptotic-NYr}.
Then the desired result follows.
$\blacksquare$

\begin{figure}[ht]
\centering
\psfrag{Density}{ \Huge{\bf{Density}}}
\rotatebox{-90}{ \resizebox{2. in}{!}{ \includegraphics{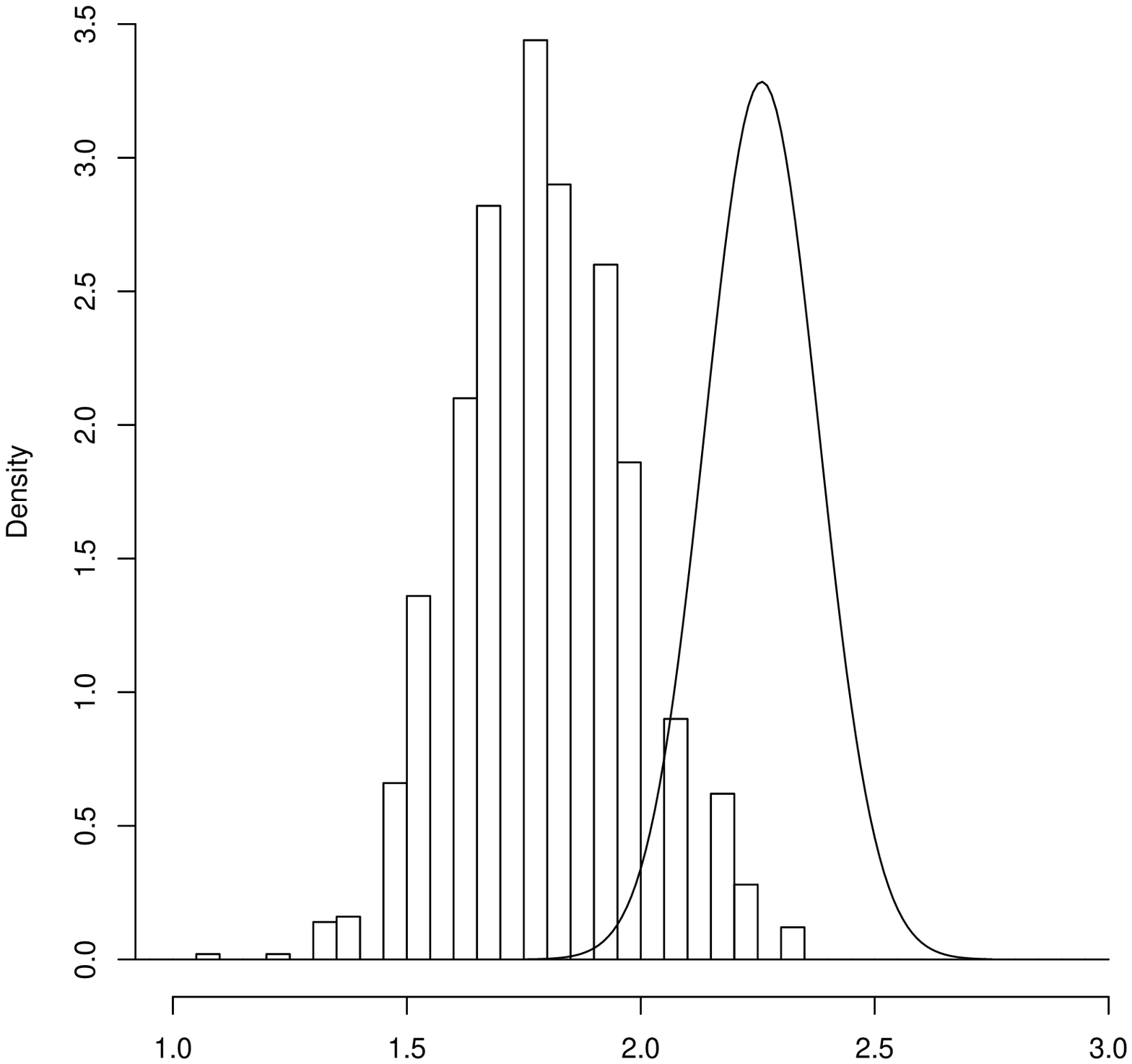} } }
\rotatebox{-90}{ \resizebox{2. in}{!}{ \includegraphics{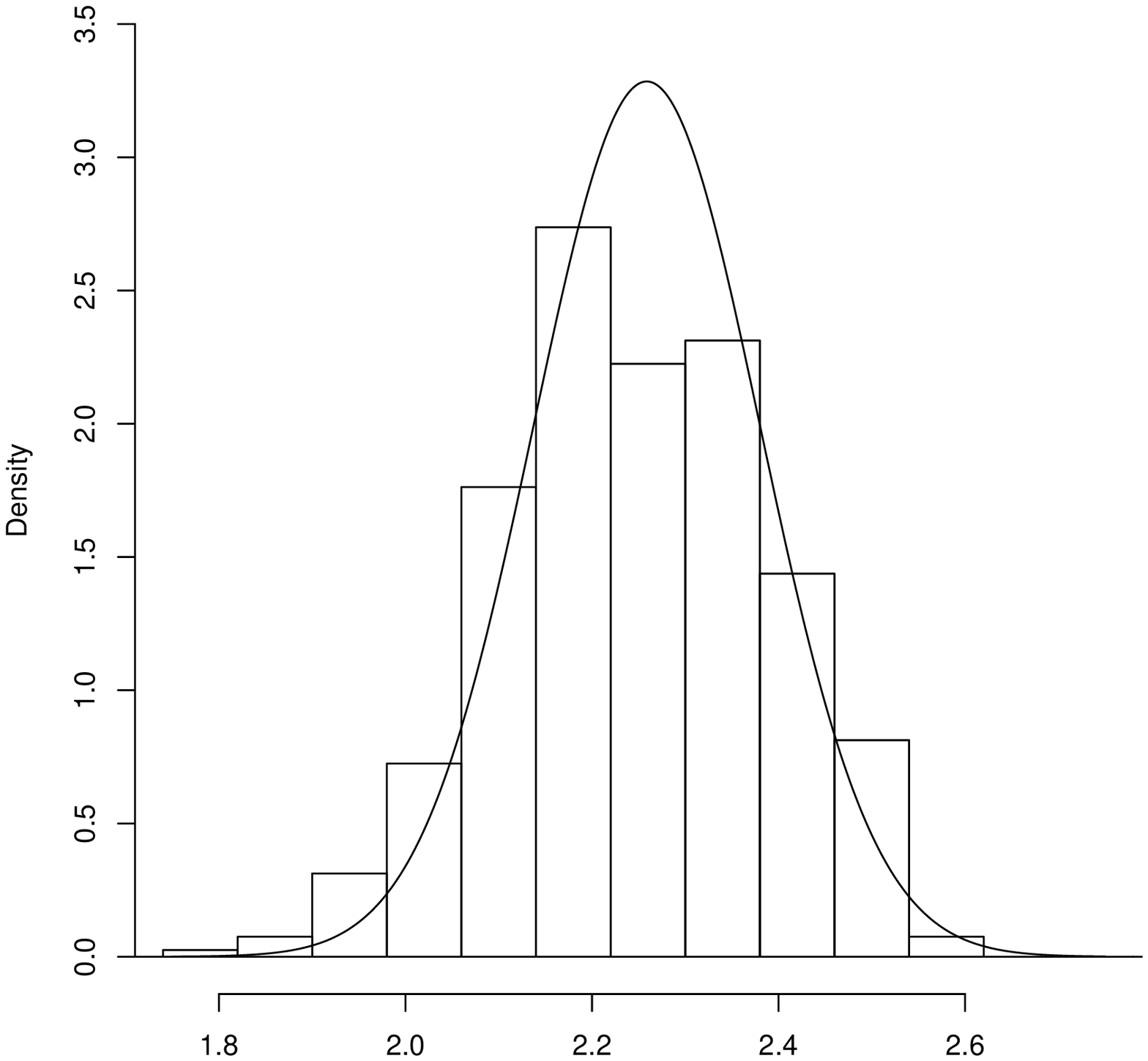} } }
\rotatebox{-90}{ \resizebox{2. in}{!}{ \includegraphics{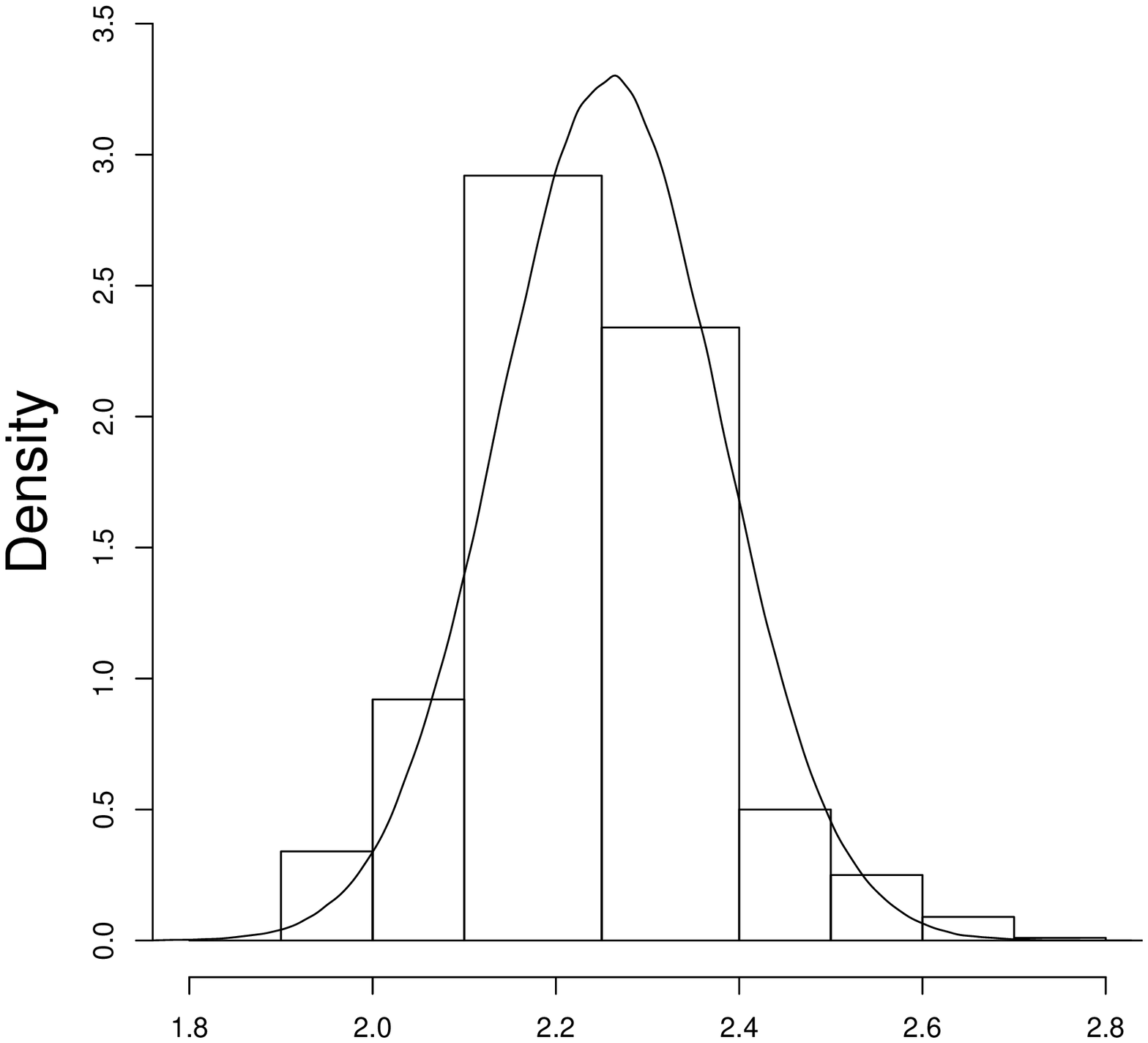} } }
\rotatebox{-90}{ \resizebox{2. in}{!}{ \includegraphics{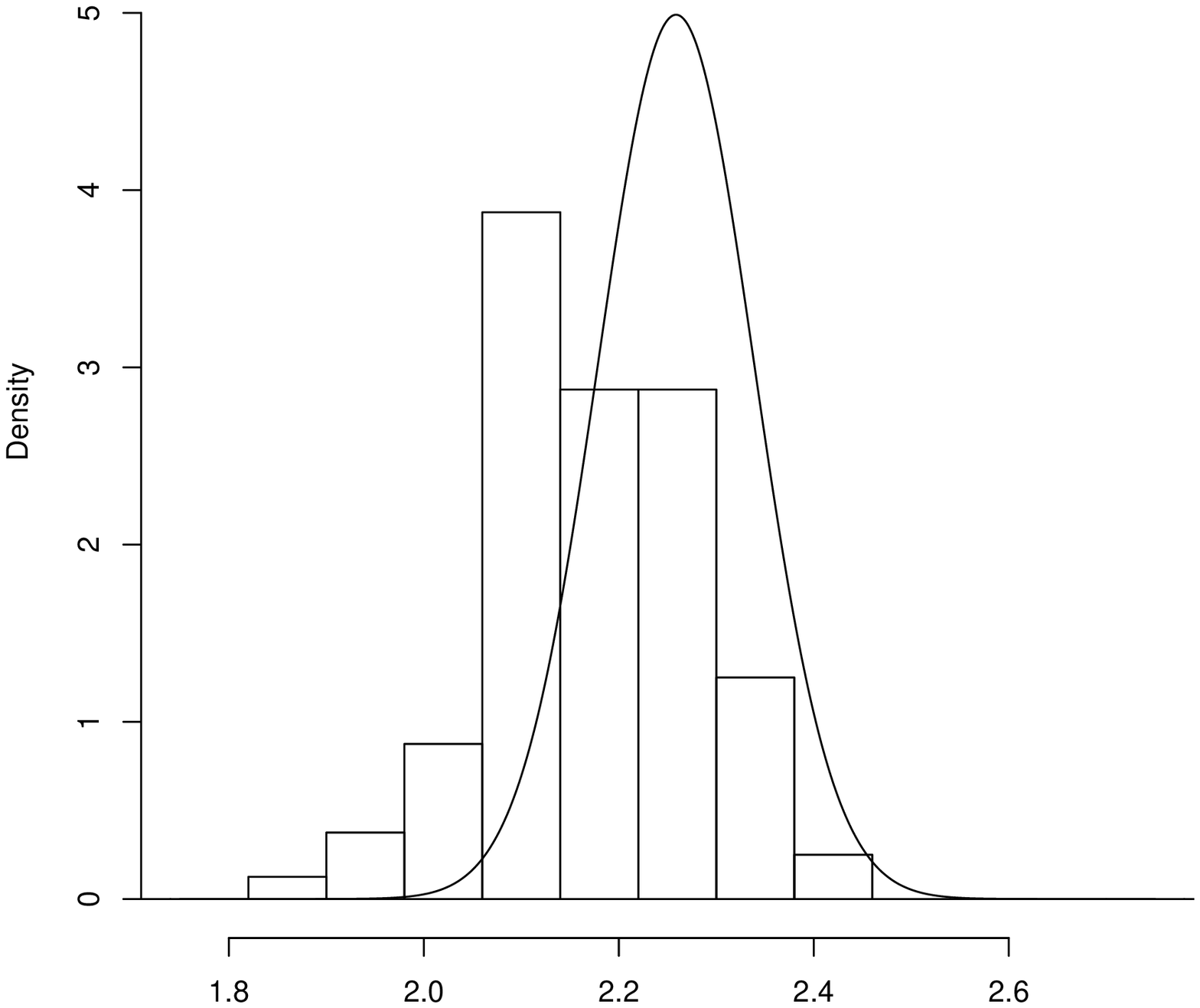} } }
\rotatebox{-90}{ \resizebox{2. in}{!}{ \includegraphics{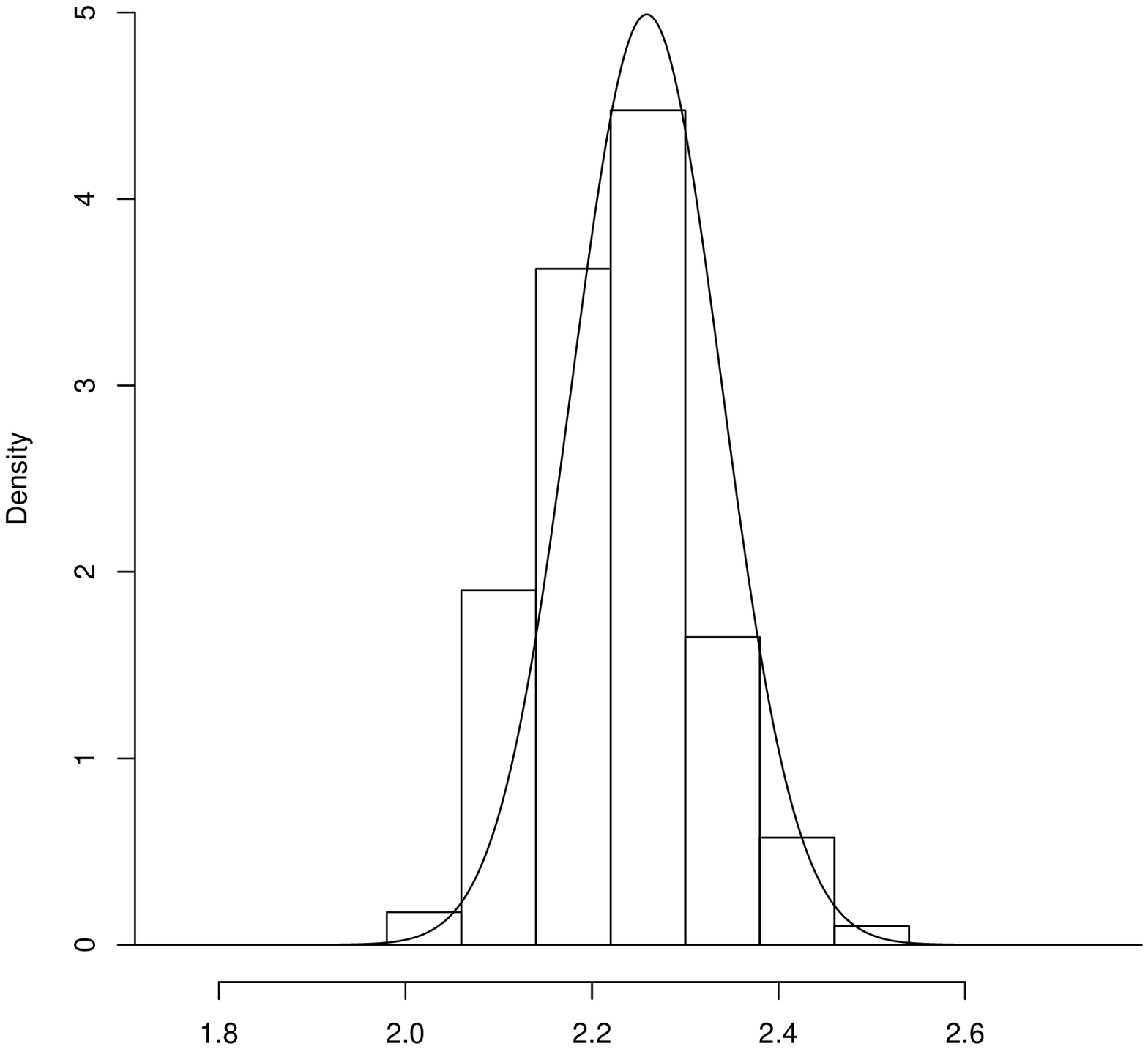} } }
\rotatebox{-90}{ \resizebox{2. in}{!}{ \includegraphics{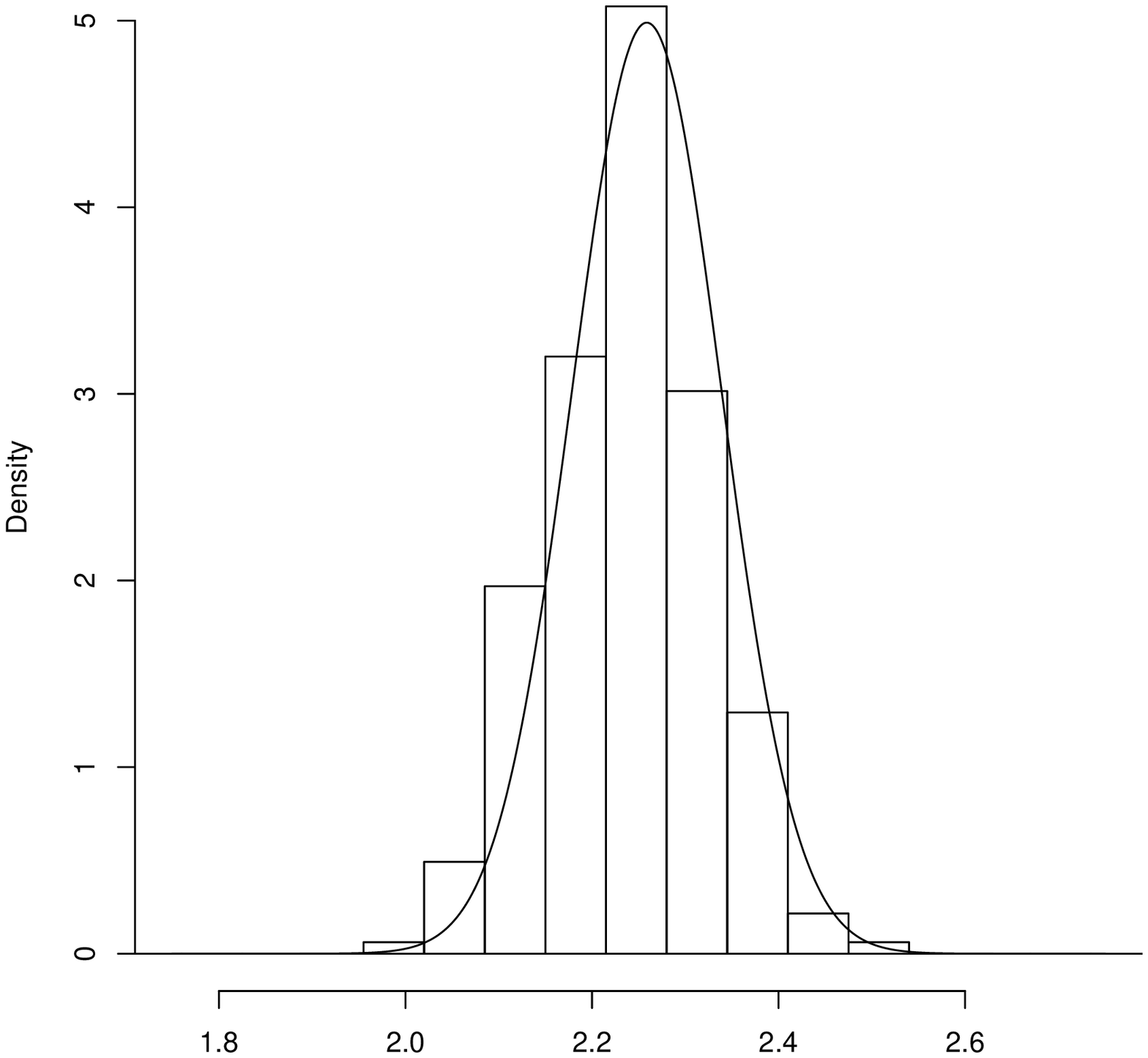} } }
\caption{
\label{fig:NormSkewMC}
Depicted in the top row are $\overline G(r=3/2,M=M_C) \stackrel{\text{approx}}{\sim}
\mathcal{N}(\mu \approx 2.2587,\sigma^2/J_{10} \approx .1918/J_{10} )$
for $J_{10}=13$ and $n=100$ (left), $n=1000$ (middle), and $n=5000$ (right).
In the bottom row, depicted are $\overline G(r=3/2,M=M_C) \stackrel{\text{approx}}{\sim}
\mathcal{N}(\mu \approx 2.2587,\sigma^2/J_{20} \approx.1918/J_{20} )$
for $J_{20}=30$ and $n=100$ (left), $n=1000$ (middle), and $n=5000$ (right).
Histograms are based on $1000$ Monte Carlo replicates and the curves are the associated approximating normal curves.
}
\end{figure}

\begin{figure}[ht]
\centering
\psfrag{Density}{ \Huge{\bf{Density}}}
\rotatebox{-90}{ \resizebox{2. in}{!}{ \includegraphics{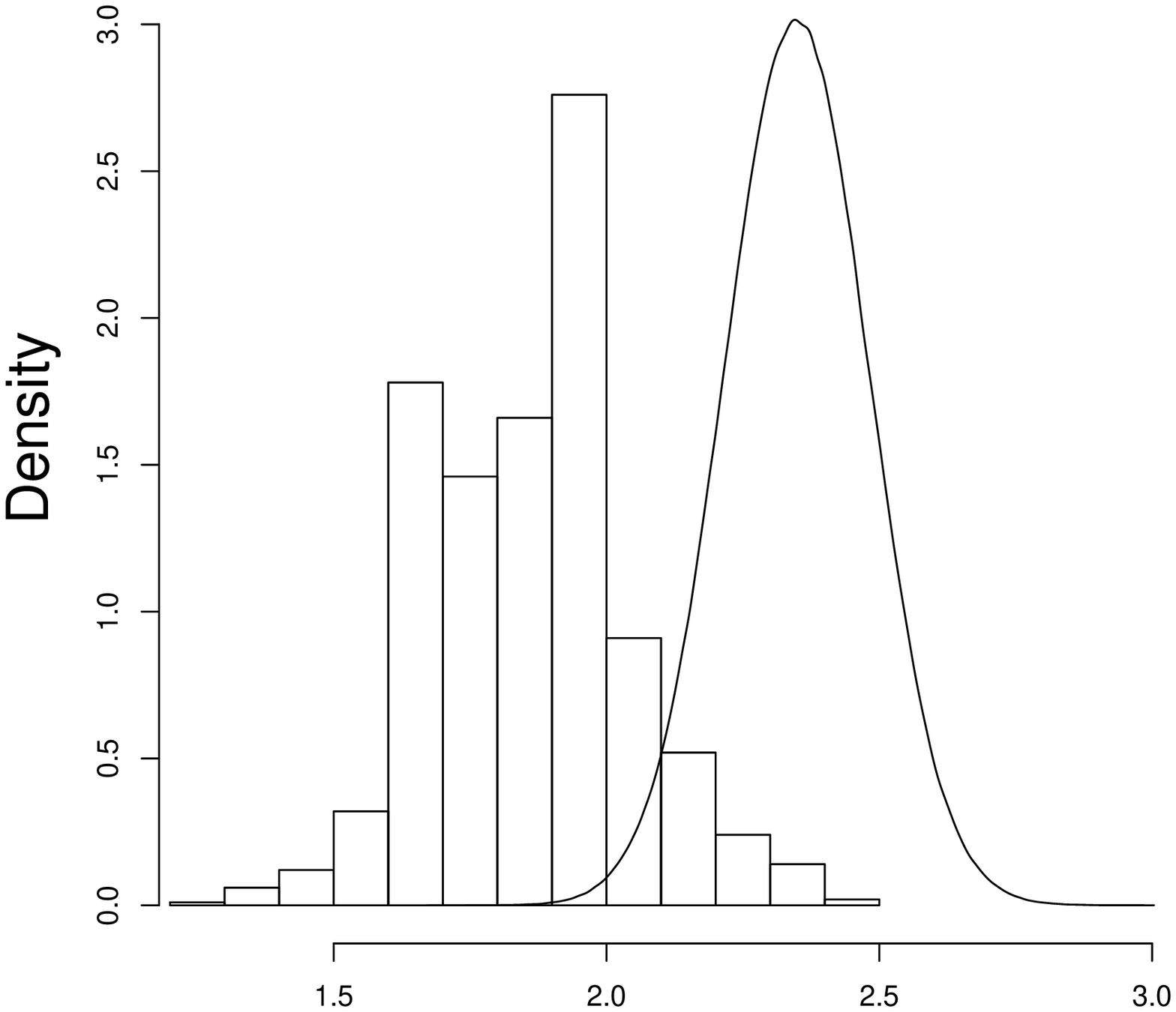} } }
\rotatebox{-90}{ \resizebox{2. in}{!}{ \includegraphics{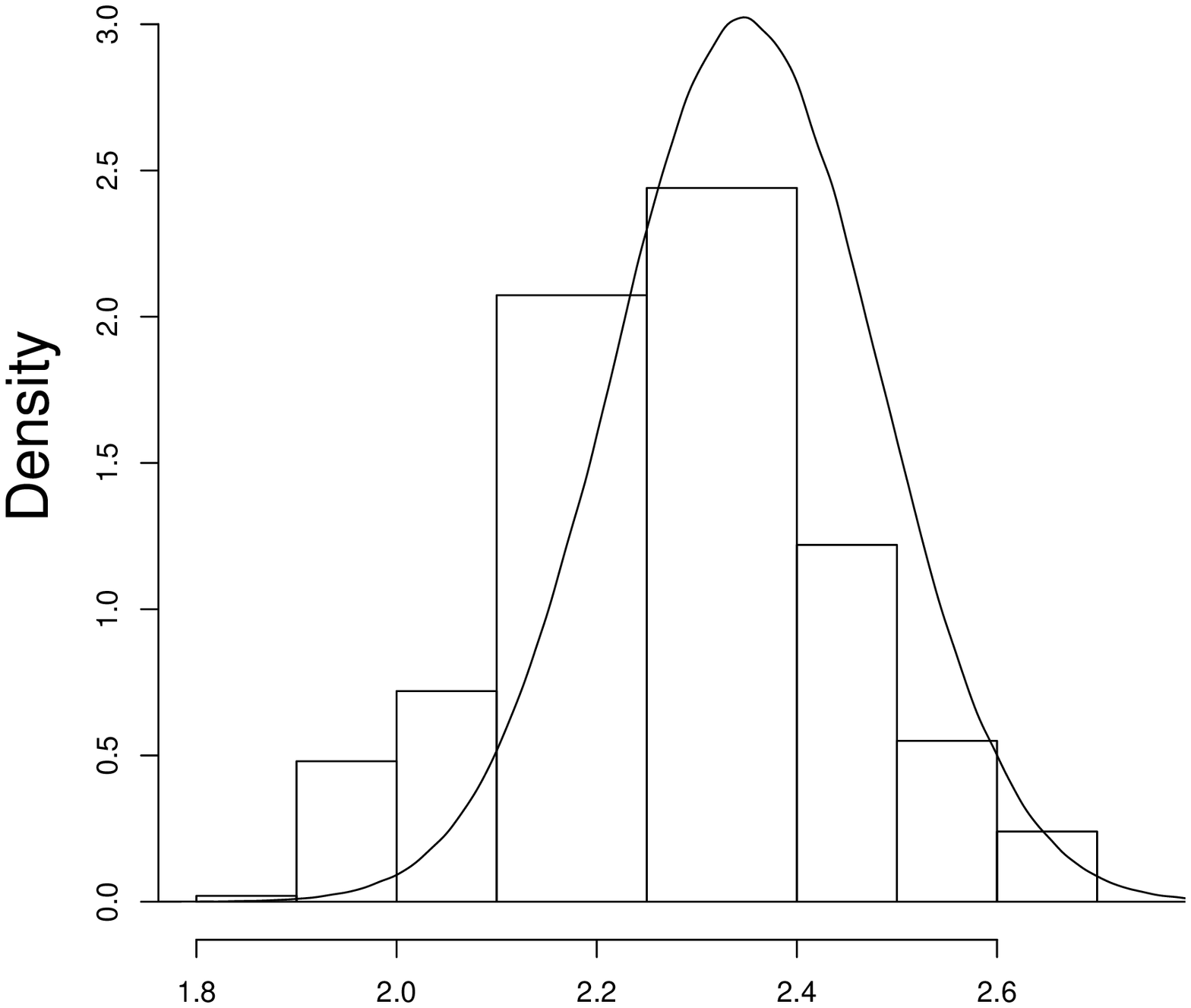} } }
\rotatebox{-90}{ \resizebox{2. in}{!}{ \includegraphics{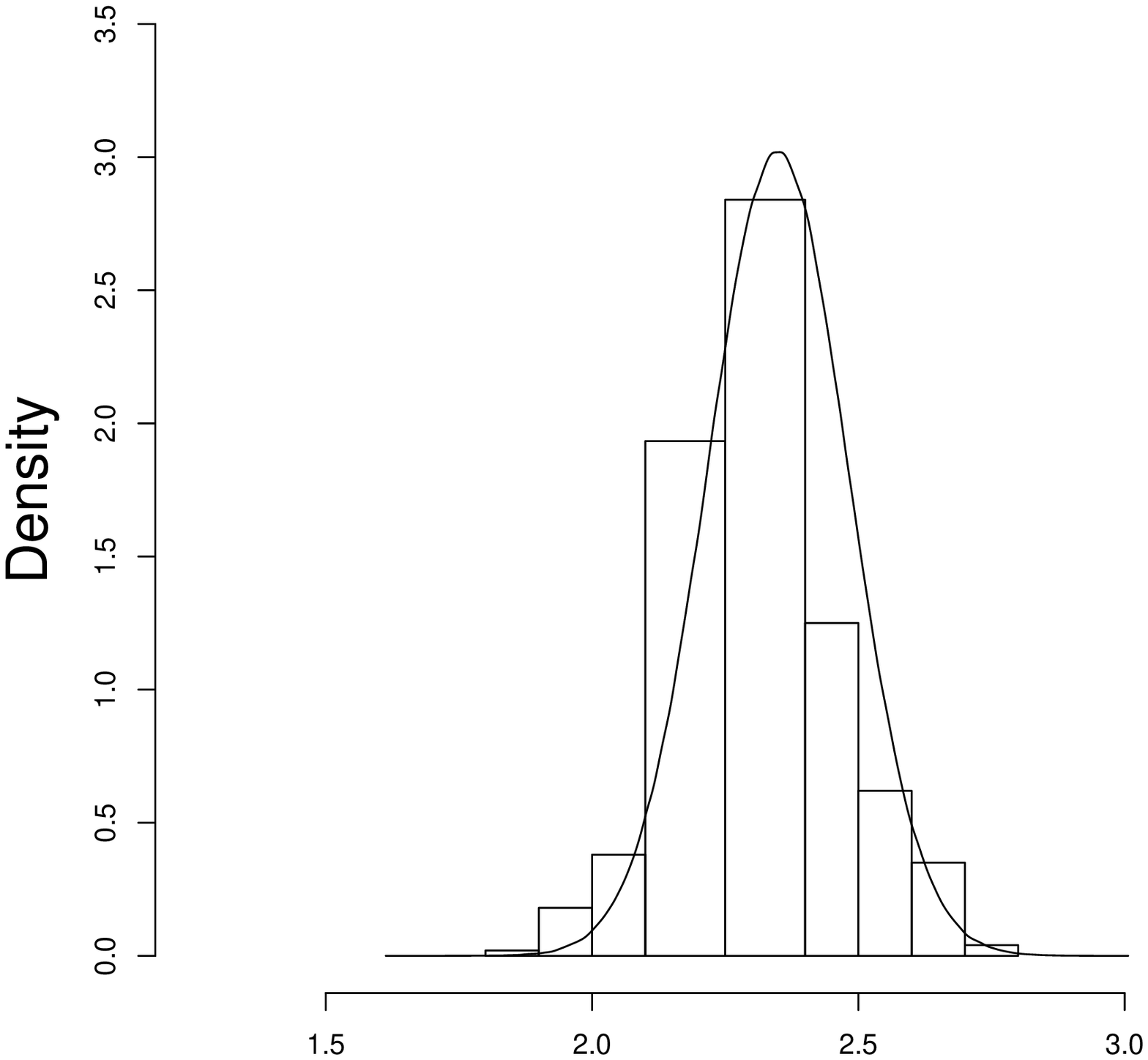} } }
\caption{
\label{fig:NormSkewMr54}
Depicted are $\overline G(r=5/4,M=\left( 7/10,\sqrt{3}/10 \right) \stackrel{\text{approx}}{\sim}
\mathcal{N}(\mu \approx 2.3486,\sigma^2/J_{10} \approx .2271/J_{10} )$
for $J_{10}=13$ and $n=100$ (left), $n=1000$ (middle), and $n=5000$ (right).
Histograms are based on $1000$ Monte Carlo replicates and the curves are the associated approximating normal curves.
}
\end{figure}

In Figure \ref{fig:NormSkewMC} (top),
we plot the histograms and the approximating normal curves for $\overline G(r,M)$
with $r=3/2$ and $M=M_C$ for $n=100, 1000$, and 5000 $\X$ points generated iid $\U(C_H(\Y_m))$
where $\Y_m$ (which yields $J_m=13$ triangles) is given in Figure \ref{fig:deltri}.
Notice that, even though the distribution looks symmetric with $n=100$,
the normal approximation is not appropriate, since not all $n_j$ are sufficiently large
to make the binomial distribution hold as in Equation \eqref{eqn:asymptotic-NYr-Jm},
but as $n$ increases (see $n=1000$ and $n=5000$ cases)
the histograms and the corresponding normal curves become more similar indicating
that the asymptotic normal approximation gets better,
since all $n_j$ are sufficiently large.
However, larger $J_m$ values require larger sample sizes in order to obtain approximate normality.
With $J_{20}=30$ triangles based on the Delaunay triangulation of $20$ $\Y$ points iid uniform on the unit square (not presented),
we plot the histograms and the approximating normal curves for $r=3/2$ and $M=M_C$ in Figure \ref{fig:NormSkewMC} (bottom).
Observe that with more triangles (i.e., as $J_m$ increases), the normal approximation gets better.
We also present the histograms of the mean domination number and the approximating
normal curves for $r=5/4$ and $M=(7/10,\sqrt{3}/10)$ in Figure \ref{fig:NormSkewMr54},
where the trend is similar to the one in Figure \ref{fig:NormSkewMC} (top).

For finite $n$, let $\overline{G}(r,M)$ be the mean domination number (per triangle)
associated with the digraph based on $\NPE^r$.
Then as a corollary to Theorem \ref{thm:stoch-order-1tri},
it follows that for $r_1<r_2$, we have $\overline{G}(r_2,M) <^{ST}\overline{G}(r_1,M)$.

\section{Alternative Patterns: Segregation and Association}
\label{sec:alternatives}
In a two class setting,
the phenomenon known as {\em segregation} occurs when members of
one class have a tendency to repel members of the other class.
For instance, it may be the case that one type of plant
does not grow well in the vicinity of another type of plant,
and vice versa.
This implies, in our notation,
that $X_i$ are unlikely to be located near any elements of $\Y$.
Alternatively, association occurs when members of one class
have a tendency to attract members of the other class,
as in symbiotic species, so that the $X_i$ will tend to
cluster around the elements of $\Y$, for example.
See, for instance, \cite{dixon:1994} and \cite{coomes:1999}.

These alternatives can be parametrized as follows:
In the one triangle case, without loss of generality
let $\Y_3=\left\{ (0,0),(1,0),(c_1,c_2) \right\}$ and $T_b=T(\Y_3)$
with $\y_1=(0,0),\y_2=(1,0)$, and $\y_3=(c_1,c_2)$.
For the basic triangle $T_b$,
let $Q_{\theta}:=\{x \in T_b: d(x,\Y_3) \le \theta\}$ for $\theta \in (0,(c_1^2+c_2^2)/2]$
and $S(F)$ be the support of $F$.
Then consider
$$\mathscr H_S:=\{F: S(F) \subseteq T_b \text{ and } P_F(X \in Q_{\theta}) < P_U(X \in Q_{\theta}) \}$$
and
$$\mathscr H_A:=\{F: S(F) \subseteq T_b \text{ and } P_F(X \in Q_{\theta}) > P_U(X \in Q_{\theta}) \}$$
where $P_F$ and $P_U$ are probabilities with respect to distribution
function $F$ and the uniform distribution on $T_b$, respectively.
So if $X_i \stackrel{iid}{\sim} F \in \mathscr H_S$,
the pattern between $\X$ and $\Y$ points is segregation,
but if $X_i \stackrel{iid}{\sim} F \in \mathscr H_A$,
the pattern between $\X$ and $\Y$ points is association.
For example the distribution family
$$\mathscr F_S:=\{F: S(F) \subset T_b \text{ and the associated pdf
$f$ increases as $d(x,\Y_3)$ increases}\}$$
is a subset of $\mathscr H_S$ and
yields samples from the segregation alternatives.
Likewise, the distribution family
$$\mathscr F_A:=\{F: S(F) \subset T_b \text{ and the associated pdf
$f$ increases as $d(x,\Y_3)$ decreases}\}$$
is a subset of $\mathscr H_A$ and
yields samples from the association alternatives.

In the basic triangle, $T_b$,
we define the $H^S_{\ve}$ and $H^A_{\ve}$ with $\ve \in \left( 0,\sqrt{3}/3 \right)$,
for segregation and association alternatives, respectively.
Under $H^S_{\ve}$, $4 \ve^2/3\times 100$ \% of the area of $T_b$
is chopped off around each vertex so that the $\X$ points are restricted
to lie in the remaining region.
That is, for $\y_j \in \Y_3$,
let $e_j$ denote the edge of $T_b$ opposite vertex $\y_j$ for $j=1,2,3$,
and for $x \in T_b$
let $\ell_j(x)$ denote the line parallel to $e_j$ through $x$.
Then define
$T_j(\ve) = \{x \in T_b: d(\y_j,\ell_j(x)) \le \ve_j\}$ where
$\displaystyle \ve_1=\frac{2c_2\ve}{3\sqrt{c_2^2+(1-c_1)^2}}$,
$\displaystyle \ve_2=\frac{2c_2\ve}{3\sqrt{c_1^2+c_2^2}}$, and
$\displaystyle \ve_3=\frac{2c_2\ve}{3}$.
Let $\mathcal T_{\ve}:=\bigcup_{j=1}^{3} T_j(\ve)$.
Then under $H^S_{\ve}$ we have
$X_i \stackrel{iid}{\sim} \U\left(T_b \setminus \mathcal T_{\ve}\right)$.
Similarly under $H^A_{\ve}$ we have
$X_i \stackrel{iid}{\sim} \U\left(\mathcal T_{\sqrt{3}/3 - \ve} \right)$.
Thus the segregation model excludes the possibility of
any $X_i$ occurring around a $\y_j$,
and the association model requires
that all $X_i$ occur around $\y_j$'s.
The $\sqrt{3}/3 - \ve$ is used in the definition of the
association alternative so that $\ve=0$
yields $H_o$ under both classes of alternatives.
Thus, we have the below distribution families under this parametrization.
\begin{equation}
\label{eqn:epsilon-alternatives}
\mathscr U_{\ve}^S:=\{F: F = \U(T_b \setminus \mathcal T_{\ve}) \}
\text{ ~and~  }
\mathscr U_{\ve}^A:=\{F: F = \U(\mathcal T_{\sqrt{3}/3 - \ve}) \}.
\end{equation}
Clearly $\mathscr U_{\ve}^S \subsetneq \mathscr H_S$
and $\mathscr U_{\sqrt{3}/3 - \ve}^A \subsetneq \mathscr H_A$,
but
$\mathscr U_{\ve}^S \nsubseteq \mathscr F_S$
and $\mathscr U_{\sqrt{3}/3 - \ve}^A \nsubseteq \mathscr F_A$.

These alternatives $H^S_{\ve}$ and $H^A_{\ve}$ with $\ve \in \left( 0,\sqrt{3}/3 \right)$,
can be transformed into the equilateral triangle as in
(\cite{ceyhan:arc-density-PE} and \cite{ceyhan:arc-density-CS}).

\begin{figure} []
\centering
 \scalebox{.35}{\input{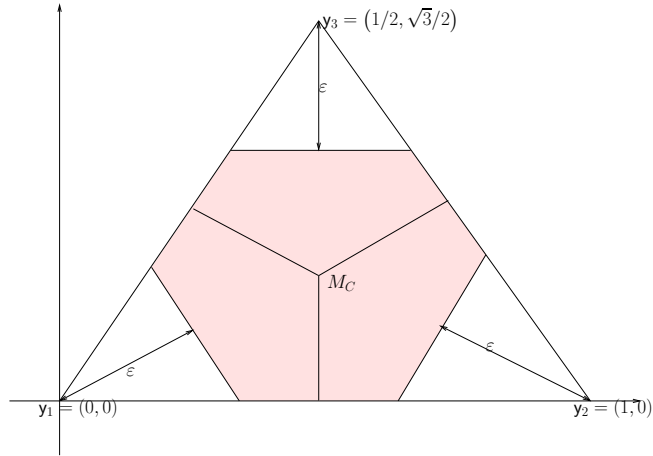}}
 \caption{An example for the segregation alternative for a particular $\varepsilon$ (shaded region),
 and its complement is for the association alternative (unshaded region) on the standard
 equilateral triangle.}
 \label{fig:seg-alt}
\end{figure}

For the standard equilateral triangle,
in $T_j(\ve) = \{x \in T_e: d(\y,\ell_j(x)) \le \ve_j\}$
we have $\ve_1=\ve_2=\ve_3=\ve$.
Thus $H^S_{\ve}$ implies
$X_i \stackrel{iid}{\sim} \U\left(T_e \setminus \mathcal T_{\ve} \right)$
and $H^A_{\ve}$ be the model under which
$X_i \stackrel{iid}{\sim} \U\left(\mathcal T_{\sqrt{3}/3 - \ve}\right)$.
See Figure \ref{fig:seg-alt} for a depiction of the above
segregation and the association alternatives in $T_e$.

\begin{remark}
These definitions of the alternatives $H^S_{\ve}$ and $H^A_{\ve}$
are given for the standard equilateral triangle.
The geometry invariance result of Theorem \ref{thm:geo-inv-NYr}
still holds under the alternatives $H^S_{\ve}$ and $H^A_{\ve}$.
In particular, the segregation alternative with $\ve \in \left( 0,\sqrt{3}/4 \right)$ in the standard equilateral triangle
corresponds to the case that in an arbitrary triangle, $\delta \cdot 100\%$ of the area is carved
away as forbidden from the vertices using line segments parallel to the opposite edge
where $\delta = 4\ve^2$ (which implies $\delta \in (0,3/4)$).
But the segregation alternative with $\ve \in \left( \sqrt{3}/4,\sqrt{3}/3 \right)$ in the standard equilateral triangle
corresponds to the case that in an arbitrary triangle, $\delta \cdot 100\%$ of the area is carved
away as forbidden from each vertex using line segments parallel to the opposite edge
where $\delta = 1-4 \left(1-\sqrt{3}\ve \right)^2$ (which implies $\delta \in (3/4,1)$).
This argument is for the segregation alternative;
a similar construction is available for the association alternative. $\square$
\end{remark}

\subsection{Asymptotic Distribution under the Alternatives}
\label{sec:asy-dist-alt}

Let $\g^S_n(F,r,M)$, $F \in \mathscr H_\theta^S$
be the domination number under segregation.
Under this alternative with $M=M_C$,
the domination number will have a discrete distribution as
$p^F_j:=P(\g_n=j)$ for $j=1,2,3$ and $p^F_1+p^F_2+p^F_3=1$.
Clearly $p^F_j$ values depend on the distribution $F$
and their explicit forms for finite $n$ or in the asymptotics
are not always analytically tractable.
The same holds for the domination number under association
$\g^A_n(F,r,M)$, $F \in \mathscr H_\theta^A$.

However, under the alternatives $H^S_{\ve}$ and $H^A_{\ve}$,
the asymptotic distribution of the domination number is much
easier to find.
Let $\g^S_n(\ve,r,M)$ and $\g^A_n(\ve,r,M)$ be the
domination numbers under segregation
and association alternatives, respectively.
Under $H^S_{\ve}$ with $M=M_C$,
the distribution of the domination number is nondegenerate
when $r=3/2-\ve\sqrt{3}/2$ which implies $r \in (9/8,3/2)$
for $\ve \in \left( 0,\sqrt{3}/4 \right)$,
and $r \in (1,9/8)$ for $\ve \in \left( \sqrt{3}/4,\sqrt{3}/3 \right)$.
In particular, the asymptotic distribution of the domination number for uniform data
in one triangle is as follows.
As $n \rightarrow \infty$, under $H^S_{\ve}$ with $M=M_C$ and $\ve \in \left( 0,\sqrt{3}/4 \right)$,
\begin{equation}
\label{eqn:asymptotic-NYr-seg1}
\g^S_n(\ve,r,M_C) \stackrel{\mathcal L}{\longrightarrow}
\left\lbrace \begin{array}{ll}
       2+\BER\left(1-p^S_{r,\ve}\right)& \text{for $r = 3/2 - \ve \sqrt{3}/2$,}\\
       1                  & \text{for $r>3/2$,}\\
       2                  & \text{for $3/2 - \ve \sqrt{3}/2< r < 3/2$,}\\
       3                  & \text{for $9/8< r < 3/2 - \ve \sqrt{3}/2$,}\\
\end{array} \right.
\end{equation}
where $p^S_{r,\ve}$ can be calculated similarly as in \eqref{eqn:p_r-form} for fixed numeric $\ve$.

Furthermore, as $n \rightarrow \infty$, under $H^S_{\ve}$ with $M=M_C$ and $\ve \in \left( \sqrt{3}/4,\sqrt{3}/3 \right)$),
\begin{equation}
\label{eqn:asymptotic-NYr-seg2}
\g^S_n(\ve,r,M_C) \stackrel{\mathcal L}{\longrightarrow}
\left\lbrace \begin{array}{ll}
       2+\BER\left(1-p^S_{r,\ve}\right) & \text{for $r = 3/2 - \ve \sqrt{3}/2$,}\\
       1                     & \text{for $r> 2- \sqrt{3}\ve$,}\\
       2                     & \text{for $3/2 - \ve \sqrt{3}/2< r < 2- \sqrt{3}\ve$,}\\
       3                     & \text{for $1< r < 3/2 - \ve \sqrt{3}/2$.}\\
\end{array} \right.
\end{equation}

Under $H^A_{\ve}$ with $M=M_C$,
the domination number $\g_n$ is nondegenerate when $r=\sqrt{3}/(2\,\ve)$ which implies $r>2$
for $\ve \in \left( 0,\sqrt{3}/4 \right)$, and $\ve \in (3/2,2)$ for $\ve \in \left( \sqrt{3}/4,\sqrt{3}/3 \right)$.
In particular, the asymptotic distribution of the domination number for uniform data
in one triangle is as follows.
As $n \rightarrow \infty$, under $H^A_{\ve}$ with $M=M_C$ and $\ve \in \left( 0,\sqrt{3}/3 \right)$,
\begin{equation}
\label{eqn:asymptotic-NYr-assoc}
\g^A_n(\ve,r,M_C) \stackrel{\mathcal L}{\longrightarrow}
\left\lbrace \begin{array}{ll}
       2+\BER\left(1-p^A_{r,\ve}\right) & \text{for $r = \sqrt{3}/(2\,\ve)$,}\\
       1                     & \text{for $r > \sqrt{3}/(2\,\ve)$,}\\
       3                     & \text{for $ r < \sqrt{3}/(2\,\ve)$,}\\
\end{array} \right.
\end{equation}
where $p^A_{r,\ve}$ can be calculated similarly as in \eqref{eqn:p_r-form} for fixed numeric $\ve$.
However, for finite $n$,
$\g^A_n(\ve,r,M_C)$ is also nondegenerate for $\sqrt{3}/(2\,\ve)-1 < r< \sqrt{3}/(2\,\ve)$.

Under segregation with general $M$,
suppose $M \in T_e \setminus \bigcup_{\y \in \Y_e} T(\y,\ve)$
(i.e., $M$ is in the support of $\X$ points under $H^S_{\ve}$).
Then for fixed $r=r_o$ for which $\g_n$ is nondegenerate under CSR
(i.e., $r_o$ is a value such that $M \in \{t_1(r_o),t_2(r_o),t_3(r_o)\}$),
then $\g_n$ is nondegenerate under $H^S_{\ve}$  if $r=r_o \left( 2-4/\left( \sqrt{3}\ve \right) \right)$.
For $r_o \in (4/3,3/2)$,
if $M \not\in T_e \setminus \bigcup_{\y \in \Y_e} T(\y,\ve)$ and
$\displaystyle \ve > \frac{3}{2}\left(1-\frac{1}{2r}\right)$, then $\g_n \rightarrow 1$ in probability as $n \rightarrow \infty$;
and the same also holds if $\displaystyle \sqrt{3}\left(1-\frac{1}{r}\right) < \ve < \frac{3}{2}\left(1-\frac{1}{2r}\right)$.
$\g_n$ is nondegenerate when $r=r_o\left( 2-4\,\ve/\sqrt{3} \right)$.
For general $M$, if $\ve \in \left( 0,\sqrt{3}/4 \right)$,
then $\g_n$ is nondegenerate when $r=r_o\left( 1-\ve/\sqrt{3} \right)$.

Under association with general $M$,
when $M \not\in \bigcup_{\y \in \Y_e} T(\y,\ve)$ then $\g_n$ is nondegenerate when $r=r_o$
(i.e., $M$ is not in the support of $\X$ points under $H^A_{\ve}$).
If $M \in \bigcup_{\y \in \Y_e} T(\y,\ve)$ then $\g_n$ is nondegenerate
when $\displaystyle r=\frac{\sqrt{3}(r_o-2)}{2\ve(r_o-1)+\sqrt{3}(2 r_o-3)}$.

\begin{theorem}
\label{thm:stoch-order-seg}
\textbf{(Stochastic Ordering)}
Let $\g^S_{n}(\ve,r,M)$ be the domination number under the segregation alternative with $\ve>0$.
Then with $\ve_j \in \left( 0,\sqrt{3}/3 \right)$, $j=1,2$, $\ve_1 > \ve_2$
implies that $\g^S_{n}(\ve_1,r,M) <^{ST} \g^S_{n}(\ve_2,r,M)$.
\end{theorem}

\noindent
{\bf Proof:}
Note that for $\ve_1 > \ve_2$ and finite $n$,
$P(\g^S_{n}(\ve_1,r,M)=1)>P(\g^S_{n}(\ve_2,r,M)=1)$
and $P(\g^S_{n}(\ve_1,r,M)=2)>P(\g^S_{n}(\ve_2,r,M)=2)$,
hence the desired result follows.
$\blacksquare$

Note that for Theorem \ref{thm:stoch-order-seg} to hold in the limiting case when
$r \in [1,3/2]$ and $M \in \{t_1(r),t_2(r),t_3(r)\}$,
$\ve_1 \in I_i(r)$ and $\ve_2 \in I_j(r)$ should hold for $i<j$ where
$I_1(r)=\left( (2-r)/\sqrt{3},\sqrt{3}/3 \right)$,
$I_2(r)=\left( (3-2r)/\sqrt{3},(2-r)/\sqrt{3} \right)$, and
$I_3(r)=\left( 0,(3-2r)/\sqrt{3} \right)$.
For $\ve \in \big(0,\sqrt{3}/4\big]$,
$\g^S_{n}(\ve,r,M) \rightarrow 2$ in probability as $n \rightarrow \infty$,
and for $\ve \in \left( \sqrt{3}/4,\sqrt{3}/3 \right)$,
$\g^S_{n}(\ve,r,M) \rightarrow 1$ in probability as $n \rightarrow \infty$.

Similarly, the stochastic ordering result of Theorem \ref{thm:stoch-order-seg}
holds for association for all $\ve$ and $n<\infty$,
with the inequalities being reversed.

Notice that under segregation with $\ve \in \left( 0,\sqrt{3}/4 \right)$,
$\g_{n,\ve}(r,M_C)$ is degenerate in the limit except for $\displaystyle r=(3-\sqrt{3}\,\ve)/2 $.
With $\ve \in \left( \sqrt{3}/4,\sqrt{3}/3 \right)$,
$\g_{n,\ve}(r,M_C)$ is degenerate in the limit
except for $r=3-\ve/\sqrt{3}$.
Under association with $\ve \in \left( 0,\sqrt{3}/4 \right)$,
$\g_{n,\ve}(r,M_C)$ is degenerate in the limit
except for $r=\sqrt{3}/(2\,\ve)$.

\begin{remark}
\label{rem:alt-multi-tri}
\textbf{The Alternatives with Multiple Triangles:}
In the multiple triangle case,
the segregation and association alternatives,
$H^S_{\ve}$ and $H^A_{\ve}$ with $\ve \in \left( 0,\sqrt{3}/3 \right)$,
are defined as in the one-triangle case,
in the sense that,
when each triangle (together with the data in it)
is transformed to the standard equilateral triangle as in Theorem \ref{thm:geo-inv-NYr},
we obtain the same alternative pattern described above.

Let $\g^S_{n,m}(\ve,r,M)$ and $\g^A_{n,m}(\ve,r,M)$ be the domination numbers under segregation
and association alternatives in the multiple triangle case with $m$ triangles, respectively.
The extensions of their distributions from Equations \eqref{eqn:asymptotic-NYr-seg1},
\eqref{eqn:asymptotic-NYr-seg2}, and \eqref{eqn:asymptotic-NYr-assoc}
are similar to the extension of the distribution of the domination number
from one-triangle to multiple-triangle case under the null hypothesis in Section \ref{sec:NYr-mult-tri}.
Furthermore,
the stochastic ordering result of Theorem \ref{thm:stoch-order-seg}
extends in a straightforward manner.
$\square$
\end{remark}

\subsection{The Test Statistics and Their Distributions}
\label{sec:test-stat-distr}
A translated form of the domination number of the PCD
is a test statistic for the segregation/association alternative:
\begin{equation}
\label{eqn:Bnm-test-stat}
B_{n,m}:=
\left\lbrace \begin{array}{ll}
\g_n(r,M)-2J_m=\sum_{j=1}^{J_m} \g_{n_j}(r,M)-2J_m & \text{if  $\g_n(r,M)>2J_m$,}\\
       0            & \text{otherwise.}\\
\end{array} \right.
\end{equation}
Rejecting for extreme values of $B_{n,m}$ is appropriate,
since under segregation we expect $B_{n,m}$ to be small,
while under association we expect $B_{n,m}$ to be large.
Using this test statistic the critical value for finite $J_m$
and large $n$ for the one-sided level $\alpha$ test against segregation
is given by $b_{\alpha}$, the $(\alpha)100$th percentile of  $BIN(J_m,1-p_r)$
(i.e., the test rejects for $B_{n,m} \leq b_{\alpha}$),
and against association, the test rejects for $B_{n,m} \geq b_{1-\alpha}$.

Similarly, the mean domination number (per triangle) of the PCD,
$\displaystyle \overline G (r,M):=\frac{1}{J_m}\,\sum_{j=1}^{J_m} \g_{n_j}(r,M)$,
can also be used as a test statistic for the segregation/association alternative
when $n \gg J_m$ and both $n$ and $J_m$ are sufficiently large.
Rejecting for extreme values of $\overline G(r,M)$ is appropriate,
since under segregation we expect $\overline G(r,M)$ to be small,
while under association we expect $\overline G(r,M)$ to be large.
Using the standardized test statistic
\begin{equation}
\label{eqn:Snm-test-stat}
S_{n,m} = \sqrt{J_m} (\overline G(r,M) - \mu)/\sigma,
\end{equation}
where $\mu=3-p_r$ and $\sigma^2=p_r(1-p_r)/J_m$,
the asymptotic critical value
for the one-sided level $\alpha$ test against segregation
is given by $z_{\alpha} = \Phi^{-1}(\alpha)$
where $\Phi(\cdot)$ is the standard normal distribution function.
The test rejects for $S_{n,m}<z_{\alpha}$.
Against association,
the test rejects for $S_{n,m}>z_{1-\alpha}$.

\begin{figure}[ht]
\centering
\rotatebox{-90}{ \resizebox{2. in}{!}{\includegraphics{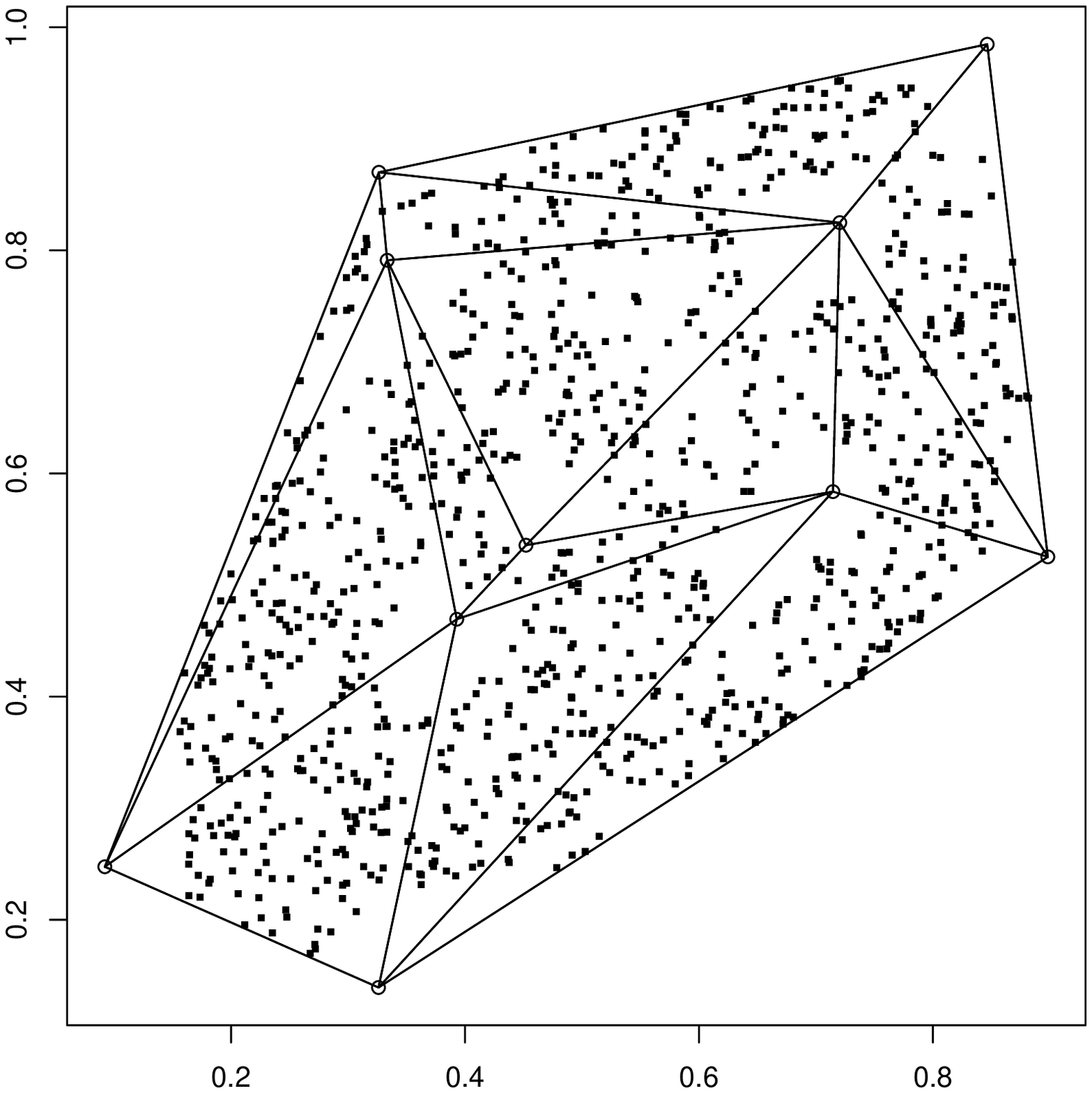}}}
\rotatebox{-90}{ \resizebox{2. in}{!}{\includegraphics{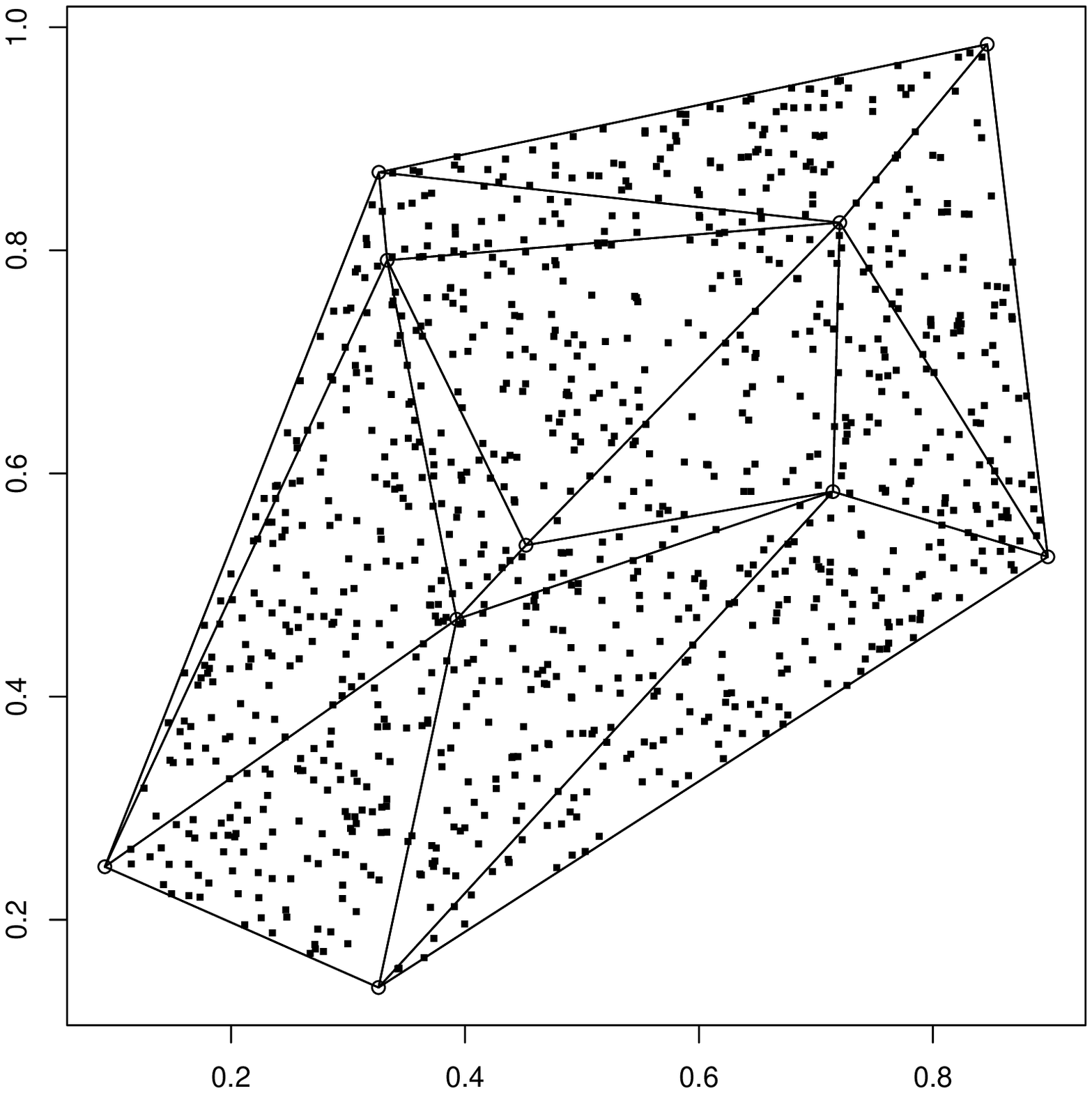}}}
\rotatebox{-90}{ \resizebox{2. in}{!}{\includegraphics{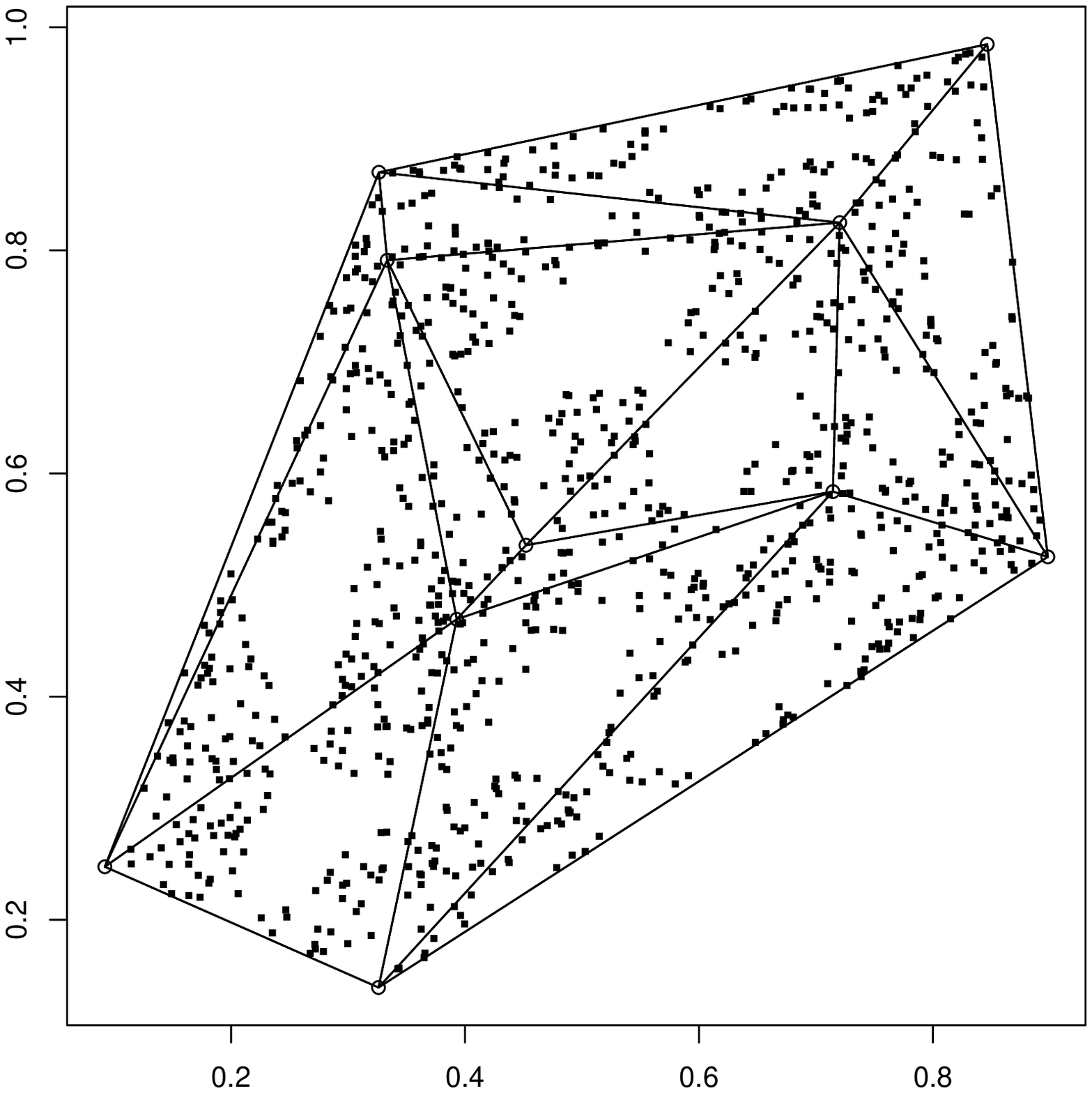}}}
\caption{
\label{fig:deldata}
A realization of segregation (left), CSR (middle), and association (right) for $|\Y|=10$, $J_{10}=13$, and $n=1000$.
}
\end{figure}

Depicted in Figure \ref{fig:deldata} are the segregation with $\delta=3/16$, CSR, and
association with $\delta=1/4$ realizations for $m=10$ and $J_m=13$, and $n=1000$.
The associated mean domination numbers with $r=3/2$ are $2.000,\; 2.1538$, and $3.000$,
for the segregation alternative, null realization, and the association alternatives, respectively,
yielding $p$-values $\approx 0.000$, $0.6139$, and $\approx 0.000$ based on binomial approximation,
and $p$-values $0.0166$, $0.3880$, and $< 0.0001$ based on normal approximation.
We also present a Monte Carlo power investigation in Section \ref{sec:monte-carlo-sim} for these cases.

\begin{theorem}
\label{thm:consistency-I}
\textbf{(Consistency-I)}
Let $\g^S_{n,m}(F,r,M)$ and $\g^A_{n,m}(F,r,M)$ be the domination numbers under segregation
and association alternatives in the multiple triangle case with $m$ triangles, respectively.
The test against segregation with $F \in \mathscr H_S$
which rejects for $S_{n,m} < z_{\alpha}$
and
the test against association with $F \in \mathscr H_A$
which rejects for $S_{n,m}>z_{1-\alpha}$ are consistent.
\end{theorem}

\noindent
{\bf Proof:}
Given $F \in \mathscr H_S$.
Let $\g_{n,m}(\U,r,M)$ be the domination number for $\X_n$
being a random sample from $\U(\TY)$.
Then $P(\g^S_{n,m}(F,r,M)=1) \ge P(\g_{n,m}(\U,r,M)=1)$;
$P(\g^S_{n,m}(F,r,M) \le 2) \ge P(\g^S_{n,m}(\U,r,M) \le 2)$;
and
$P(\g^S_{n,m}(F,r,M)=3) \le P(\g^S_{n,m}(\U,r,M)=3)$.
Hence $S_{n,m}<0$ with probability 1, as $n \gg m \rightarrow \infty$.
Hence consistency follows from the consistency
of tests which have asymptotic normality.
The consistency against the association alternative can be proved similarly.
$\blacksquare$

Below we provide a result which is stronger,
in the sense that it will hold for finite $m$ and $n \rightarrow \infty$.

\begin{theorem}
\label{thm:consistency-II}
\textbf{(Consistency-II)}
Let $\g^S_{n,m}(\ve,r,M)$ and $\g^A_{n,m}(\ve,r,M)$ be the domination numbers under segregation
and association alternatives $H^S_{\ve}$ and $H^A_{\ve}$ in the multiple triangle case with $m$ triangles, respectively.
Let $J^*(\alpha,\ve):=\left \lceil \Bigl(\frac{\sigma \cdot z_{\alpha}}{\overline G(r,M)-\mu} \Bigr)^2 \right \rceil$
where $\lceil \cdot \rceil$ is the ceiling function and
$\ve$-dependence is through $\overline G(r,M)$ under a given alternative.
Then the test against $H^S_{\ve}$ which rejects for
$S_{n,m} < z_{\alpha}$ is consistent for all $\ve \in \left( 0,\sqrt{3}/3 \right)$
and $J_m \ge J^*(\alpha,\ve)$,
and
the test against $H^A_{\ve}$ which rejects for $S_{n,m}>z_{1-\alpha}$
is consistent for all $\ve \in \left( 0,\sqrt{3}/3 \right)$ and $J_m \ge J^*(1-\alpha,\ve)$.
\end{theorem}

\noindent
{\bf Proof:}
Let $\ve>0$.
Under $H^S_{\ve}$, $\g^S_{n}(\ve,r,M)$ is
degenerate in the limit as $n \rightarrow \infty$,
which implies $\overline G(r,M)$ is a constant a.s.
In particular, for $\ve \in (0,\sqrt{3}/4]$,
$\overline G(r,M)=2$ and for $\ve \in \left( \sqrt{3}/4,\sqrt{3}/3 \right)$,
$\overline G(r,M) \leq 2$ a.s. as $n \rightarrow \infty$.
Then the test statistic $S_{n,m} = \sqrt{J_m} (\overline G(r,M) - \mu)/\sigma$ is a constant a.s.
and $J_m \ge J^*(\alpha,\ve)$ implies that $S_{n,m}< z_{\alpha}$ a.s.
Hence consistency follows for segregation.

Under $H^A_{\ve}$, as $n \rightarrow \infty$,
$\overline G(r,M)=3$ for all $\ve \in \left( 0,\sqrt{3}/3 \right)$, a.s.
Then $J_m \ge J^*(1-\alpha,\ve)$ implies that $S_{n,m} > z_{1-\alpha}$ a.s.,
hence consistency follows for association.
$\blacksquare$

Consistency in the sense of Theorems \ref{thm:consistency-I} and \ref{thm:consistency-II}
follows for $B_{n,m}$ similarly.

\begin{remark}
(\textbf{Asymptotic Efficiency})
Pitman asymptotic efficiency (PAE)
provides for an investigation of ``local (around $H_o$) asymptotic power''.
This involves the limit as $n \rightarrow \infty$ as well as
the limit as $\ve \rightarrow 0$.
A detailed discussion of PAE is available in \cite{kendall:1979} and \cite{eeden:1963}.
For segregation or association alternatives $H^S_{\ve}$ and $H^A_{\ve}$
the PAE is not applicable because the Pitman conditions
(\cite{eeden:1963}) are not satisfied by the test statistic, $\overline{G}(r,M)$.

Hodges-Lehmann asymptotic efficiency analysis (\cite{hodges:1956}) and
asymptotic power function analysis (\cite{kendall:1979}) are not applicable here either.
However, when $M=M_C$ (which also implies $r=3/2$),
for $\ve$ small and $n$ large enough,
this test is very sensitive for both alternatives
because $\g^S_{n}(\ve,3/2,M_C) \rightarrow 2$ in probability as $n \rightarrow \infty$ for segregation
and $\g^A_{n}(\ve,3/2,M_C) \rightarrow 3$ in probability as $n \rightarrow \infty$ for association.
That is, the test statistic becomes degenerate in the limit for all $\ve>0$
but in the right direction for both alternatives.
On the other hand, when $M\not=M_C$ (i.e., $r \not=3/2$)
this test is very sensitive for the segregation alternative
since $\g^S_n(\ve,r,M) \rightarrow 2$ in probability as $n \rightarrow \infty$;
the same holds for the association alternative,
but the test is not as sensitive as in the segregation case,
since we only have $\g^A_n(\ve,r,M)<^{ST}\g_n(r,M)$.
$\square$
\end{remark}

However, the variance of $\g_n(r,M)$ is minimized when $p_r=1/2$,
which happens when $r \approx 1.395$ (obtained numerically).
Hence, we expect the test to have higher power under the alternatives
for $r$ around 1.40.

\begin{remark}
The choice of the null pattern in Section \ref{sec:NYr-mult-tri}
and the conditions in Theorem \ref{thm:asy-normality} seem to be somewhat stringent;
i.e., $\X$ points are assumed to be uniformly distributed in the convex hull of $\Y$ points,
which might not be realistic in practice.
However, if the supports of distributions of $\X$ and $\Y$ points do not intersect,
or mildly intersect, then it is clear that the null hypothesis is violated (i.e., two classes are segregated)
which is easily detected by the test statistics $B_{n,m}$ or $S_{n,m}$
(see Equations \eqref{eqn:Bnm-test-stat} and \eqref{eqn:Snm-test-stat})
as they tend to be smaller under segregation than expected under CSR.
When their supports have non-empty intersection,
then either the $\X$ points are segregated from the $\Y$ points,
or follow CSR, or are associated with the $\Y$ points
in this intersection.
Then we only consider the $\Y$ points in this support intersection,
then our inference will be local (i.e., restricted to this intersection).
If one takes all of the $\Y$ points,
then our inference will be a global one
(i.e., for the entire support of $\Y$ points).
$\square$
\end{remark}

\section{Monte Carlo Simulation Analysis}
\label{sec:monte-carlo-sim}

\subsection{Empirical Size Analysis under CSR}
\label{sec:emp-size-CSR}
\begin{figure}[ht]
\centering
\rotatebox{-90}{ \resizebox{2. in}{!}{\includegraphics{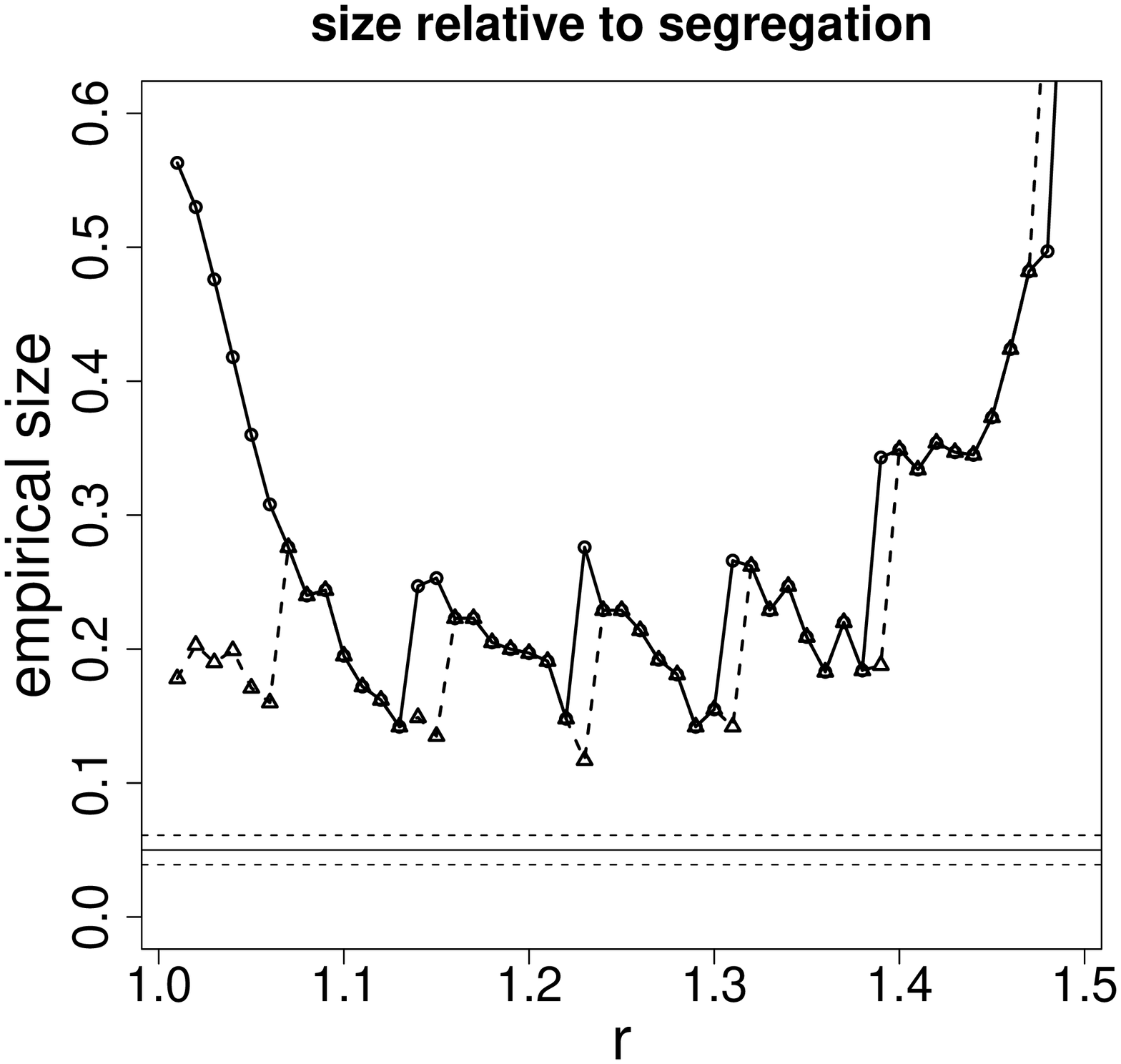} }}
\rotatebox{-90}{ \resizebox{2. in}{!}{\includegraphics{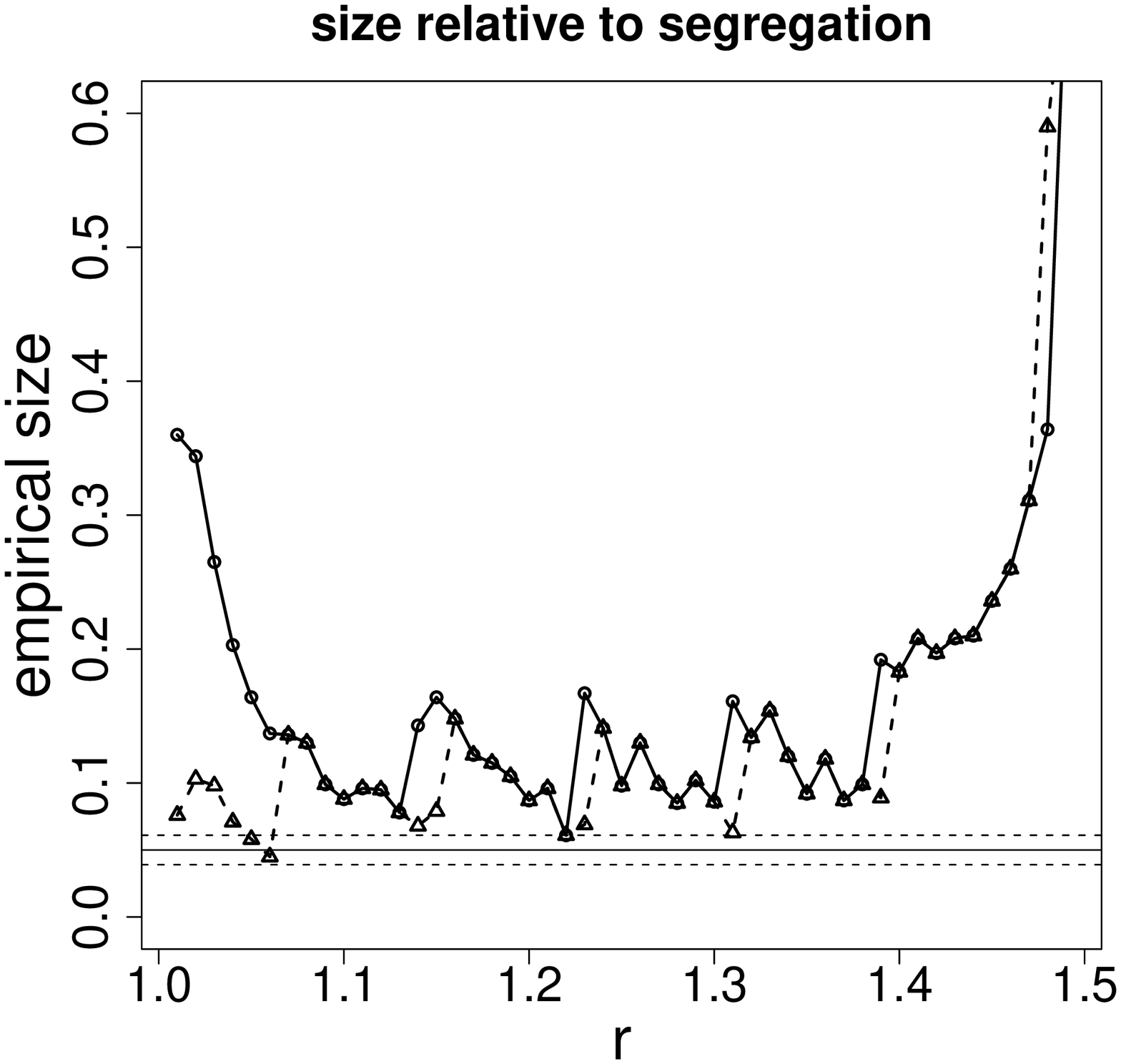} }}
\rotatebox{-90}{ \resizebox{2. in}{!}{\includegraphics{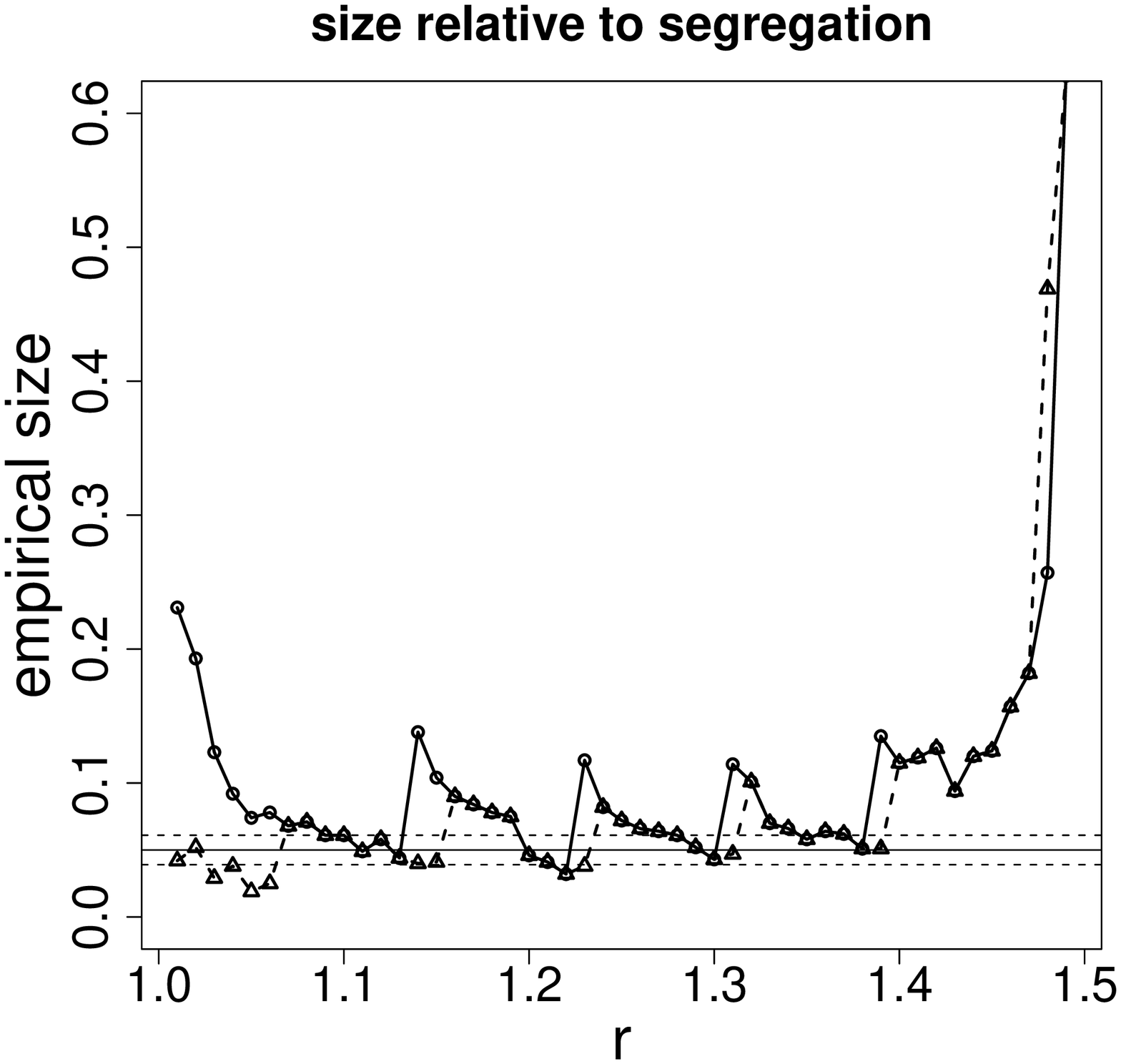} }}
\rotatebox{-90}{ \resizebox{2. in}{!}{\includegraphics{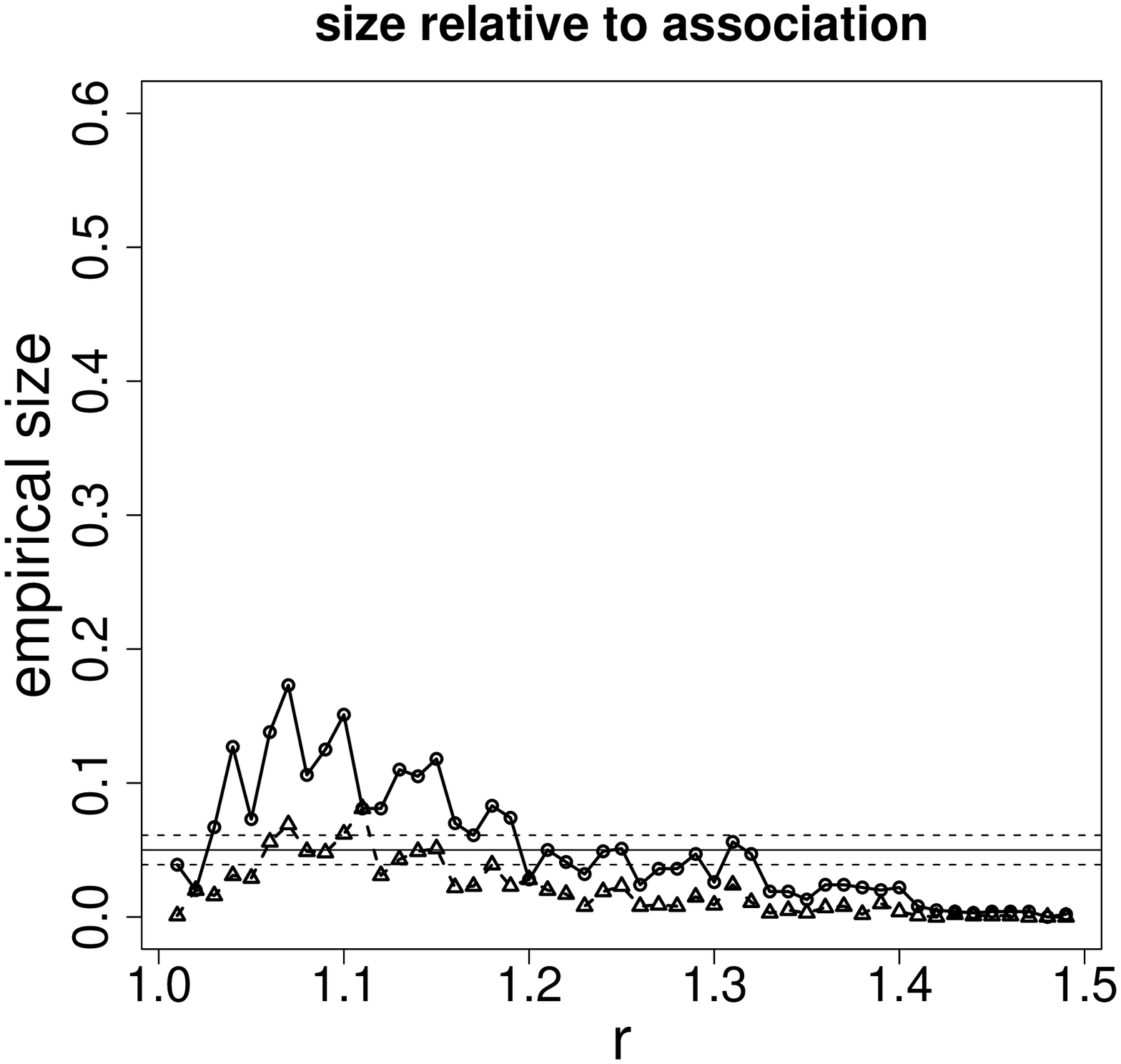} }}
\rotatebox{-90}{ \resizebox{2. in}{!}{\includegraphics{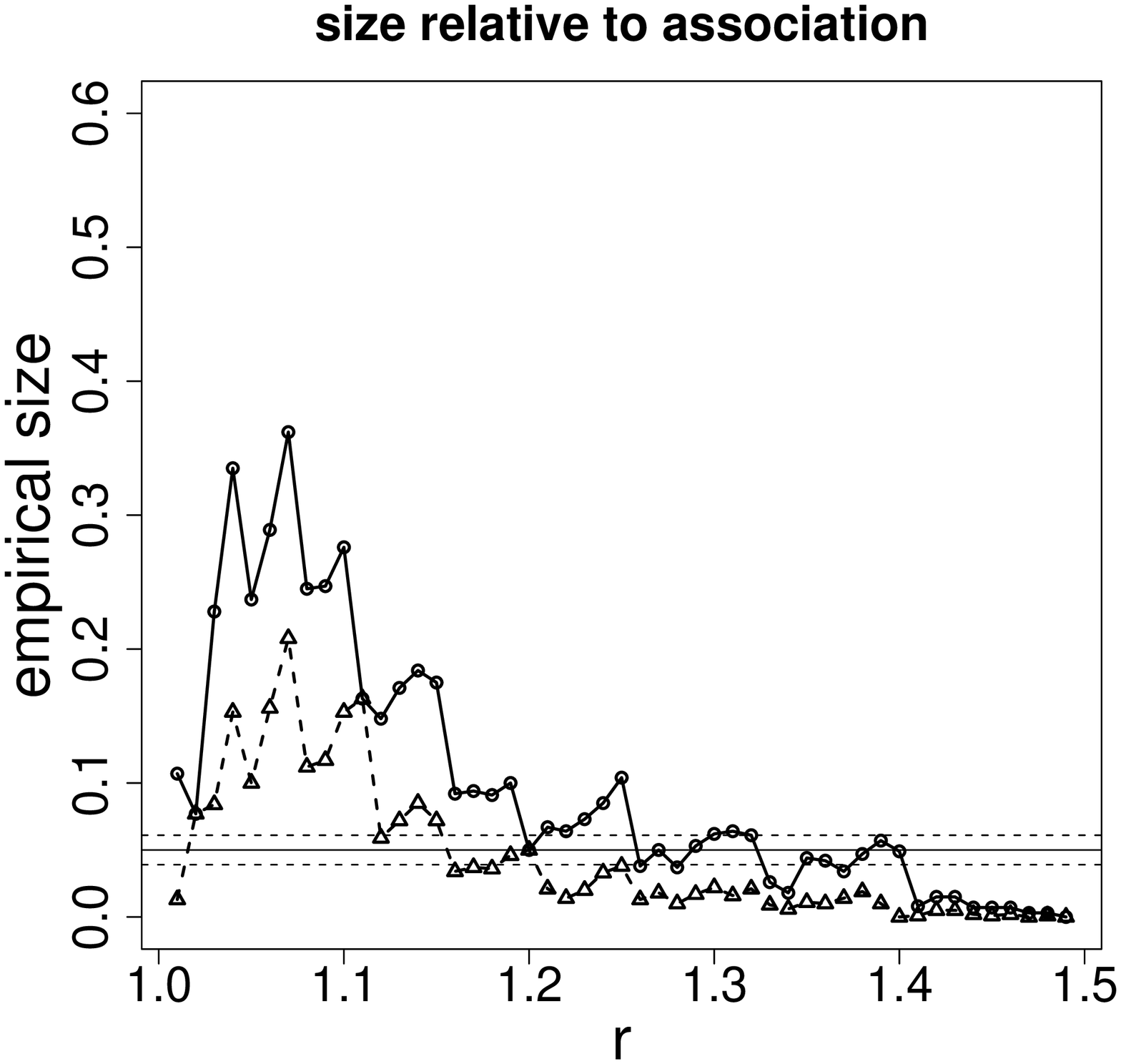} }}
\rotatebox{-90}{ \resizebox{2. in}{!}{\includegraphics{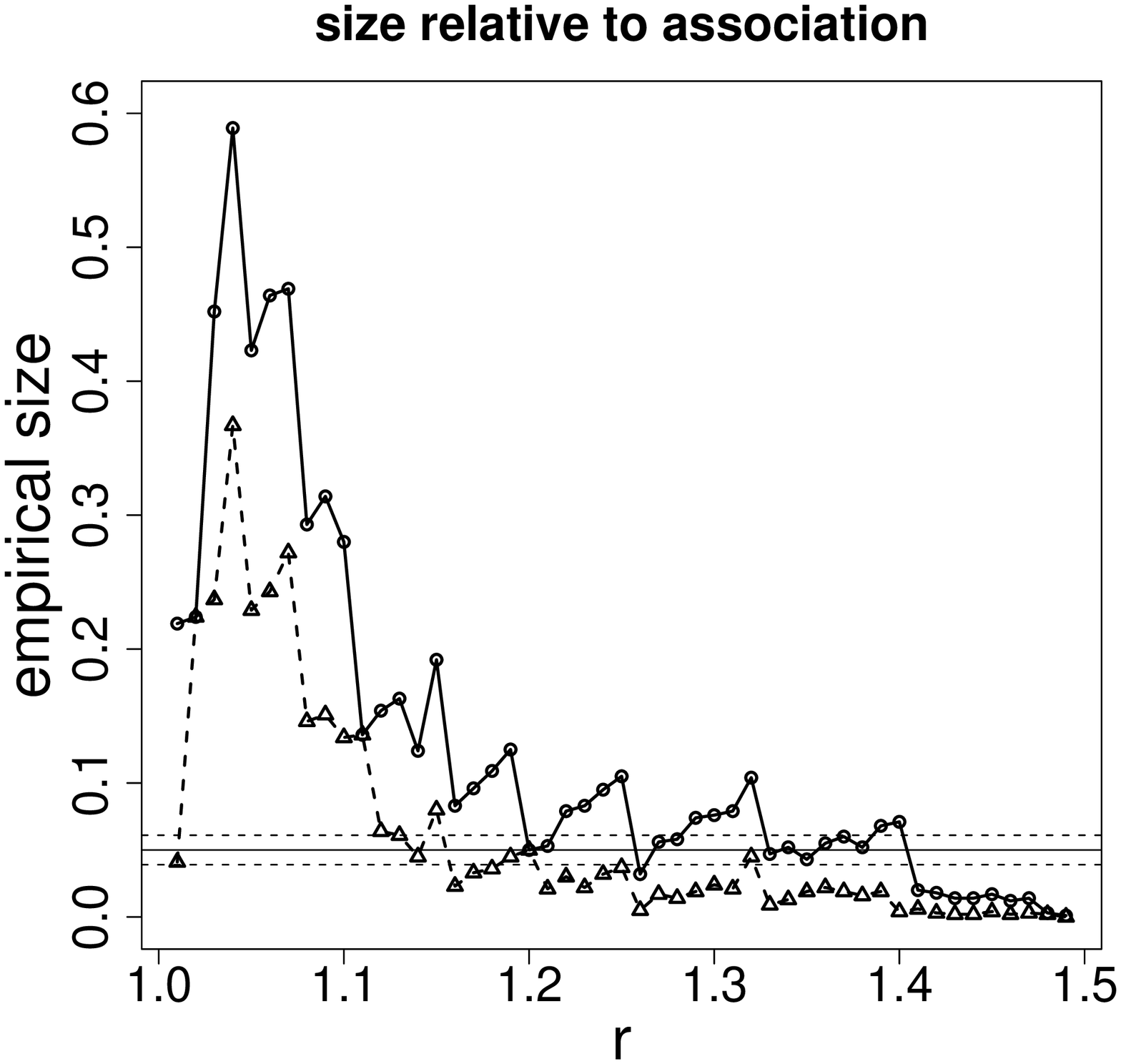} }}
\caption{
\label{fig:emp-size-CSR}
The empirical size estimates for the left-sided alternative (i.e., relative to segregation)
and the right-sided alternative (i.e., relative to association)
with $n=500$ (left), $n=1000$ (middle), and $n=2000$ (right) under the CSR pattern.
The empirical sizes based on the binomial distribution are plotted in circles ($circ$) and joined with solid lines,
and those based on the normal approximation are plotted in triangles ($\triangle$) and joined with dashed lines.
The horizontal lines are located at .039 (upper threshold for conservativeness),
.050 (nominal level), and .061 (lower threshold for liberalness).
}
\end{figure}

For the null pattern of CSR,
we generate $n$ $\X$ points iid $\U(C_H(\Y_{10}))$ where $\Y_{10}$ is the set of
the 10 $\Y$ points in Figure \ref{fig:deldata}.
We calculate and record the domination number $\g_n(r,M)$ and the mean domination number (per triangle),
$\overline G(r,M)$ for $r=1.00,1.01,1.02,\ldots,1.49$ at each Monte Carlo replicate.
We repeat the Monte Carlo procedure $N_{mc}=1000$ times for each of $n=500, 1000, 2000$.
Using the critical values based on the binomial distribution for the domination number and the normal approximation for
$\overline G(r,M)$, we calculate the empirical size estimates for both right- and left-sided tests.
The empirical sizes significantly smaller (larger) than .05 are deemed conservative (liberal).
The asymptotic normal approximation to proportions is used in determining the significance of
the deviations of the empirical sizes from .05.
For these proportion tests, we also use $\alpha=.05$ as the significance level.
With $N_{mc}=1000$, empirical sizes less than .039 are deemed conservative,
greater than .061 are deemed liberal at $\alpha=.05$ level.
The empirical sizes together with upper and lower limits of liberalness and conservativeness
are plotted in Figure \ref{fig:emp-size-CSR}.
Observe that right-sided tests are liberal with being less liberal when sample size $n$ increases,
and it has about the nominal level for most $r$ values between 1.1 and 1.4.
The left-sided test tends to be liberal for small $r$, and conservative for large $r$,
but has about the desired nominal level for $r$ around 1.2 and 1.3.

Since $p_r$ has a different form when $r=1.50$,
we estimate the empirical sizes for $r=1.50$ separately.
The size estimates for $n=500,1000,$ and $2000$ relative to segregation and
association alternatives are presented in Table \ref{tab:emp-size-power-r=3/2}.
Based on the Monte Carlo simulations under CSR,
the use of domination number for $r \in (1.45,1.50)$ is not recommended,
as the test is extremely liberal for the segregation (i.e., left-sided) alternative,
while it is extremely conservative for the association (i.e., right-sided) alternative.
This deviation from the nominal level for the test is due to the fact that
for $r \in (1.45,1.50)$ much larger sample sizes are required for the binomial
and the normal approximations to hold.
Instead of $r \in (1.45,1.50)$, we recommend the use of $r=3/2$ with the
asymptotic distribution provided in \cite{ceyhan:dom-num-NPE-SPL}.

\subsection{Empirical Power Analysis under the Alternatives}
\label{sec:monte-carlo-power}
\begin{figure}[ht]
\centering
\psfrag{kernel density estimate}{ \Huge{\bf{kernel density estimate}}}
\psfrag{relative density}{ \Huge{\bf{relative density}}}
\rotatebox{-90}{ \resizebox{2. in}{!}{ \includegraphics{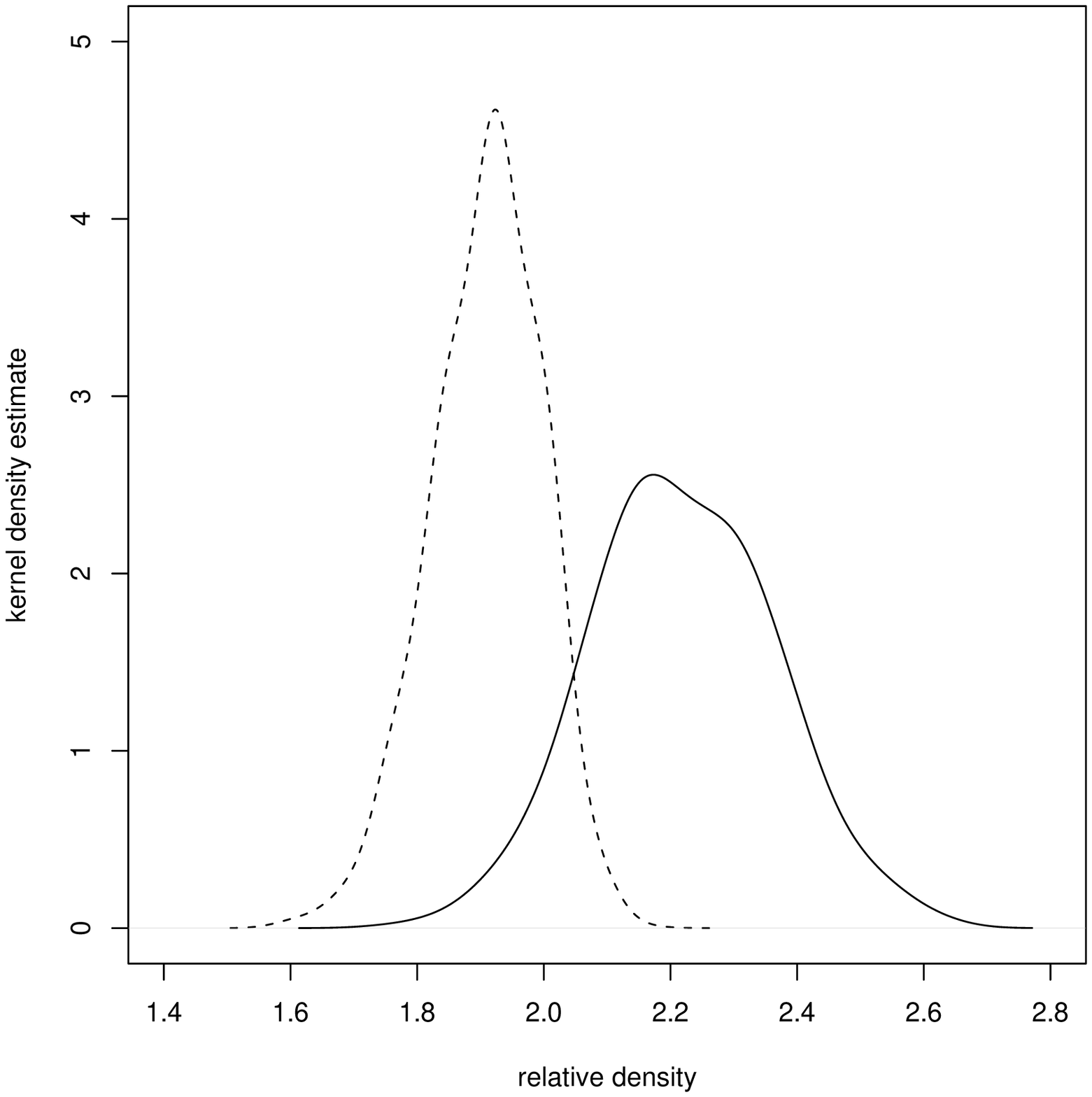} } }
\rotatebox{-90}{ \resizebox{2. in}{!}{ \includegraphics{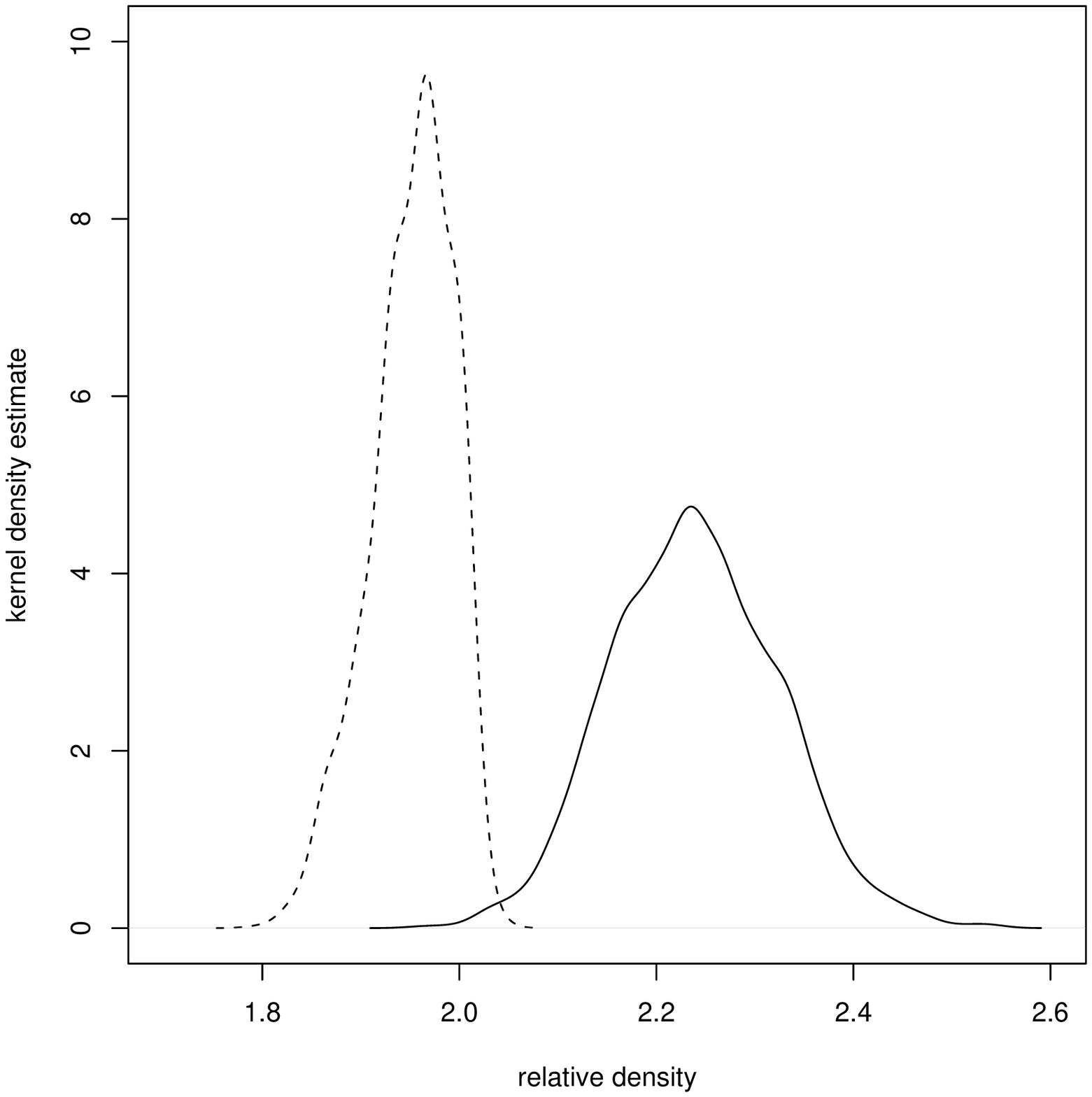} } }
\caption{
\label{fig:CSRvsSeg}
Two Monte Carlo experiments against the segregation alternatives $H^S_{\sqrt{3}/8}$
with $\delta=1/16$.
Depicted are kernel density estimates of $\overline G(r=3/2,M=M_C)$ for $J=13$ and $n=1000$
with $1000$ replicates (left) and $J_{20}=30$ and $n=5000$ with $1000$ replicates (right)
under the null (solid) and segregation alternative (dashed).
}
\end{figure}

\begin{figure}[ht]
\centering
\psfrag{kernel density estimate}{ \Huge{\bf{kernel density estimate}}}
\psfrag{relative density}{ \Huge{\bf{relative density}}}
\rotatebox{-90}{ \resizebox{2. in}{!}{ \includegraphics{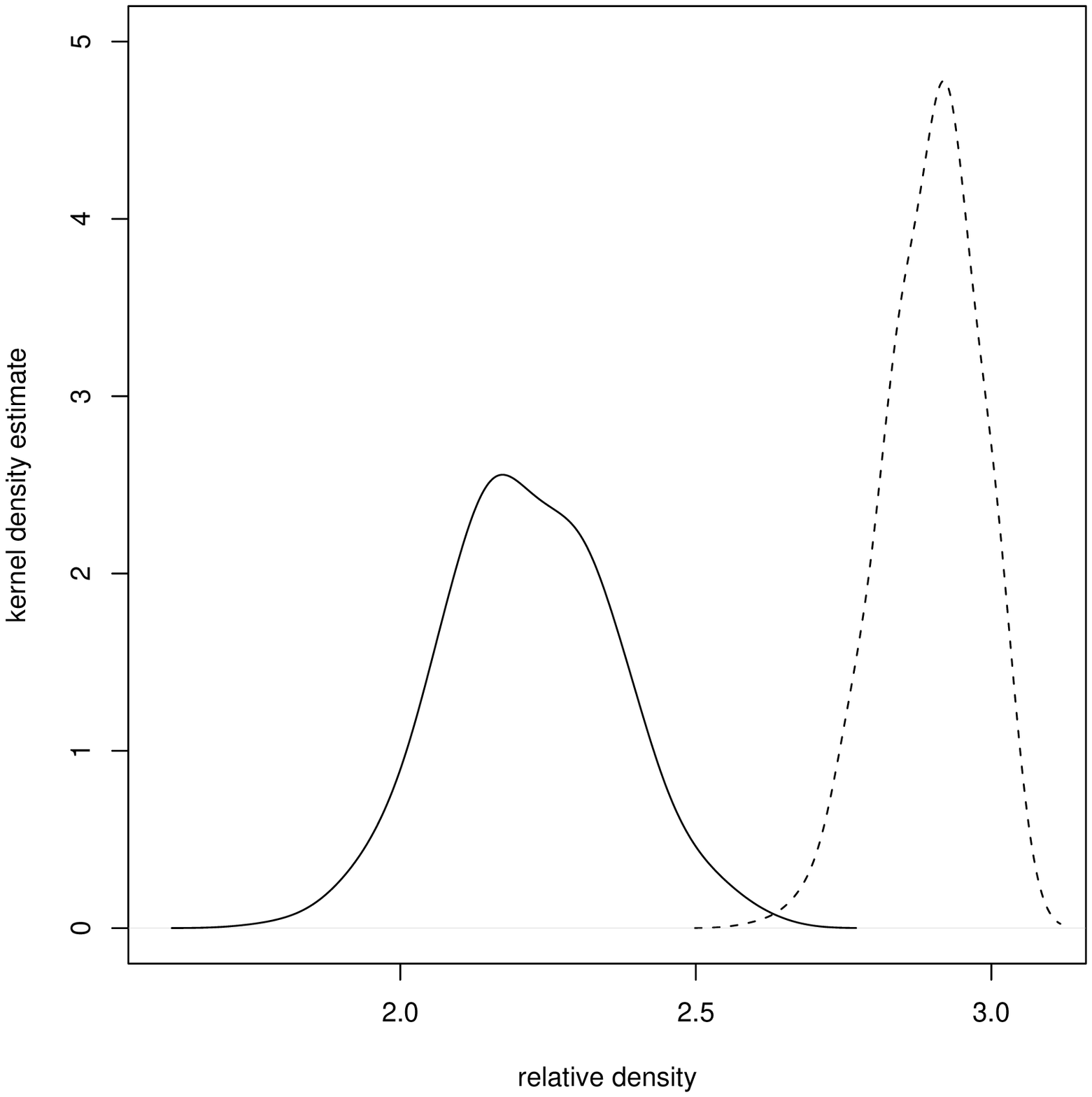} } }
\rotatebox{-90}{ \resizebox{2. in}{!}{ \includegraphics{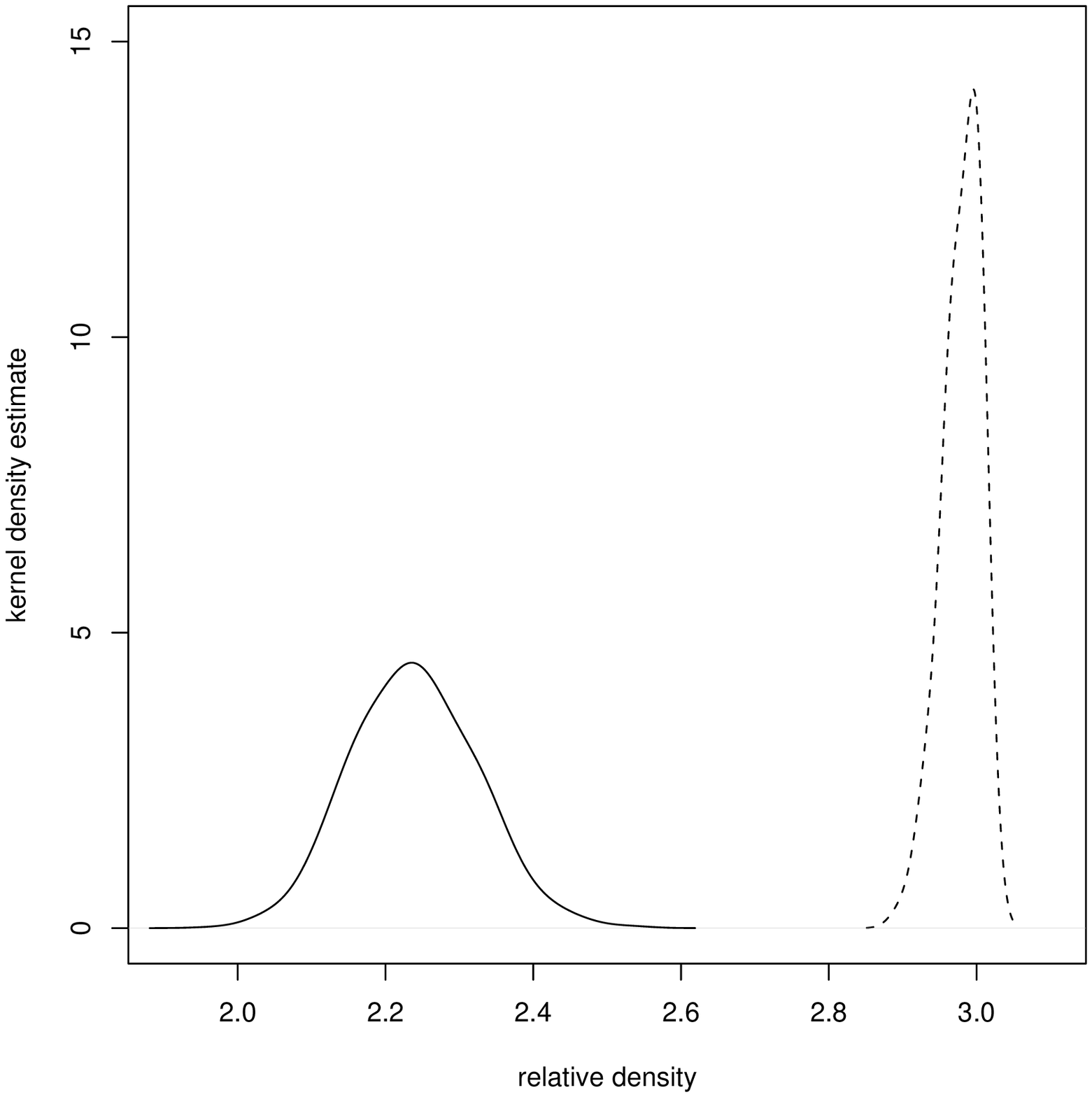} } }
\caption{
\label{fig:CSRvsAgg}
Two Monte Carlo experiments against the association alternatives $H^A_{\sqrt{3}/21}$, i.e.
$\delta=16/49$. Depicted are kernel density estimates of $\overline G(r=3/2,M=M_C)$ for $J=13$ and
$n=1000$ with $500$ replicates (left) and $J_{20}=30$ and $n=5000$ with $100$
replicates under the null (solid) and association alternative (dashed).
}
\end{figure}

To compare the distribution of the test statistic under CSR, and the segregation
and association alternatives,
we generate $n$ points iid $\U(C_H(\Y_m))$ under CSR,
iid uniformly on the support that corresponds to $H^S_{\sqrt{3}/8}$
for each triangle based on the same $\Y_m$ points,
and iid uniformly on the support that corresponds to $H^A_{\sqrt{3}/21}$
for each triangle based on the same $\Y_m$ points.
Under each case, we generate $n=1000$ points with
$J_{10}=13$ and $n=5000$ points with $J_{20}=30$ for 500 Monte Carlo replicates.
The kernel density estimates of $\overline G(r=3/2,M=M_C)$
are presented in Figures \ref{fig:CSRvsSeg} and \ref{fig:CSRvsAgg}.
In Figure \ref{fig:CSRvsSeg},
we observe empirically that even under mild segregation
we obtain considerable separation between the kernel density estimates
under null and segregation cases for moderate $J_m$ and $n$
values suggesting high power at $\alpha=.05$.
A similar result is observed for association.
With $J_{10}=13$ and $n=1000$, under $H_o$,
the estimated significance level is $\widehat{\alpha}=.09$ relative to segregation,
and $\widehat{\alpha}=.07$ relative to association.
Under $H^S_{\sqrt{3}/8}$, the empirical power (using the asymptotic critical value)
is $\widehat{\beta}=.97$, and under $H^A_{\sqrt{3}/21}$, $\widehat{\beta}=1.00$.
With $J_{20}=30$ and $n=5000$, under $H_o$,
the estimated significance level is $\widehat{\alpha}=.06$ relative to segregation,
and $\widehat{\alpha}=.04$ relative to association.
The empirical power is $\widehat{\beta}=1.00$ for both alternatives.

We also estimate the empirical power by using the empirical critical values.
With $J_{10}=13$ and $n=1000$, under $H^S_{\sqrt{3}/8}$,
the empirical power is $\widehat{\beta}_{mc}=.72$ at empirical level
$\widehat{\alpha}_{mc}=.033$ and under $H^A_{\sqrt{3}/21}$ the
empirical power is $\widehat{\beta}_{mc}=1.00$ at empirical
level $\widehat{\alpha}_{mc}=.03$.
With $J_{20}=30$ and $n=5000$, under $H^S_{\sqrt{3}/8}$,
the empirical power is $\widehat{\beta}_{mc}=1.00$
at empirical level $\widehat{\alpha}_{mc}=.034$ and
under $H^A_{\sqrt{3}/21}$ the empirical power is
$\widehat{\beta}_{mc}=1.00$ at empirical level
$\widehat{\alpha}_{mc}=.04$.

In Figure \ref{fig:CSRvsAgg}, we observe that even in mild association we obtain considerable separation
for moderate $J_m$ and $n$ values suggesting high power (with $J_{10}=13$ and $n=1000$,
the empirical critical value is $2.46$, $\widehat{\alpha}=.034$ and
empirical power is $\widehat{\beta}=1.0$ and with $J_{20}=30,\;\;n=5000$,
the empirical critical value is $2.36$, $\widehat{\alpha}=.04$
and empirical power is $\widehat{\beta}=1.0$).

\begin{figure}[ht]
\centering
\rotatebox{-90}{ \resizebox{2. in}{!}{\includegraphics{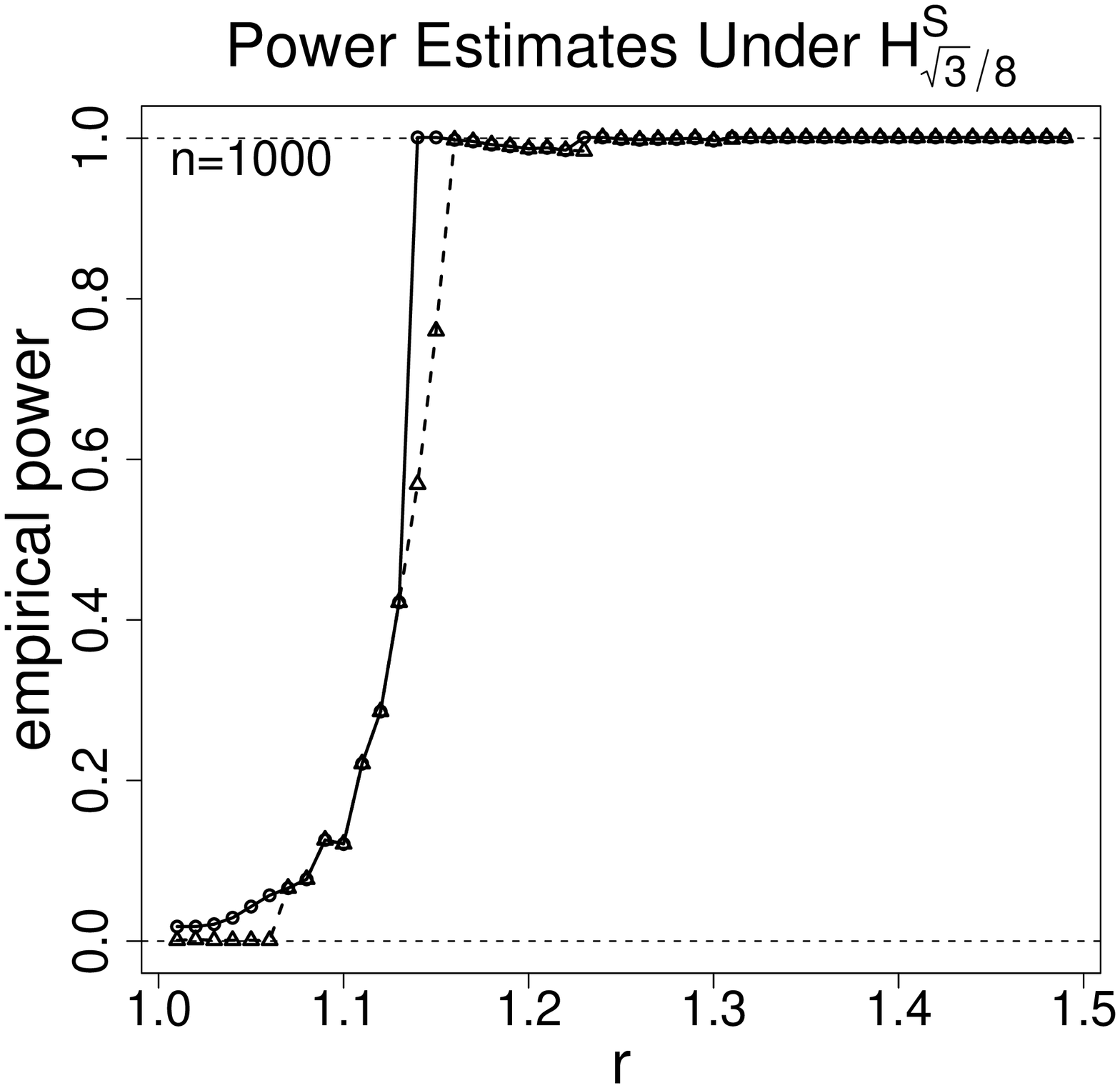} }}
\rotatebox{-90}{ \resizebox{2. in}{!}{\includegraphics{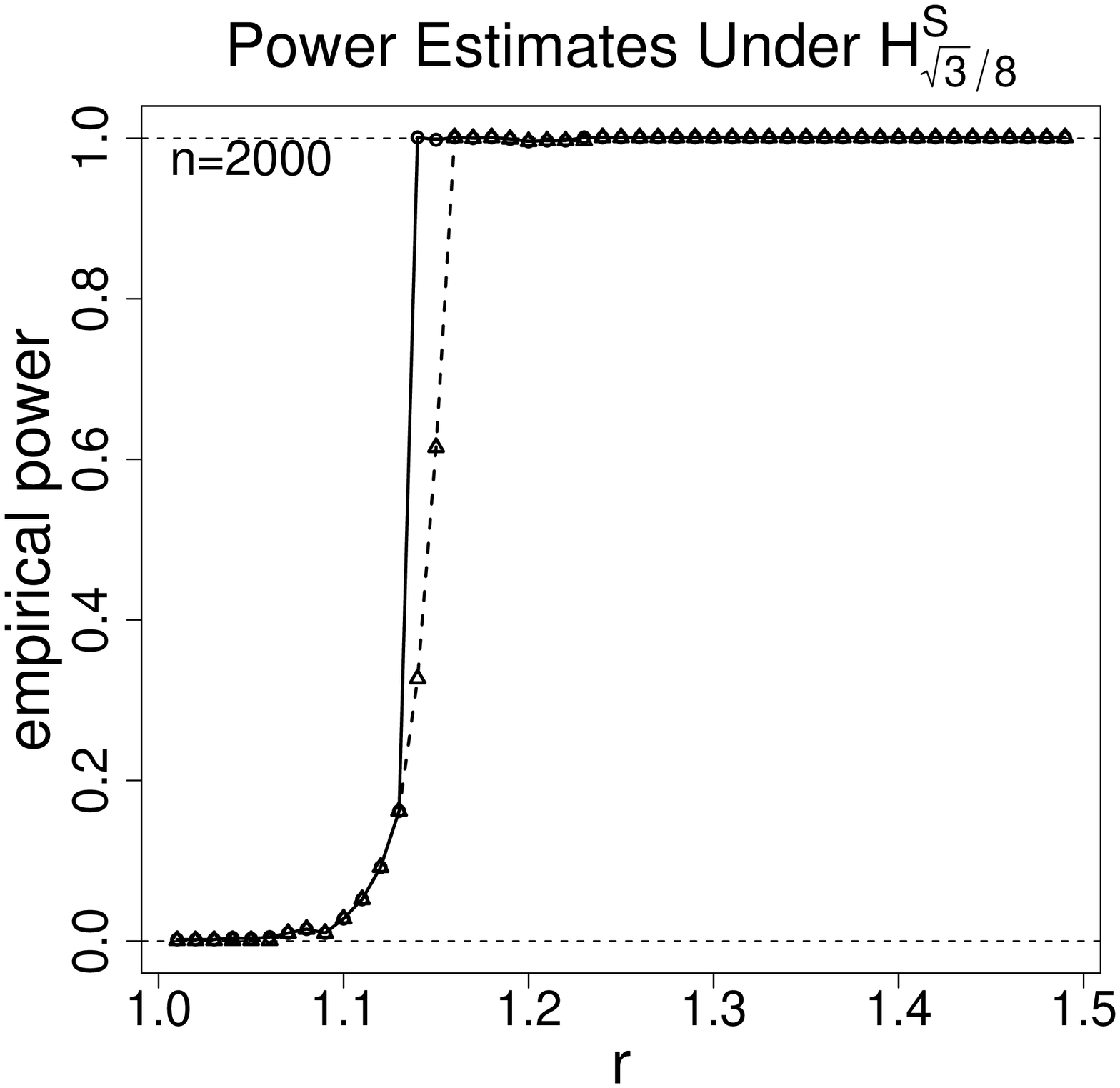} }}
\rotatebox{-90}{ \resizebox{2. in}{!}{\includegraphics{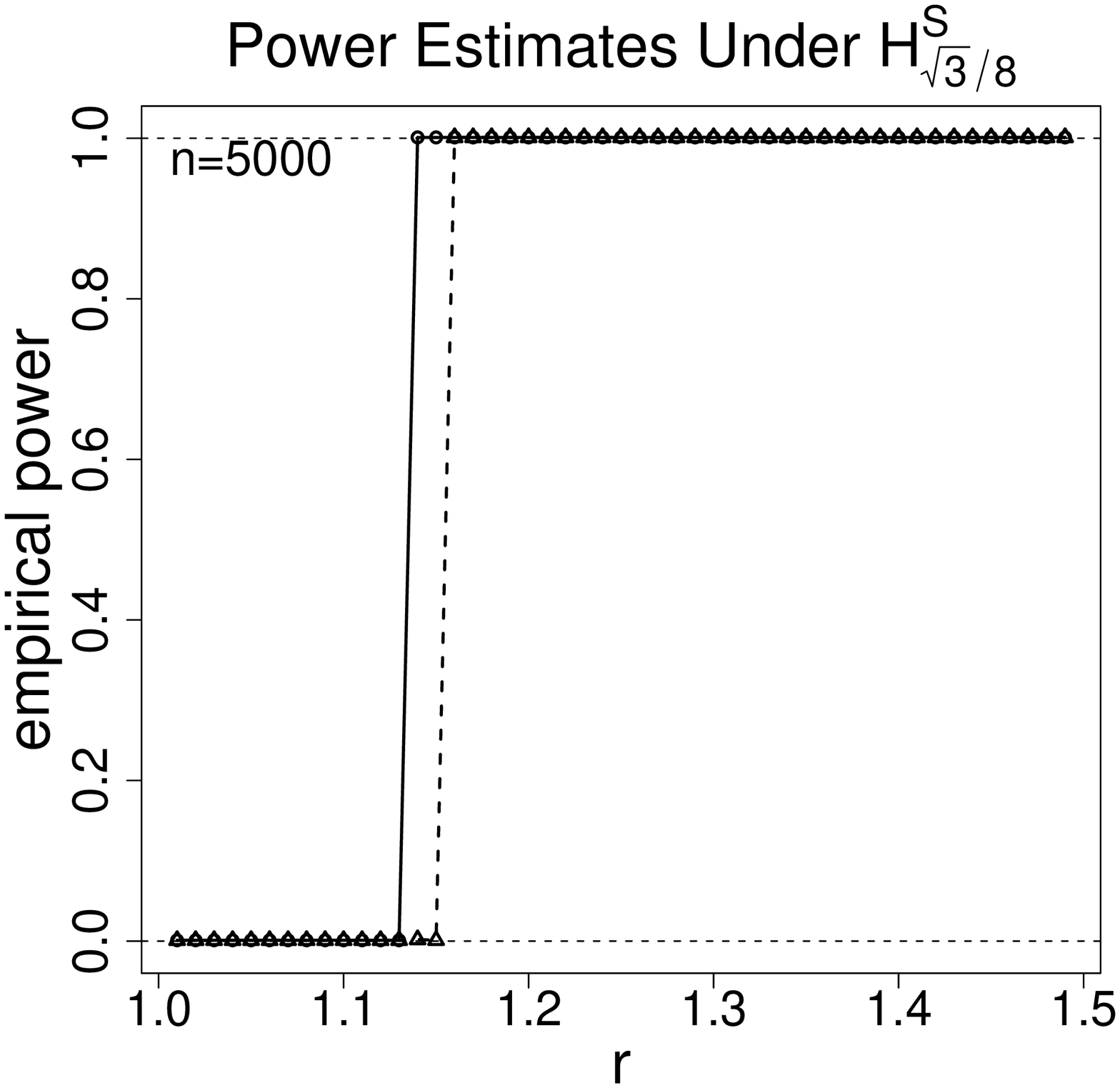} }}
\rotatebox{-90}{ \resizebox{2. in}{!}{\includegraphics{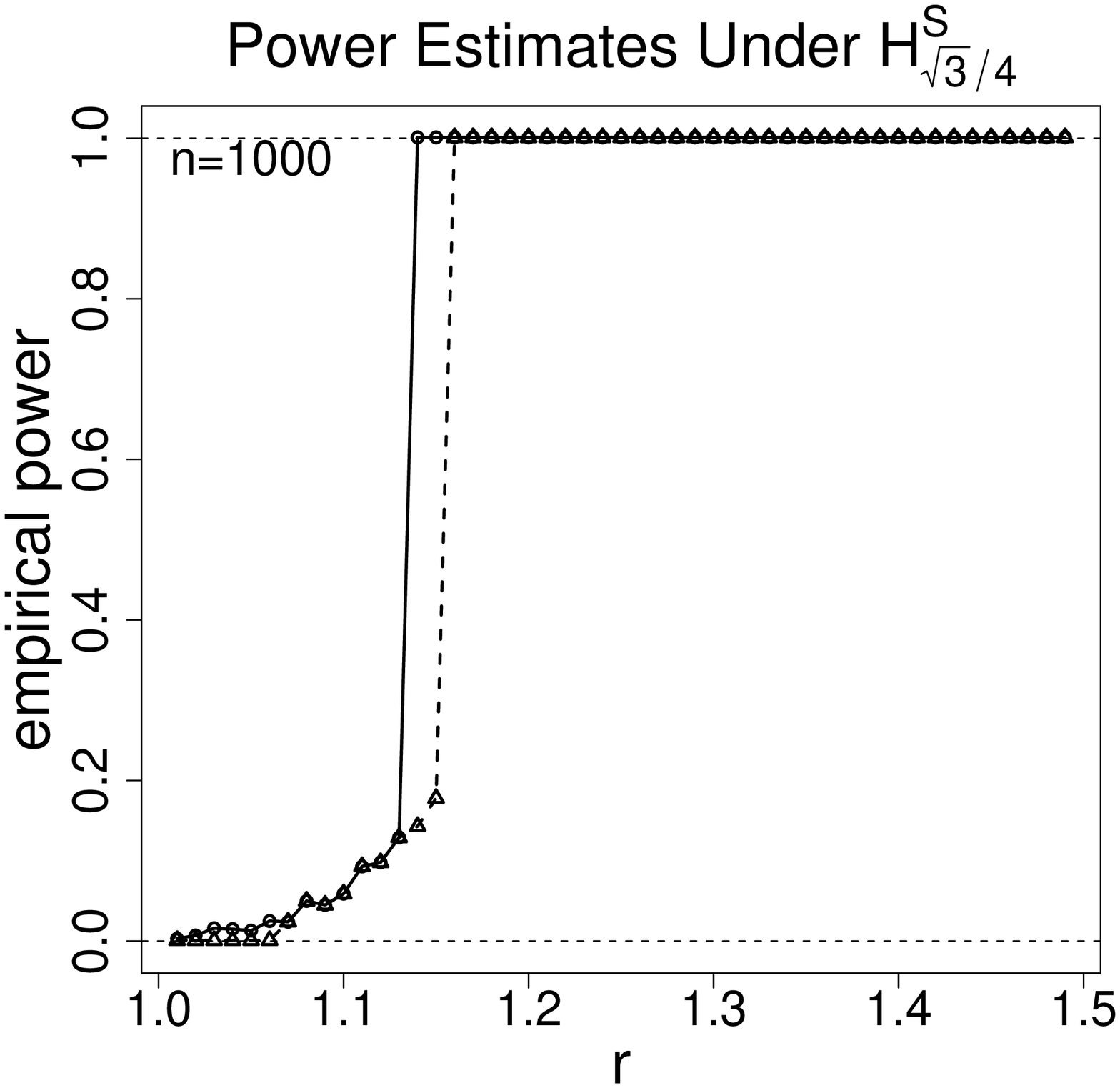} }}
\rotatebox{-90}{ \resizebox{2. in}{!}{\includegraphics{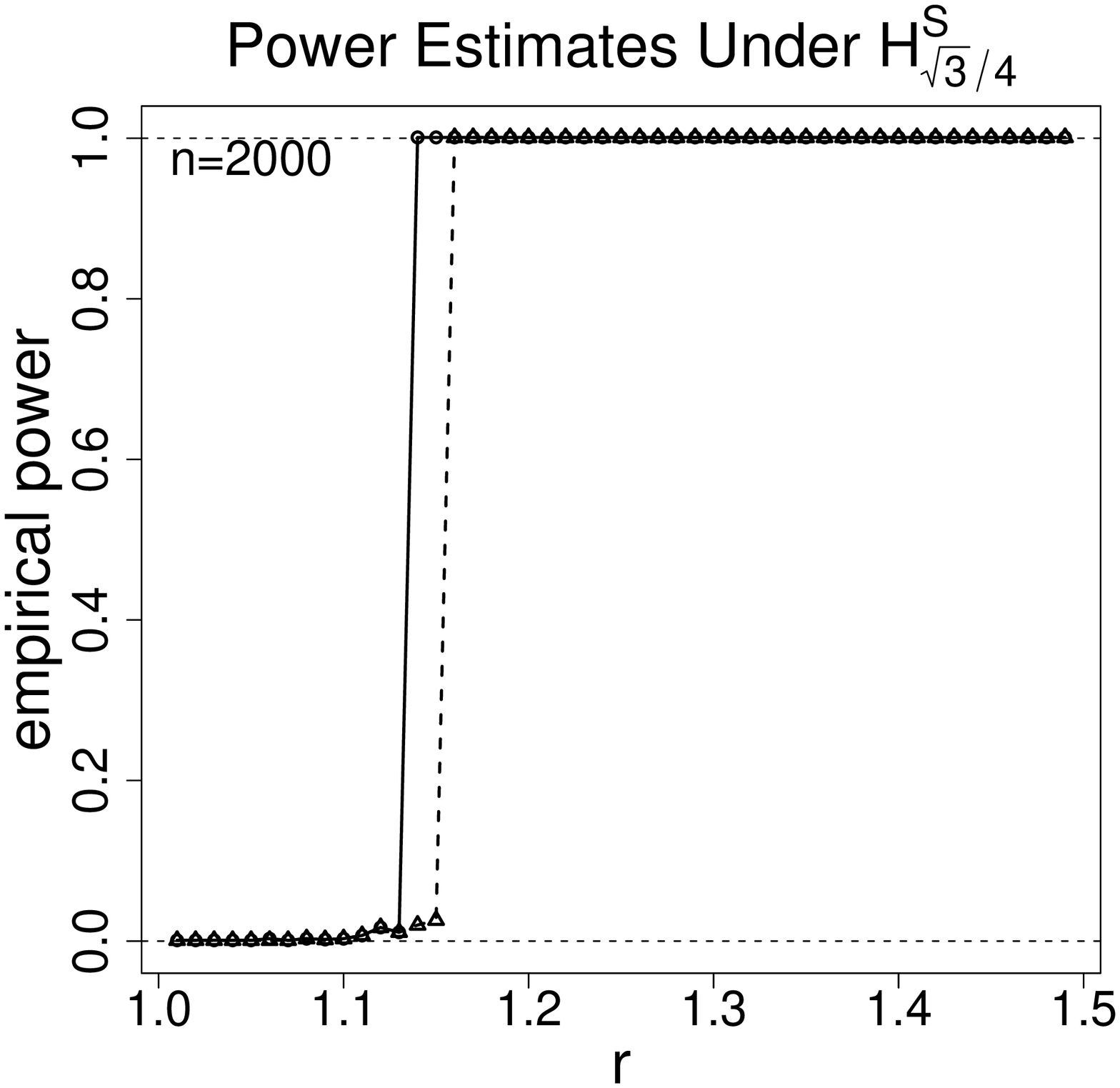} }}
\rotatebox{-90}{ \resizebox{2. in}{!}{\includegraphics{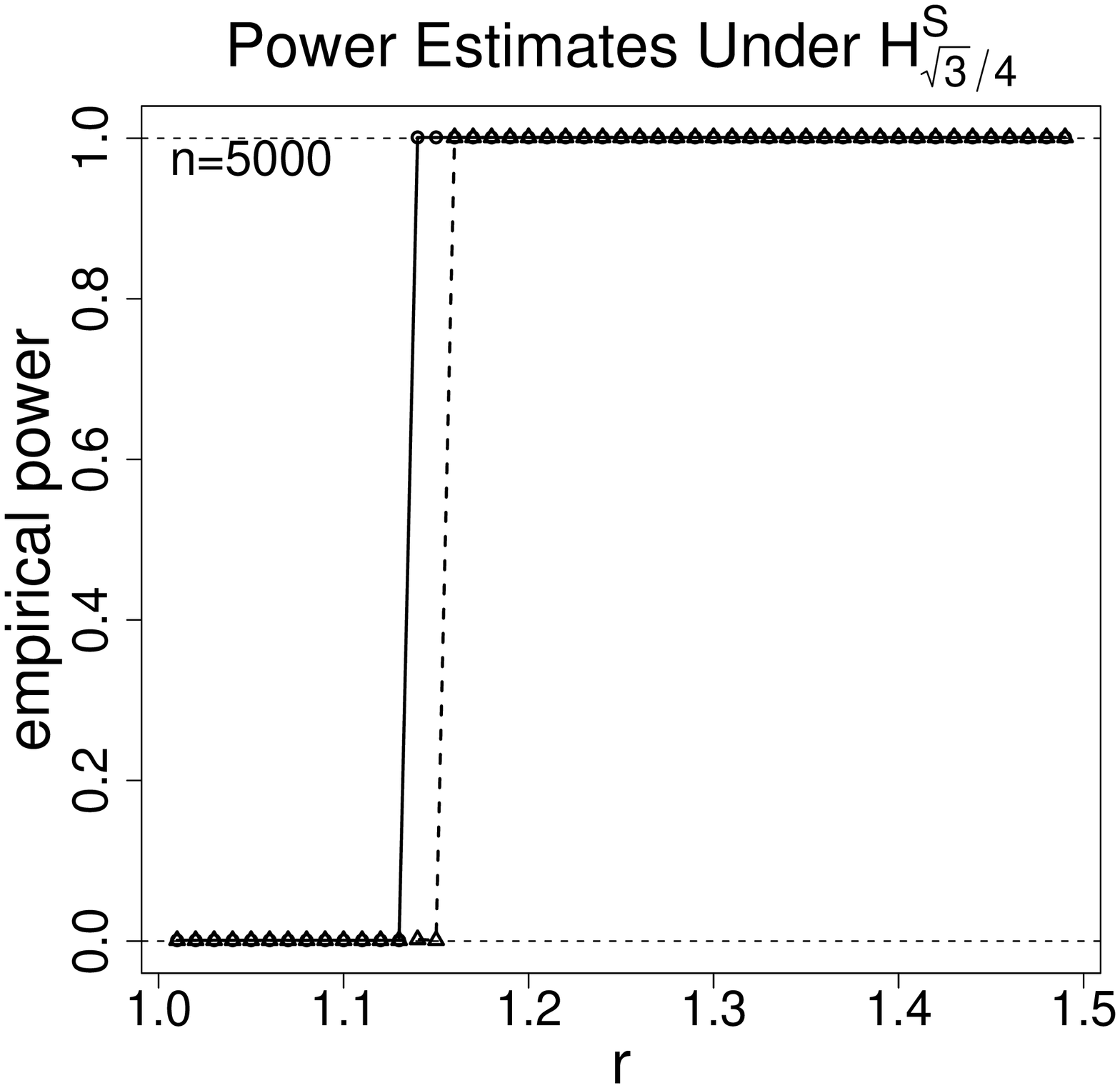} }}
\caption{
\label{fig:power-seg}
The empirical power estimates under segregation with $\ve=\sqrt{3}/8,\ve=\sqrt{3}/4$
and $n=1000$ (left), $n=2000$ (middle), and $n=5000$ (right).
The power estimates based on the binomial distribution are plotted in circles ($\circ$) and joined with solid lines,
and those based on the normal approximation are plotted in triangles ($\triangle$) and joined with dashed lines.
}
\end{figure}

For the segregation alternatives,
we consider the following three cases:
$\ve=\sqrt{3}/8,\ve=\sqrt{3}/4,\ve=2\sqrt{3}/7$ in the 13 Delaunay triangles
obtained by the 10 $\Y$ points in Figure \ref{fig:deltri}.
We generate $n=500,1000,2000,5000$ in the convex hull of $\Y_{10}$ at each Monte Carlo replication.
We estimate the empirical power of the tests for $r=1.00,1.01,1.02,\ldots,1.49$ values using $N_{mc}=1000$ replicates.
The power estimates based on the binomial distribution and normal approximation under $H^S_{\sqrt{3}/8}$
for $n=1000,2000,5000$ are plotted in Figure \ref{fig:power-seg}.
Observe that the power estimates are about 1.0 for $r \gtrsim 1.15$.
Considering the empirical size and power estimates together,
we recommend $r$ values around 1.22 or 1.30 for the segregation alternatives.

\begin{figure}[ht]
\centering
\rotatebox{-90}{ \resizebox{2. in}{!}{\includegraphics{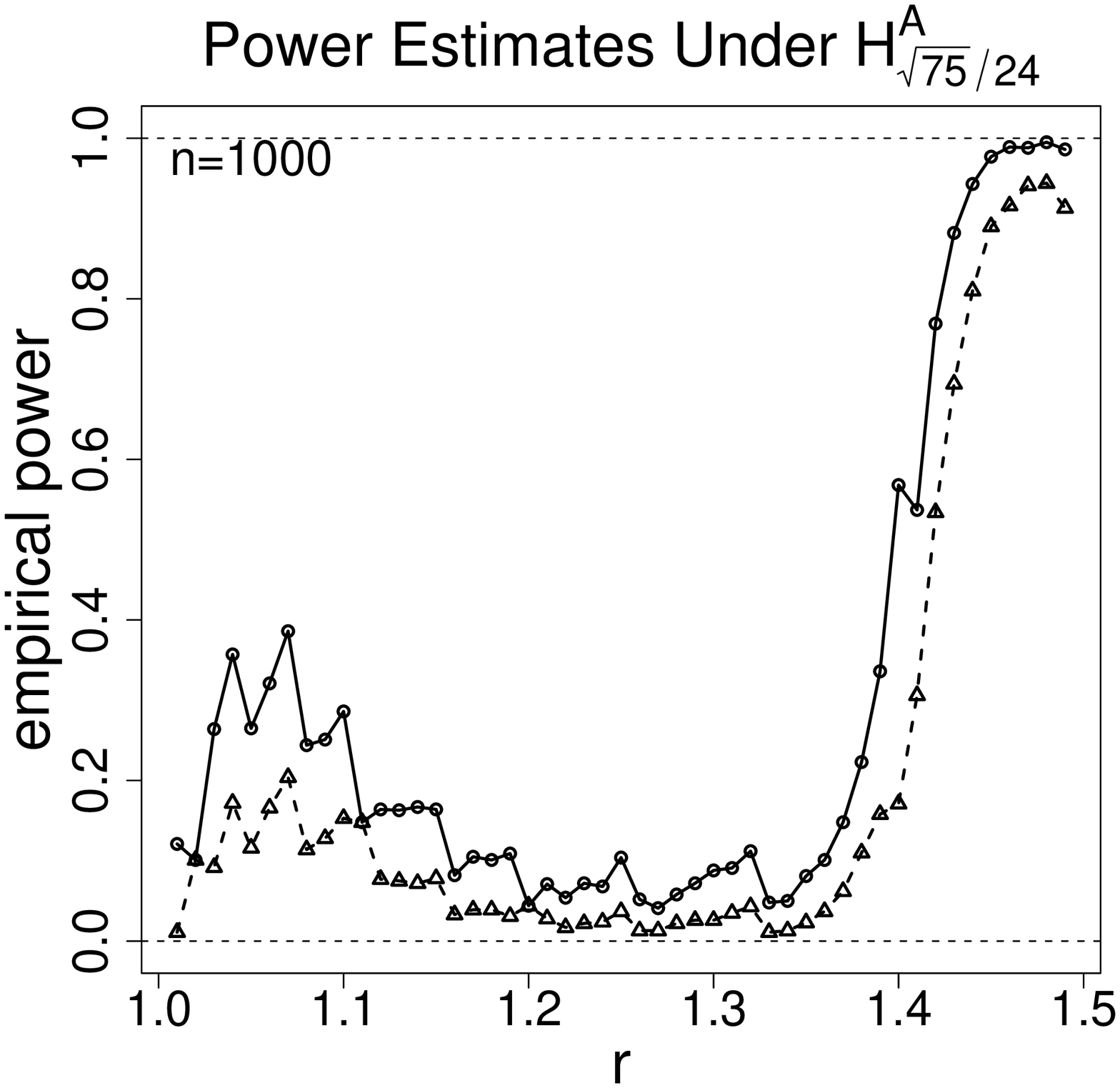} }}
\rotatebox{-90}{ \resizebox{2. in}{!}{\includegraphics{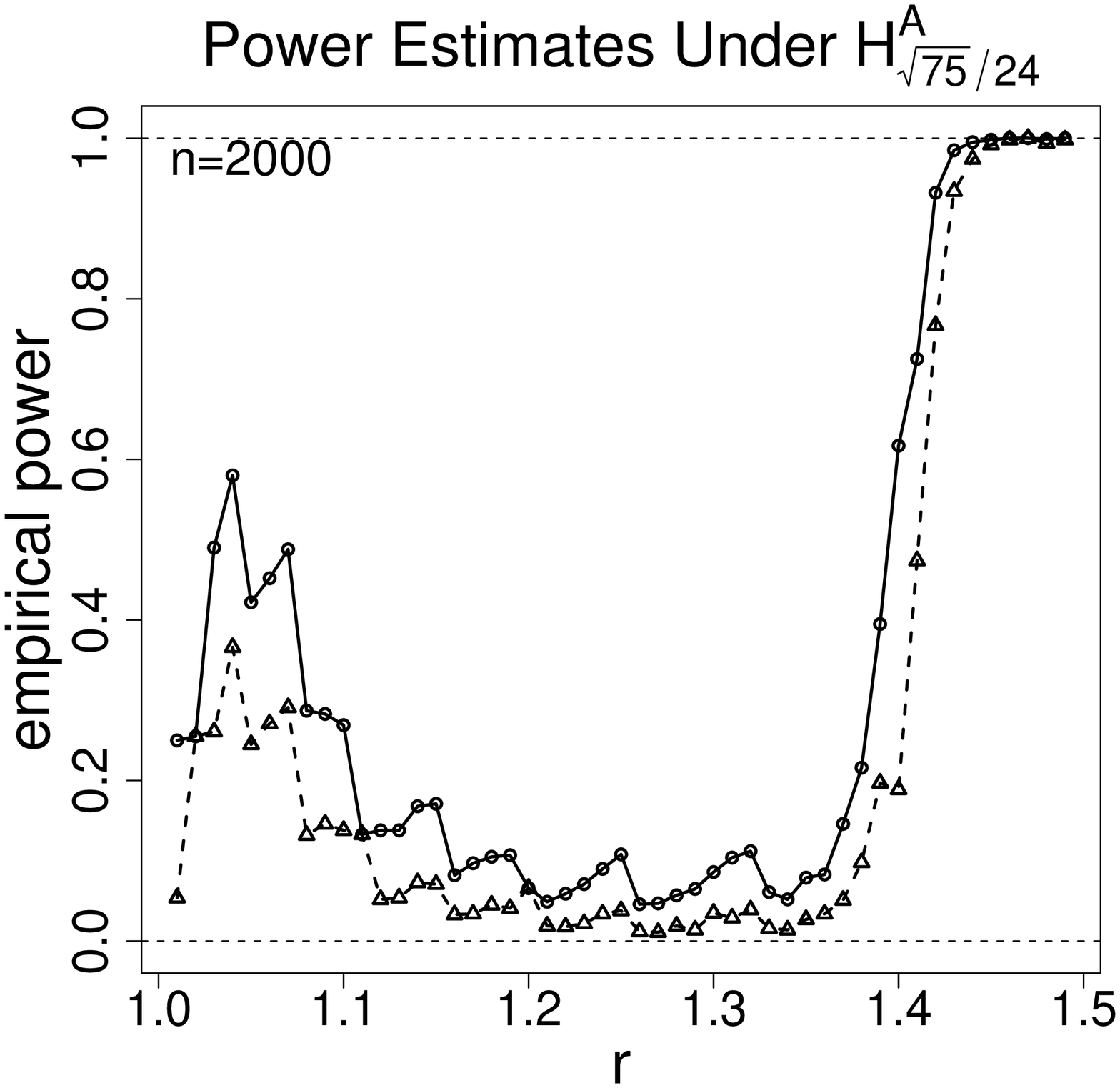} }}
\rotatebox{-90}{ \resizebox{2. in}{!}{\includegraphics{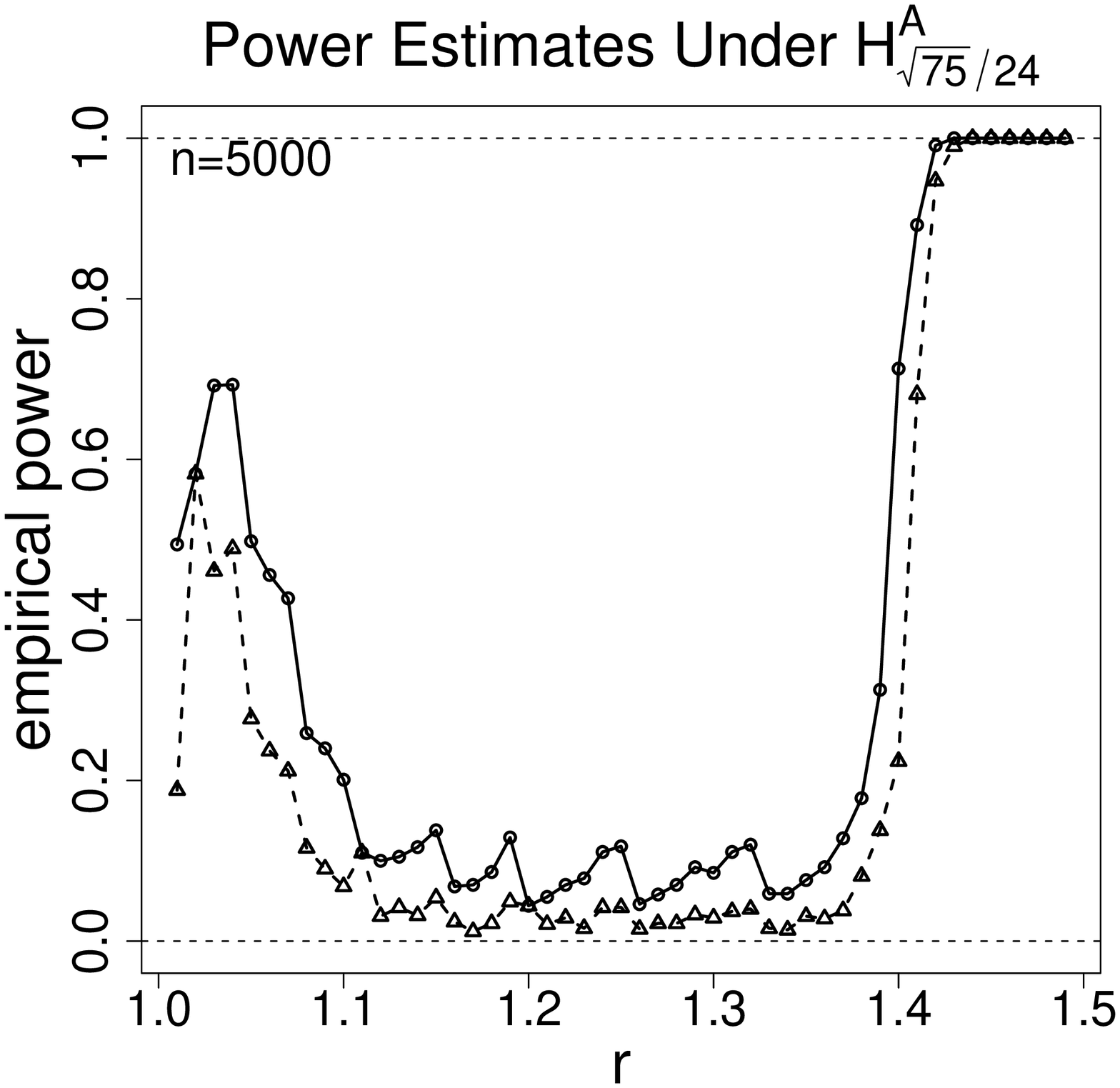} }}
\rotatebox{-90}{ \resizebox{2. in}{!}{\includegraphics{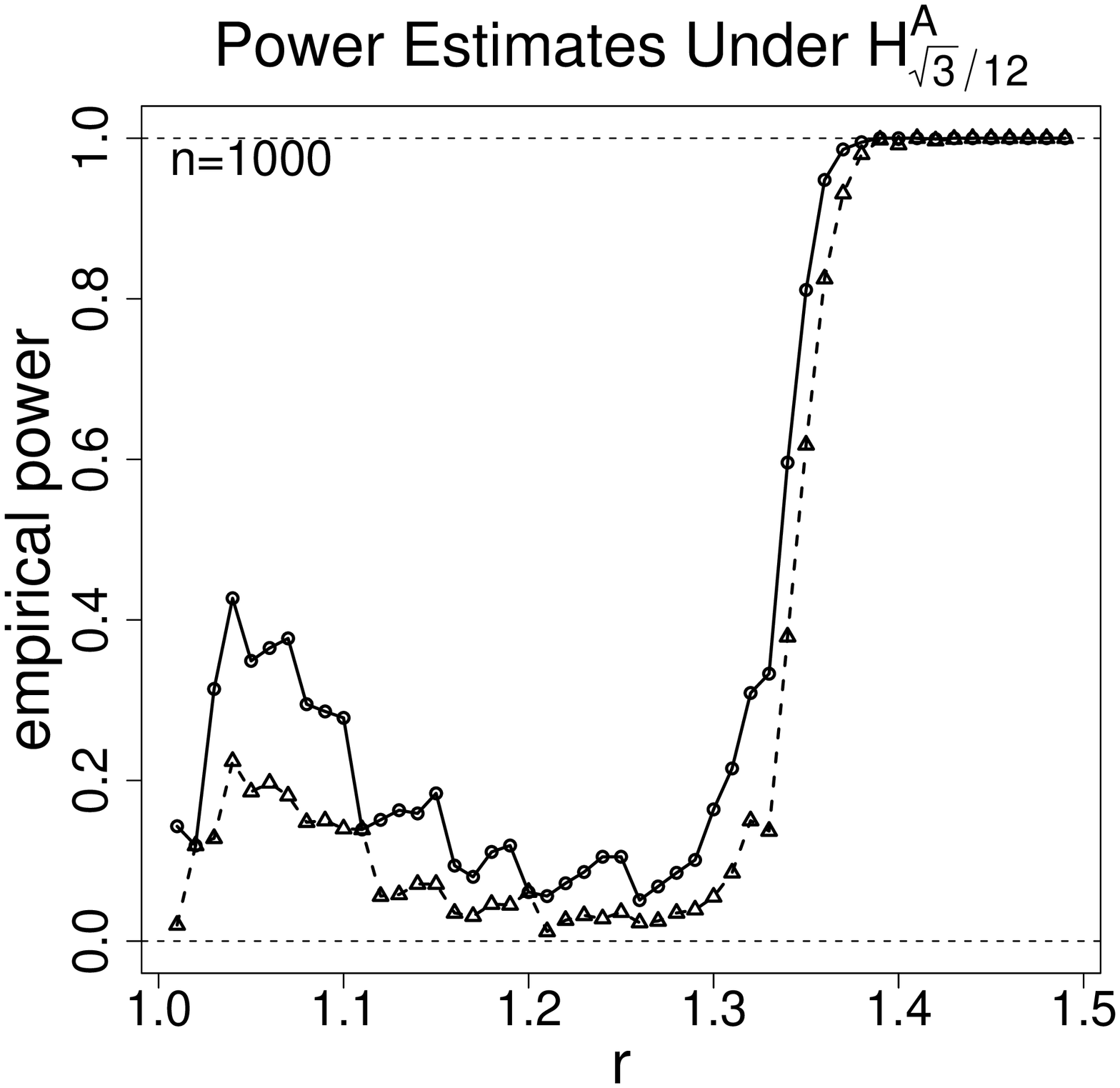} }}
\rotatebox{-90}{ \resizebox{2. in}{!}{\includegraphics{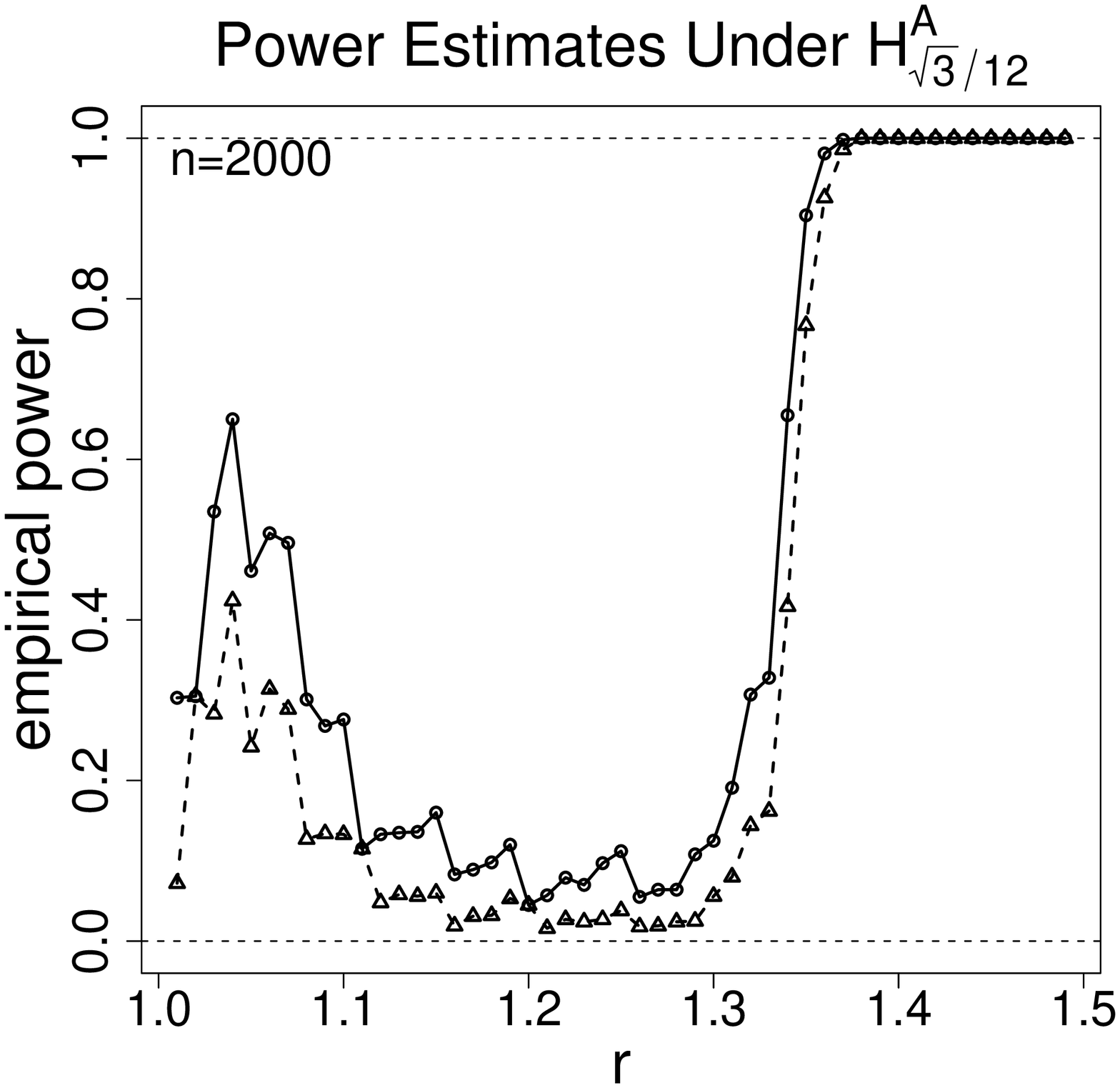} }}
\rotatebox{-90}{ \resizebox{2. in}{!}{\includegraphics{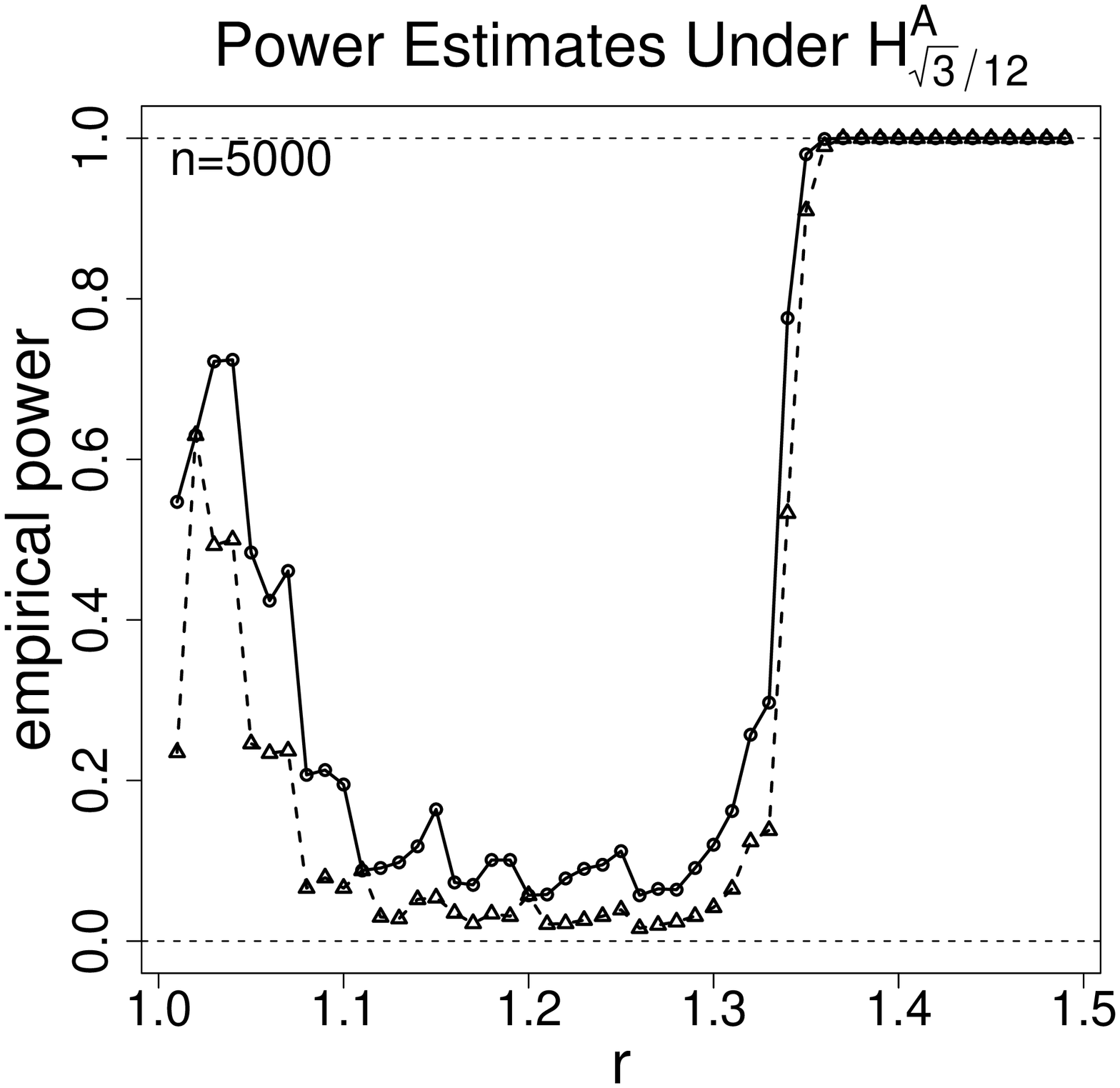} }}
\caption{
\label{fig:power-assoc}
The empirical power estimates under association with $\ve=5\sqrt{3}/24,\ve=\sqrt{3}/12$
and $n=1000$ (left), $n=2000$ (middle), and $n=5000$ (right).
The power estimates based on the binomial distribution are plotted in circles ($\circ$) and joined with solid lines,
and those based on the normal approximation are plotted in triangles ($\triangle$) and joined with dashed lines.
}
\end{figure}

For the association alternatives,
we consider the following three cases:
$\ve=5\sqrt{3}/24,\ve=\sqrt{3}/12,\ve=\sqrt{3}/21$ in the 13 Delaunay triangles
obtained by the 10 $\Y$ points in Figure \ref{fig:deltri}.
We generate $n=500,1000,2000,5000$ in the convex hull of $\Y_{10}$ at each Monte Carlo replication.
We estimate the empirical power of the tests for $r=1.00,1.01,1.02,\ldots,1.49$ values using $N_{mc}=1000$ replicates.
The power estimates based on the binomial distribution and normal approximation under $H^A_{5\sqrt{3}/24}$
for $n=1000,2000,5000$ are plotted in Figure \ref{fig:power-assoc}.
Observe that the power estimates are about 1.0 for $r \gtrsim 1.33$,
but the power performance is poor for $r$ between 1.1 and 1.33.
Considering the empirical size and power estimates together,
we recommend $r$ values around 1.35 for the association alternatives.

The empirical power estimates for $r=3/2$ and $M=M_C$ are presented in Table \ref{tab:emp-size-power-r=3/2}.

\begin{table}[ht]
\begin{center}
\begin{tabular}{|c|c|c||c|c|c|c||c|c|}
\hline
\multicolumn{9}{|c|}{Empirical Size and Power Estimates for $r=3/2$ and $M=M_C$} \\
\hline
$n$  & $\ah^S$ & $\ah^A$ & $\bh^S_1$ &  $\bh^S_2$ & $\bh^S_3$ & $\bh^A_1$ & $\bh^A_2$ & $\bh^A_3$\\
\hline
500  & 0.161 & 0.062 & 0.961 & 1.000 & 1.000 & 1.000 & 1.000 & 0.997\\
\hline
1000 &  0.071 & 0.082 & 0.975 & 1.000 & 1.000 & 1.000 & 1.000 & 1.000\\
\hline
2000  & 0.049 & 0.081 & 0.995 & 1.000 & 1.000 & 1.000 & 1.000 & 1.000\\
\hline
\end{tabular}
\end{center}
\caption{\label{tab:emp-size-power-r=3/2}
The empirical size and power estimates for $r=3/2$ and $M=M_C$
under the null and alternatives.
$n$ stands for size of $\X$ points, $\ah^S$ for empirical size relative segregation,
$\ah^A$ for empirical size relative to association,
$\bh^S_1$, $\bh^S_2$, and $\bh^S_3$ for empirical power estimates under $H^S_{\ve}$
with $\ve=\sqrt{3}/8,\ve=\sqrt{3}/4,$ and $\ve=2\sqrt{3}/7$, respectively,
$\bh^A_1$, $\bh^A_2$, and $\bh^A_3$ for empirical power estimates under $H^A_{\ve}$
with $\ve=5\sqrt{3}/24,\ve=\sqrt{3}/12,$ and $\ve=\sqrt{3}/21$, respectively.}
\end{table}

\section{Correction for $\X$ Points Outside the Convex Hull of $\Y_m$}
\label{sec:conv-hull-correction}
Our null hypothesis in \eqref{eqn:null-pattern-mult-tri} is rather restrictive,
in the sense that,
it might not be that realistic to assume the support of $\X$
being $C_H(\Y_m)$ in practice.
Up to now, our inference is restricted to the $C_H(\Y_m)$.
However, crucial information from the data (hence power) might be lost
since a substantial proportion of $\X$ points,
denoted $\pi_{out}$,
might fall outside the $C_H(\Y_m)$.
We investigate the effect of $\pi_{out}$ (or restriction to the $C_H(\Y_m)$)
on our tests and propose an empirical correction to mitigate this
based on an extensive Monte Carlo simulation study.

We consider the following 6 cases to investigate how the removal of points outside
$C_H(\Y_m)$ affects the empirical size and power performance of the tests.
We only consider $r=1.35$ and $r=1.5$ which have better size and power performances
compared to others.
In each case,
at each Monte Carlo replication,
we generate $\X_n$ and $\Y_m$ independently
as random samples from  $\U(\mS_X)$ and $\U(\mS_Y)$, respectively,
for various values of $n$ and $m$
where $\mS_X$ and $\mS_Y$ are the support sets of $\X$ and $\Y$ points, respectively.
We take  $\mS_Y=(0,1)\times(0,1)$ and
manipulate $\mS_X$ in each case to simulate CSR and various forms of deviations from CSR.
We repeat the generation procedure $N_{mc}$ times for each combination of $m$ and $n$.
At each Monte Carlo replication,
we record the number of $\X$ points outside $C_H(\Y_m)$ and the domination number, $\g_{m,n}(r)$.

\begin{itemize}
\item[Case 1:]
In this case, we also set $\mS_X=(0,1)\times(0,1)$, 

\item[Case 2:]
$\mS_X=(-\delta,1+\delta)\times(-\delta,1+\delta)$ for $\delta \in \{.1,.25,.5\}$,

\item[Case 3:]
$\mS_X=(0,1)\times(0,1+\delta)$ for $\delta \in \{.1,.25,.5\}$,

\item[Case 4:]
$\mS_X=(0,1)\times(\delta,1+\delta)$ for $\delta \in \{.1,.25,.5\}$.

\item[Case 5:]
Given a realization of $\Y$ points, $\Y_m = \{\y_1,\y_2,\ldots,\y_m\}$, from $\U(\mS_Y=(0,1)\times(0,1))$,
$\mS_X=[(-\delta,1+\delta)\times(-\delta,1+\delta)]\setminus \bigcup_{i=1}^m B(\y_i,\ve)$
with $\delta=\frac{1}{2\,\sqrt{\lambda}}=\frac{1}{2\,\sqrt{m}}$ which the expected interpoint distance
in a homogeneous Poisson process with intensity (expected number of points per unit area) $\lambda$ (\cite{dixon:EncycEnv2002})
and $\ve=\delta/k$ for $k=1.5,2.0$,

\item[Case 6:]
Given a realization of $\Y$ points, $\Y_m = \{\y_1,\y_2,\ldots,\y_m\}$, from $\U(\mS_Y)$,
$\mS_X=\bigcup_{i=1}^m B(\y_i,\ve)$ with $\ve=\delta/k$, $\delta=\frac{1}{2\,\sqrt{m}}$,
and $k=1.0,1.5$.
\end{itemize}

Notice that in Case 1
both $\X_n$ and $\Y_m$ have the same support.
By construction the two classes follow CSR independence with
very different relative abundances
(i.e., number of $\X$ points being larger than number of $\Y$ points).
In Cases 2 and 3
the support of $\X_n$ contains (but larger than) the support of $\Y_m$,
which suggests segregation of $\X$ points from $\Y$ points,
at least when we move away from the support of $\Y$ points (which is the unit square).
However, when we restrict our attention to $C_H(\Y_m)$ or the unit square,
we have CSR or CSR independence, respectively.
Furthermore, the larger the $\delta$ value,
the larger the level of segregation of $\X$ from $\Y$.
In Case 4 the support of $\X_n$ and $\Y_m$ overlap,
but neither is a subset of the other,
which suggests segregation between $\X$ and $\Y$ points.
When we restrict our attention to $C_H(\Y_m)$,
there is still segregation between $\X$ and $\Y$ points.
Furthermore, the larger the $\delta$ value,
the larger the level of segregation between $\X$ and $\Y$ points.
In Case 5,
$\X$ points are segregated from $\Y$ points
both in and outside $C_H(\Y_m)$.
Furthermore, the larger the $\delta$ value,
the larger the level of segregation of $\X$ points from $\Y$ points.
Finally, in Case 6
$\X$ points are associated with $\Y$ points.
Furthermore, the smaller the $\delta$ value,
the larger the level of association of $\X$ points with $\Y$ points.

In Case 1 (i.e., the benchmark case),
we consider $n=100,200,\ldots,900,1000,2000,\ldots,9000,10000$
for each of $m=10,20,\ldots,50$.
We generate $N_{mc}=1000$ replication for each $n,m$ combination.
In the other cases, we consider $n=100,500,1000$ for $m=10$
and $n=500,1000$ for $m=20$;
and we generate $N_{mc}=10000$ replication for each $n,m$ combination.

In Cases 1-6,
we estimate the proportion of $\X$ points outside the $C_H(\Y_m)$.
For each $m,n$ combination
we average (over $n$) this proportion
which is denoted as $\widehat \pi_{out}$.
We present the estimated (mean) proportions $\widehat \pi_{out}$
for Case 1 in Table \ref{tab:prop-out-CH}
and for Cases 2-6 in Table \ref{tab:prop-out-CH-cases2-6}.
Observe that in Cases 2-5,
$\widehat \pi_{out}$ values are larger than that in Case 1,
while in Case 6,
$\widehat \pi_{out}$ values are smaller than that in Case 1.

\begin{table}[ht]
\begin{center}
\begin{tabular}{|c|c|c|c|c|c|}
\hline
$m$  & 10 & 20 & 30 &  40 & 50\\
\hline
$\widehat \pi_{out}$  & 0.56 & 0.37 & 0.29 & 0.23 & 0.20\\
\hline
$\widehat \pi_{fit}$  & 0.57 & 0.36 & 0.28 & 0.24 & 0.21\\
\hline
\end{tabular}
\end{center}
\caption{\label{tab:prop-out-CH}
The (mean) proportion of $\X$ points outside the $C_H(\Y_m)$
which is denoted as $\widehat \pi_{out}$
and the fitted values $\widehat \pi_{fit}$ for various $m$ values in Case 1.
}
\end{table}

\begin{table}[ht]
\begin{center}

\begin{tabular}{|c|c|c|c|}
\hline
\multicolumn{4}{|c|}{$\widehat \pi_{out}$ values for Case 2} \\
\hline
$\delta$  & 0.1 & 0.25 & 0.50\\
\hline
$m=10$ & 0.697 & 0.806 & 0.891\\
\hline
$m=20$ & 0.566 & 0.722 & 0.843\\
\hline
\end{tabular}
\hspace{1 cm}
\begin{tabular}{|c|c|c|c|}
\hline
\multicolumn{4}{|c|}{$\widehat \pi_{out}$ values for Case 3} \\
\hline
$\delta$  & 0.1 & 0.25 & 0.50\\
\hline
$m=10$ & 0.604 & 0.652 & 0.740\\
\hline
$m=20$ & 0.431 & 0.499 & 0.582\\
\hline
\end{tabular}
\vspace{0.5 cm}

\begin{tabular}{|c|c|c|c|}
\hline
\multicolumn{4}{|c|}{$\widehat \pi_{out}$ values for Case 4} \\
\hline
$\delta$  & 0.1 & 0.25 & 0.50\\
\hline
$m=10$ & 0.573 & 0.629 & 0.782\\
\hline
$m=20$ & 0.395 & 0.488 & 0.687\\
\hline
\end{tabular}
\hspace{0.5 cm}
\begin{tabular}{|c|c|c|}
\hline
\multicolumn{3}{|c|}{$\widehat \pi_{out}$ values for Case 5} \\
\hline
$k$  & 1.5 & 2.0\\
\hline
$m=10$ & 0.806 & 0.783\\
\hline
$m=20$ & 0.652 & 0.611\\
\hline
\end{tabular}
\hspace{0.5 cm}
\begin{tabular}{|c|c|c|}
\hline
\multicolumn{3}{|c|}{$\widehat \pi_{out}$ values for Case 6} \\
\hline
$k$  & 1.0 & 1.5\\
\hline
$m=10$ & 0.535 & 0.479\\
\hline
$m=20$ & 0.358 & 0.310\\
\hline
\end{tabular}

\end{center}
\caption{\label{tab:prop-out-CH-cases2-6}
The (mean) proportion of $\X$ points outside the $C_H(\Y_m)$
for various $\delta$ and $m$ values in Cases 2-4 and
various $k$ and $m$ values in Cases 5-6.
}
\end{table}

\begin{figure}[ht]
\centering
\rotatebox{-90}{ \resizebox{2. in}{!}{\includegraphics{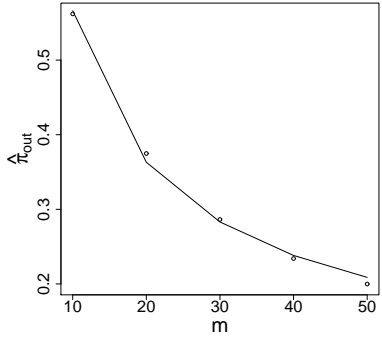} }}
\caption{
\label{fig:prop-out-CH}
The proportion of $X$ points outside $C_H(\Y_m)$
as a function of $m$.
The solid line is the fitted line based on $\pi_{out} \approx 1.7932/m+1.2229/\sqrt{m}$.
}
\end{figure}

For Case 1, we model the relationship between $\widehat \pi_{out}$ and $m$.
Our simulation results suggest that
$\widehat \pi_{out} \approx 1.7932/m+1.2229/\sqrt{m}$.
We present the actual fitted values denoted $\widehat \pi_{fit}$ based on this model
in Table \ref{tab:prop-out-CH}.
See also Figure \ref{fig:prop-out-CH} for the plot of estimated $\widehat \pi_{out}$ values
versus fitted values based on our model.
Notice that as $m \rightarrow \infty$,
$\widehat \pi_{out} \rightarrow 0$.

Based on our Monte Carlo simulation results
we propose a coefficient to adjust for the proportion of $\X$ points outside $C_H(\Y_m)$,
namely,
\begin{equation}
\label{eqn:conv-hull-correct}
C_{ch}:=1-(p_{out}-\E[\widehat \pi_{out}])
\end{equation}
where $p_{out}$ is the observed and $\E[\widehat \pi_{out}] \approx 1.7932/m+1.2229/\sqrt{m}$
is the expected (under the conditions stated in Case 1)
proportion of $\X$ points outside $C_H(\Y_m)$.
For the binomial test statistic in Equation \eqref{eqn:Bnm-test-stat},
we suggest
\begin{equation}
\label{eqn:Bnm-test-stat-adj}
B^{ch}_{n,m}:=
\left\lbrace \begin{array}{ll}
(\g_n(r,M)-2J_m)\cdot C_{ch}=(\sum_{j=1}^{J_m} \g_{{}_{[j]}}(r)-2J_m)\cdot C_{ch} & \text{if  $\g_n(r,M)\cdot C_{ch}>2J_m$,}\\
       0            & \text{otherwise.}\\
\end{array} \right.
\end{equation}
For the mean domination number (per triangle) of the PCD,
we suggest
\begin{equation}
\label{eqn:Snm-test-stat-adj}
S^{ch}_{n,m} = S_{n,m} \cdot C_{ch}.
\end{equation}
This (convex hull) adjustment slightly affects the empirical size estimates in Case 1,
since $p_{out}$ and $\E[\widehat \pi_{out}]$ values are very similar.
In Cases 2-5, there is segregation when all data points are considered,
and $p_{out}$ values tend to be larger than $\E[\widehat \pi_{out}]$ values,
and in Case 6 (which is the simulation of the association case),
$p_{out}$ values tend to be smaller than $\E[\widehat \pi_{out}]$ values.
Hence in Cases 2-6,
the adjustment seems to correct the power estimates in the desired direction,
thereby increasing the power estimates.

\section{Correction for Small Samples}
\label{sec:small-sample-correction}
The distributional results in Equations \eqref{eqn:asymptotic-NYr} and \eqref{eqn:asymptotic-NYr-Jm}
might require large $n$ for the convergence to hold.
In particular, it might be necessary for
the number of $\X$ points per Delaunay triangle to be larger than 100
as a practical guide which implies very large samples from $\X$ are needed
for a large number of $\Y$ points.
Hence it might be necessary to propose a correction in the test statistics
for small $n$ also.
Based on our extensive Monte Carlo simulations (of Case 1 above)
we suggest that the test statistic $S_{n,m}$ in Equation \eqref{eqn:Snm-test-stat}
can be adjusted as $S^{adj}_{n,m}:=\frac{S_{n,m}-a_{n,m}}{b_{n,m}}$.
We provide the explicit forms of $a_{n,m}$ and $b_{n,m}$ for $m=10,20,\ldots,50$
in Table \ref{tab:small-sample-correction}.
For example for $m=10$,
$S_{n,m}$ in Equation \eqref{eqn:Snm-test-stat}
can be adjusted as $S^{adj}_{n,m}:=\frac{S_{n,m}-a_{n,m}}{b_{n,m}}$
where $a_{n,m}=-8.80/(n/J_m)-30.94/\sqrt{n/J_m}+9.09/\sqrt[3]{n/J_m}$
and $b_n=1-18.81/(n/J_m)+16.26/\sqrt{n/J_m}-4.42/\sqrt[3]{n/J_m}$.
Observe that as expected $S^{adj}_{n,m}$ converges to $S_{n,m}$ as $n \rightarrow \infty$
for each $m$ value considered provided $n/J_m \rightarrow \infty$
which is a requirement in our testing framework.

\begin{table}[ht]
\begin{center}
\begin{tabular}{|c|c|c|}
\hline
\multicolumn{3}{|c|}{$r=1.5$} \\
\hline
$m$  & $a_{n,m}$ & $b_{n,m}$\\
\hline
10 & $-8.80/(n/J_m)-30.94/\sqrt{n/J_m}+9.09/\sqrt[3]{n/J_m}$ & $1-18.81/(n/J_m)+16.26/\sqrt{n/J_m}-4.42/\sqrt[3]{n/J_m}$ \\
\hline
20 & $10.19/(n/J_m)-58.15/\sqrt{n/J_m}+20.27/\sqrt[3]{n/J_m}$ & $1-11.16/(n/J_m)+11.71/\sqrt{n/J_m}-3.24/\sqrt[3]{n/J_m}$ \\
\hline
30 & $18.72/(n/J_m)-77.36/\sqrt{n/J_m}+28.46/\sqrt[3]{n/J_m}$ & $1-6.85/(n/J_m)+7.56/\sqrt{n/J_m}-1.62/\sqrt[3]{n/J_m}$ \\
\hline
40 & $28.11/(n/J_m)-99.66/\sqrt{n/J_m}+38.73/\sqrt[3]{n/J_m}$ & $1-5.23/(n/J_m)+5.81/\sqrt{n/J_m}-0.92/\sqrt[3]{n/J_m}$ \\
\hline
50 & $33.37/(n/J_m)-115.58/\sqrt{n/J_m}+46.03/\sqrt[3]{n/J_m}$ & $1-3.93/(n/J_m)+3.88/\sqrt{n/J_m}+0.03/\sqrt[3]{n/J_m}$ \\
\hline
\hline
\multicolumn{3}{|c|}{$r=1.35$} \\
\hline
$m$  & $a_{n,m}$ & $b_{n,m}$\\
\hline
10 & $-0.13/(n/J_m)-34.35/\sqrt{n/J_m}+8.79/\sqrt[3]{n/J_m}$ & $1-16.29/(n/J_m)+13.43/\sqrt{n/J_m}-3.43/\sqrt[3]{n/J_m}$ \\
\hline
20 & $16.05/(n/J_m)-58.95/\sqrt{n/J_m}+18.01/\sqrt[3]{n/J_m}$ & $1-10.49/(n/J_m)+10.70/\sqrt{n/J_m}-3.04/\sqrt[3]{n/J_m}$ \\
\hline
30 & $24.22/(n/J_m)-77.98/\sqrt{n/J_m}+25.78/\sqrt[3]{n/J_m}$ & $1-5.59/(n/J_m)+5.52/\sqrt{n/J_m}-0.82/\sqrt[3]{n/J_m}$ \\
\hline
40 & $30.66/(n/J_m)-95.07/\sqrt{n/J_m}+32.91/\sqrt[3]{n/J_m}$ & $1-4.02/(n/J_m)+3.57/\sqrt{n/J_m}-0.06/\sqrt[3]{n/J_m}$ \\
\hline
50 & $34.49/(n/J_m)-107.87/\sqrt{n/J_m}+38.18/\sqrt[3]{n/J_m}$ & $1-3.07/(n/J_m)+2.55/\sqrt{n/J_m}+0.42/\sqrt[3]{n/J_m}$ \\
\hline
\end{tabular}

\end{center}
\caption{
\label{tab:small-sample-correction}
The finite sample adjustment for $S_{n,m}$ in Equation \eqref{eqn:Snm-test-stat}
as $S^{adj}_{n,m}:=\frac{S-a_{n,m}}{b_{n,m}}$
with $m=10,20,\ldots,50$ and $n=100,200,\ldots,1000,2000,\ldots,10000$
for $r=1.5$ (top) and $r=1.35$ (bottom).
}
\end{table}

\section{Extension of $\NPE^r$ to Higher Dimensions:}
\label{sec:NYr-higher-D}
The extension to $\R^d$ for $d > 2$ with $M=M_C$ is provided in (\cite{ceyhan:dom-num-NPE-SPL}),
the extension for general $M$ is similar:
Let $\Y = \{\y_1,\y_2,\cdots,\y_{d+1}\}$ be $d+1$ non-coplanar points.
Denote the simplex formed by these $d+1$ points as $\mS(\Y)$.
For $r \in [1,\infty]$, define the $r$-factor proximity map as follows.
Given a point $x$ in $\mS(\Y)$,
let $Q_{\y}(M,x)$ be the polytope with vertices being the $d\,(d+1)/2$ points on the edges,
the vertex $\y$ and $x$
so that the faces of $Q_{\y}(M,x)$ are formed by $d-1$ line segments each of which joining one of $\Y$
points, say $\y_i$, to $M$ and that are between $M$ and the face opposite $\y_i$.
That is, the vertex region for vertex $v$ is the polytope with vertices
given by $v$ and such points on the edges.
Let $v(x)$ be the vertex in whose region $x$ falls.
If $x$ falls on the boundary of two vertex regions, we assign $v(x)$ arbitrarily.
Let $\varphi(x)$ be the face opposite to vertex $v(x)$,
and $\eta(v(x),x)$ be the hyperplane parallel to $\varphi(x)$ which contains $x$.
Let $d(v(x),\eta(v(x),x))$ be the (perpendicular) Euclidean distance from $v(x)$ to $\eta(v(x),x)$.
For $r \in [1,\infty)$, let $\eta_r(v(x),x)$ be the hyperplane parallel to $\varphi(x)$
such that $d(v(x),\eta_r(v(x),x))=r\,d(v(x),\eta(v(x),x))$ and $d(\eta(v(x),x),\eta_r(v(x),x))< d(v(x),\eta_r(v(x),x))$.
Let $\mS_r(x)$ be the polytope similar to and with the same orientation as
$\mS(\Y)$ having $v(x)$ as a vertex and $\eta_r(v(x),x)$ as the opposite face.
Then the $r$-factor proximity region $\NY^r(x):=\mS_r(x) \cap \mS(\Y)$.
Also, let $\zeta_j(x)$ be the hyperplane such that $\zeta_j(x) \cap \mS(\Y) \not=\emptyset$ and
$r\,d(\y_j,\zeta_j(x))=d(\y_j,\eta(\y_j,x))$ for $j=1,2,\ldots,d+1$.
Then the $\G_1$-region is $\G^r_1(x)=\bigcup_{j=1}^{d+1} (\G^r_1(x)\cap R_M(\y_j))$ where
$\G^r_1(x)\cap R_M(\y_j)=\{z \in R_M(\y_j): d(\y_j,\eta(\y_j,z)) \ge d(\y_j,\zeta_j(x)\}$, for $j=1,2,\ldots,d+1$.

Let $X_{\varphi}:=\argmin_{X \in \X_n}d(X,\varphi)$ be the
closest point among $\X_n$ to face $\varphi$.
Then it is easily seen that $\G^r_1(\X_n,M)=\bigcap_{i=1}^{d+1}\G^r_1(X_{\varphi_i},M)$,
where $\varphi_i$ is the face opposite vertex $\y_i$, for $i=1,2,\ldots,d$.
So $\G^r_1(\X_n,M)\cap R_M(\y_i)=\{z \in R_M(\y_i):\; d(\y_i, \eta(\y_i,z) \ge d(\y_i,\Xi_i(X_{\varphi_i}))\}$, for $i=1,2,\ldots,d$.

Let the domination number be $\g_n(r,F,M,d):=\g_n(\X_n;F,\NPE^r,d)$ and
$X_{[i,1]}:=\argmin_{X\in \X_n \cap R_M(\y_i)}d(X,\varphi_i)$.
Then $\g_n(r,M) \le d+1$ with probability 1,
since $\X_n \cap R_M(\y_i) \subset \NPE^r\left( X_{[i,1]},M \right)$ for each of $i=1,2,\ldots,d$.

In $\mS(\Y)$, drawing the hyper-surfaces $Q_i(r,x)$ such that
$d(\y_i,\varphi_i)=r\,d(\y_i,Q_i(r,x))$ for $i\in \{1,2,\ldots,d\}$  yields another polytope,
denoted as $\msPr$, for $r<(d+1)/d$.
Let $\g_n(r,M,d):=\g(\X_n,\NPE^r,M,d)$ be the domination number of the PCD based on
the extension of $\NPE^r(\cdot,M)$ to $\R^d$.
Then we conjecture the following:
\begin{conjecture}
Suppose $\X_n$ is  set of iid random variables from the uniform
distribution on a simplex in $\R^d$. Then as $n \rightarrow \infty$,
the domination number $\g_n(r,M,d)$ in the simplex satisfies
{\small
\begin{equation}
\label{eqn:asymptotic-NYr-d}
\g_n(r,M,d) \stackrel{\mathcal L}{\longrightarrow}
\left\lbrace \begin{array}{ll}
       d+\BER(1-p_{r,d})& \text{for $r \in [1,(d+1)/d)$ and $M \in \{t_1(r),t_2(r),\ldots,t_{d+1}(r)\}$,}\\
       \le (d-1)        & \text{for $r>(d+1)/d$ and $M \in \mS(\Y)^o$,}\\
       d+1              & \text{for $r \in [1,(d+1)/d)$ and $M \in \msPr\setminus \{t_1(r),t_2(r),\ldots,t_{d+1}(r)\}$,}\\
\end{array} \right.
\end{equation}
}
\end{conjecture}
where $p_{r,d}$ can be calculated by intensive numerical integration as in the calculation
of Equation \eqref{eqn:p_r-form}
and for $r=(d+1)/d$ and $M=M_C$, $p_{r,d}$ will be different from
the continuous extension of Equation \eqref{eqn:asymptotic-NYr-d}.

\section{Discussion and Conclusions}
\label{sec:disc-conc}
In this article, we consider the asymptotic distribution of the
domination number of proportional-edge proximity catch digraphs (PCDs),
for testing bivariate spatial point patterns of segregation and association.
To our knowledge the PCD-based methods are the only graph theoretic
methods for testing spatial patterns in literature
(\cite{ceyhan:dom-num-NPE-SPL}, \cite{ceyhan:arc-density-PE}, and \cite{ceyhan:arc-density-CS}).
The new PCDs when compared to class cover catch digraphs (CCCDs), have some advantages.
In particular, the asymptotic distribution of the domination number
$\g_n(r,M)$ of the proportional-edge PCDs, unlike that of CCCDs, is mathematically tractable
(although computable by numerical integration).
A minimum dominating set can be found in
polynomial time for proportional-edge PCDs in $\R^d$ for all $d\ge 1$,
but finding a minimum
dominating set is an NP-hard problem for CCCDs (except for $\R$).
These nice properties of proportional-edge PCDs are due to the geometry invariance of
distribution of $\g_n(r,M)$ for uniform data in triangles.

On the other hand,
CCCDs are easily extendable to higher dimensions
and are defined for all $\X_n \subset \R^d$,
while proportional-edge PCDs are only defined for
$\X_n \subset \C_H(\Y_m)$.
Furthermore, the CCCDs based on balls use proximity regions
that are defined by the obvious metric, while the PCDs
in general do not suggest a metric.
In particular, our proportional-edge PCDs are based on some sort of dissimilarity measure,
but not a metric.

The finite sample distribution of $\g_n(r,M)$, although
computationally tedious, can be found by numerical methods, while
that of CCCDs can only be empirically estimated by Monte Carlo simulations.
Moreover, we had to introduce many auxiliary tools to compute the distribution
of $\g_n(r,M)$ in $\R^2$.
Same tools will work in higher dimensions,
perhaps with more complicated geometry.
The proportional-edge PCDs lend themselves for such a purpose,
because of the geometry invariance property for uniform data on Delaunay triangles.
Let the two samples of sizes $n$ and $m$ be from classes $\X$ and $\Y$, respectively,
with $\X$ points being used as the vertices of the PCDs and $\Y$ points
being used in the construction of Delaunay triangulation.
For the domination number approach to be appropriate,
$n$ should be much larger compared to $m$.
This implies that $n$ tends to infinity while $m$ is assumed to be fixed.
That is, the imbalance in the relative abundance of the two classes
should be large for this method.
Such an imbalance usually confounds the
results of other spatial interaction tests.
Furthermore,
we can also use the normal approximation to binomial distribution
for the domination number, provided $n$ is much larger than $m$,
but both sizes tending to infinity.
Therefore, as long as $n \gg m \rightarrow \infty$,
we can remove the conditioning on $m$.

The null hypothesis is assumed to be CSR of $\X$ points,
i.e., the uniformness of $\X$ points in the convex hull of $\Y$ points.
Although we have two classes here, the null pattern is not the CSR independence,
since for finite $m$,
we condition on $m$ and the locations of the $\Y$ points
are irrelevant as long as they are not co-circular.
That is, the $\Y$ points can result from any pattern
that results in a unique Delaunay triangulation.
When $m \rightarrow \infty$, conditioning on $m$ does not persist.

There are many types of parametrizations for the alternatives.
The particular parametrization of the alternatives in Equation \eqref{eqn:epsilon-alternatives}
is chosen so that the distribution of the domination number under the alternatives
would be geometry invariant
(i.e., independent of the geometry of the support triangles).
The more natural alternatives (i.e.,
the alternatives that are more likely to be found in practice)
can be similar to or might be approximated by our parametrization.
Because in any segregation alternative,
the $\X$ points will tend to be further away from $\Y$ points
and in any association alternative $\X$ points will tend to cluster around the $\Y$ points.
And such patterns can be detected by the test
statistics based on the domination number,
since under segregation (whether it is parametrized as in Section \ref{sec:alternatives}
or not) we expect them to be smaller,
and under association (regardless of the parametrization)
they tend to be larger.

By construction our method uses only the $\X$ points in $C_H(\Y_m)$
(the convex hull of $\Y$ points)
which might cause substantial data (hence information) loss.
To mitigate this, we propose a correction for the proportion of $\X$ points
outside $C_H(\Y_m)$,
because the pattern inside $C_H(\Y_m)$ might not be the same as the pattern
outside $C_H(\Y_m)$.
We suggest analysis with our domination number approach in two steps:
(i) analysis restricted to $C_H(\Y_m)$,
which provides inference only for $\X$ points in $C_H(\Y_m)$,
(ii) overall analysis with convex hull correction
(i.e., for all $\X$ points with respect to $\Y_m$).
When the number of Delaunay triangles based on $\Y$ points, denoted $J_m$,
is less than 30,
we recommend the use of binomial distribution as $n \rightarrow \infty$ (i.e., for large $n$);
when $J_m$ is larger than 30,
we recommend the use of normal approximation as $n \rightarrow \infty$.
For small samples,
one might use Monte Carlo simulation or randomization with our approach
or apply a finite sample correction as in Section \ref{sec:small-sample-correction}.
In the case of small samples with some $\X$ points existing outside $C_H(\Y_m)$,
convex hull correction can be implemented first,
and then the small sample correction.
Furthermore,
when testing against segregation we recommend the parameter $r \approx 1.3$,
while for testing against association we recommend the parameter $r \approx 1.35$
as they exhibit the best performance in terms of size and power.
The proportional-edge PCDs have applications in classification.
This can be performed building discriminant regions
in a manner analogous to the procedure proposed in \cite{priebe:2003b}.

\section*{Acknowledgments}
Supported by DARPA as administered by the Air Force Office of
Scientific Research under contract DOD F49620-99-1-0213 and by
ONR Grant N00014-95-1-0777 and by TUBITAK Kariyer
Project Grant 107T647.


\end{document}